\documentclass[11pt]{article}
\usepackage[a4paper,margin=1in]{geometry}
\usepackage[T1]{fontenc}
\usepackage{lmodern,amsmath,amssymb,amsthm,mathtools,bm,mathrsfs}
\usepackage{booktabs,tabularx,enumitem,xcolor,xurl,float,hyperref}
\hypersetup{colorlinks=true,linkcolor=blue!45!black,citecolor=blue!45!black,urlcolor=blue!45!black,pdftitle={Fluctuations and Deviations about Hydrodynamic Limits of Space-Time Hawkes Processes}}
\numberwithin{equation}{section}
\setlist{leftmargin=2em,itemsep=.25em}
\newtheorem{theorem}{Theorem}[section]
\newtheorem{proposition}[theorem]{Proposition}
\newtheorem{lemma}[theorem]{Lemma}
\newtheorem{corollary}[theorem]{Corollary}
\theoremstyle{definition}
\newtheorem{assumption}[theorem]{Assumption}
\newtheorem{definition}[theorem]{Definition}
\theoremstyle{remark}
\newtheorem{remark}[theorem]{Remark}
\newcommand{\R}{\mathbb R}
\newcommand{\E}{\mathbb E}
\newcommand{\Pp}{\mathbb P}
\newcommand{\dd}{\mathrm d}
\newcommand{\eps}{\varepsilon}
\newcommand{\NN}{\mathcal N}
\newcommand{\B}{\mathcal B}
\newcommand{\LL}{\mathcal L}

\newcommand{\cA}{\mathcal A}
\newcommand{\cK}{\mathcal K}
\newcommand{\cJ}{\mathcal J}
\newcommand{\cP}{\mathcal P}
\newcommand{\cQ}{\mathcal Q}
\newcommand{\cR}{\mathcal R}
\newcommand{\cI}{\mathcal I}
\newcommand{\cTG}{\mathcal T^{\mathrm G}}
\newcommand{\cS}{\mathcal S}
\newcommand{\cV}{\mathcal V}
\newcommand{\Id}{\operatorname{Id}}
\newcommand{\ip}[2]{\langle #1,#2\rangle}
\newcommand{\norm}[1]{\lVert #1\rVert}
\newcommand{\cD}{\mathcal D}
\newcommand{\cM}{\mathcal M}

\newcommand{\one}{\mathbf 1}
\newcommand{\supp}{\operatorname{supp}}
\DeclareMathOperator{\Cov}{Cov}

\title{Fluctuations, Moderate and Large Deviations about Hydrodynamic Limits of Space-Time Hawkes Processes}

\author{
Fuqing Gao\thanks{School of Mathematics and Statistics, Wuhan University, Wuhan, China; \texttt{fqgao@whu.edu.cn}}
\and
Yingli Wang\thanks{School of Mathematical Sciences, Fudan University, Shanghai, People's Republic of China; \texttt{yingliwang@fudan.edu.cn}}
\and
Lingjiong Zhu\thanks{Department of Mathematics, Florida State University, Tallahassee, Florida, United States of America; \texttt{zhu@math.fsu.edu}}}
\date{\today}

\begin{document}
\maketitle
\begin{abstract}
We study Gaussian fluctuations, moderate and large deviations around the high-density hydrodynamic limit of critical and nearly critical space-time Hawkes processes initialized by a prescribed Poisson field of past events, with no subsequent immigration. Under a finite mean reproduction delay and centered displacements with finite spatial covariance, a positive cluster-transform estimate yields a distribution-valued central limit theorem (CLT) and a moderate deviation principle (MDP) throughout the moderate range. For finite nonzero mass of the initial profile, we obtain a full large deviation principle (LDP) for the space--time event measure in the vague topology; its good rate function is a quadratic branching-control energy with a prescribed initial trace and a finite-mass domain. The lower bound uses an explicit approximation by smooth reproduction controls, and the upper bound uses a microscopic exponential martingale and a global generator approximation. We distinguish exact-mean centering from PDE centering and compute deterministic bias corrections under second delay and fourth spatial moments. For a factorized heavy-tail power-law delay, we establish the CLT, MDP, and full measure-valued LDP for the space--time measure, with finite nonzero mass of the initial profile for the LDP. Its fractional control energy retains the source generated by the prescribed past. A source-preserving regularization proves approximation of every finite-action measure. Compact-kernel bounds for initial profiles of infinite mass and the degenerate case of a zero initial profile are stated separately.
\end{abstract}
\noindent\textbf{Keywords:} Hawkes process; hydrodynamic limit; central limit theorem; moderate deviations; large deviations; prescribed past; renewal resolvent; fractional diffusion.

\noindent\textbf{2020 MSC:} 60G55; 60F05; 60F10; 60F17; 60J80.

\section{Introduction}\label{sec:intro}
Hawkes processes describe event sequences in which past events change
the rate of future occurrences. In a self-exciting model, an event raises
the intensity for a period of time and may trigger further events,
producing clusters of activity. Since the work of Hawkes~\cite{Hawkes},
these processes have been used to model earthquake sequences
\cite{Ogata1988,Ogata1998}, financial transactions \cite{BMM}, neuronal
activity \cite{Pernice}, and crime occurrences \cite{Mohler}.
They also appear in queueing models \cite{GZstationary}.
The monograph \cite{LLT} gives an account of Hawkes models and their
branching representation, and Reinhart~\cite{Reinhart} surveys
self-exciting point processes with spatial and temporal coordinates.
Bernabeu, Zhuang, and Mateu~\cite{BernabeuZhuangMateu2025} review
spatio-temporal Hawkes models, with emphasis on simulation methods
and statistical estimation.

For a space--time process, reproduction has two effects. The delay
between a parent and a child determines the propagation of activity in
time, while the displacement of the child determines its spread in
space. At the same time, random offspring numbers produce fluctuations
in the size of each cluster. Near the critical branching ratio, large
clusters can contribute appreciably to macroscopic observations.
The resulting limit depends on the reproduction law, the observation
scales, and the density of events prescribed before the observation
begins. Understanding the mean evolution therefore leads to further
questions about typical fluctuations and the probabilities of larger
departures from that evolution.

This paper addresses these questions for the space--time
Hawkes model studied by Gao, Wang, and Zhu~\cite{GWZ2026}.
The process is initialized by a prescribed Poisson field of past events and is then
followed through its positive-time descendants. In the high-density
regime, its normalized event measure converges to a deterministic
diffusion density. We establish Gaussian fluctuations about this
limit, moderate deviations between the Gaussian and macroscopic
scales, and large deviations of the event measure itself. The
prescribed Poisson past, physical time scale, and spatial scale are the same as
in the hydrodynamic model. This permits a direct comparison between
its deterministic limit, its fluctuation covariance, and its
large-deviation rate function.

\subsection{Hawkes process and space-time Hawkes process}

A (temporal) Hawkes process \cite{Hawkes} is a random counting measure $N(\dd t)$ with the predictable intensity
\begin{equation}\label{eq:intro-temporal-Hawkes}
  \lambda_t:=\mu+\int_{(-\infty,t)}g(t-s)N(\dd s),
  \qquad \mu>0,
\end{equation}
where $g:\mathbb{R}_{+}\rightarrow\mathbb{R}_{\geq0}$ is the \textit{kernel function} (or exciting function)
\cite{Hawkes} and $\mu>0$ is the \textit{baseline intensity}.  

A space--time Hawkes process replaces event times by event locations
$(t,x)$ in $\mathbb R_+\times\mathbb R^d$.  It is a random counting measure
$N(\dd t,\dd x)$ with predictable intensity density $\lambda(t,x)$: for
each bounded Borel set $A\subset\mathbb R^d$,
\begin{equation}\label{eq:intro-spacetime-compensator}
  N((0,t]\times A)
  -\int_0^t\int_A\lambda(s,x)\dd x\dd s
\end{equation}
is a local martingale. In the translation-invariant linear model \cite[Section~2.2, equation~(3)]{BernabeuZhuangMateu2025},
\begin{equation}\label{eq:intro-Hawkes}
  \lambda(t,x)
  :=\mu(t,x)+\int_{(-\infty,t)\times\mathbb R^d}
       g(t-s,x-y)N(\dd s,\dd y).
\end{equation}
Here $\mu:\mathbb{R}\times\mathbb{R}^{d}\rightarrow\mathbb{R}_{\geq 0}$ is the baseline intensity and $g:\mathbb{R}\times\mathbb{R}^{d}\rightarrow\mathbb{R}_{\geq 0}$ is the kernel function (or exciting function) with $g(r,\cdot)=0$ for $r\leq0$.

\subsection{The immigration-birth representation and related limits}

Let $h(s,z)\dd s\dd z$ be a probability law on
$(0,\infty)\times\R^d$. In the immigration--birth representation of
Hawkes and Oakes~\cite{HawkesOakes1974}, a parent at $(r,y)$ produces a Poisson
random measure of children at $(r+s,y+z)$ with intensity
$a h(s,z)\dd s\dd z$. With a general baseline $\mu(t,x)$, the
predictable intensity is 
\begin{equation}\label{eq:intro-spacetime-model}
 \lambda(t,x)=\mu(t,x)+a\int_{(-\infty,t)}\int_{\R^d}
                  h(t-r,x-y)N(\dd r,\dd y).
\end{equation}
The normalization of $h$ makes $a$ the mean number of direct children.
For $a<1$, the expected total cluster size is
$\sum_{j=0}^\infty a^j=(1-a)^{-1}$. At $a=1$, the Poisson
Galton--Watson genealogy becomes extinct almost surely, while its
expected total size is infinite. We study the descendants of a
prescribed past on a finite observation horizon, estimating the
contribution of old parents through renewal bounds.

The pair $(S,\Xi)$ with density $h$ describes the reproduction delay
and displacement along one genealogical edge. Its coordinates may
be dependent. In the finite-mean setting we assume
\[
 m:=\E[S]\in(0,\infty),\qquad \E[\Xi]=0,\qquad
 \Sigma:=\E\left[\Xi\Xi^\top\right]\text{ is finite and positive definite}.
\]
Successive generations combine an age-dependent branching process
with a branching random walk. The mean delay converts the number of
generations into elapsed physical time, and the covariance of the
centered displacements determines the spatial diffusion.

For the model considered here, $\mu\equiv 0$ after time zero and the
negative-time counting measure is prescribed separately. Its atoms
produce children only at positive times; descendants then reproduce
according to the same law. This construction retains the influence
of earlier events through their remaining reproduction delays.
The prescribed past is random, and its contribution to fluctuations
must be estimated together with the randomness generated by the
positive-time clusters. The joint delay--displacement law also makes
the renewal estimates sensitive to time and space simultaneously.

Age-dependent branching provides a second route to understanding
these processes. Lambert, Simatos, and Zwart~\cite{LSZ} use excursion
theory to relate Crump--Mode--Jagers branching processes and
processor-sharing queues. Xu~\cite{XuDiffusion} develops diffusion
approximations for self-excited systems and general branching models,
while Dram\'e and Pardoux~\cite{DP} prove many-ancestor laws of large
numbers and central limit theorems for Crump--Mode--Jagers populations.
In spatial population models, M\'el\'eard and Tran~\cite{MT} obtain
superprocess limits by separating fast age and slow trait scales.
Fleischmann and Vatutin~\cite{FV} study multiscale clustering in a
spatial branching system with finite-mean lifetimes and increasing
initial populations. For the Hawkes model, each parent can give birth
at several ages, and the measure studied here counts births over the
whole space--time observation region. Its renewal clock and
prescribed Poisson past must therefore be retained in the limiting
equations.

Spatial Hawkes asymptotics have also been studied on networks and
in spatial mean-field systems. Delattre, Fournier, and
Hoffmann~\cite{DFH} construct Hawkes processes on possibly infinite
networks and analyze mean-field and long-time limits. Spatial
neural-field equations and Gaussian fluctuations are obtained in
\cite{CDLO,CO}. MacLaurin~\cite{MacLaurin} considers a hydrodynamic
limit for adaptive stochastic networks, and
Belmabrouk~\cite{Belmabrouk} studies long-time laws of large numbers
on an infinite graph with power-law interactions. The present model
has a continuum of physical locations. Rescaling these locations
together with event times produces a spatial diffusion operator,
and changing the density of the prescribed past determines the
size of its branching fluctuations.

\subsection{Fluctuations, moderate and large deviations in the literature}

The fluctuations, moderate and large deviations
for temporal Hawkes processes have been
well studied in the literature. 
For the linear Hawkes process, 
\cite{HawkesOakes1974} establishes a central
limit theorem (CLT) for linear Hawkes processes using the immigration-birth
representation.
\cite{BacryEtAl2013} extends 
the central limit theorem to the functional level and the multivariate setting.
\cite{BordenaveTorrisi2007} proves
a large deviation principle with an explicit rate function,
and \cite{Karim2025} extends the result to the multivariate setting with random marks. 
\cite{Zhu2013MDP} obtains the moderate deviations and \cite{Yao2018} extends it to the multivariate setting.
Precise large-time deviation asymptotics are obtained in \cite{GZprecise} using
the mod-$\phi$ convergence theory.

For nonlinear Hawkes processes, 
\cite{Zhu2013CLT} proves a functional central limit theorem.
Moreover, \cite{ZhuProcess} proves a
process-level (also known as Level-3) LDP \cite{DonskerVaradhan1983} for the empirical process of time shifts, with
speed $t$ and a specific-entropy rate function.  By the contraction principle \cite{DemboZeitouni2010}, this leads to
a Level-1 LDP for $N_t/t$.  For Markovian nonlinear Hawkes processes with
exponential kernels or finite sums of exponential kernels,
\cite{Zhu2015MarkovLDP} develops an alternative spectral and variational
characterization of the Level-1 rate function.
Most recently, Wang and Zhu~\cite{WZ2026} obtain moderate deviations using a
configuration-space corrector and martingale method. 
Other than the large-time scaling, 
Gao and Zhu~\cite{GZasymp} study a large-rate, small-excitation scaling, and obtain fluctuations, 
moderate and large deviation results.
\cite{GZlimits} obtain a functional central limit theorem for Markovian Hawkes processes when the initial intensity is large, 
and the large deviations are studied
in \cite{GZldp} in the same asymptotic regime.
For the mean-field limit of Hawkes processes, 
fluctuations, 
moderate and large deviations are studied
in \cite{GGZ,GZmeanLDP}. 

\subsection{Hydrodynamic scaling and the high-density regime}
The passage from microscopic dynamics to a macroscopic evolution
equation is a central problem in the hydrodynamic theory of
interacting particle systems; see \cite{Spohn,KL}. Fluctuations and
large deviations describe the random departures from that evolution,
as in the exclusion-process theory of Kipnis, Olla, and
Varadhan~\cite{KOV}. For a Hawkes event measure, diffusion results
from the accumulation of offspring displacements, while the initial
condition results from reproduction across time zero. The
past intensity provides an additional parameter that controls
whether random branching remains visible on the hydrodynamic scale.

Let $\eps\downarrow0$. In the finite-mean model, we assume that the observation
scales are $(t,x)\mapsto(\eps^2t,\eps x)$, and the prescribed past
has intensity
\[
 \eps^{-(\alpha+1-d)}\lambda_0(\eps y)\dd r\dd y,
 \qquad r<0,
\]
where $\lambda_0\geq0$ is bounded and uniformly continuous.
The normalized event measure assigns mass $\eps^{\alpha+3}$ to
each observed event. The near-critical balance is
\[
 \delta_\eps:=\eps^2,\qquad
 \kappa_\eps:=\frac{1-a_\eps}{\delta_\eps}\longrightarrow\kappa\geq0.
\]
On a time interval of order $\eps^{-2}$, an ancestral line has
order $\eps^{-2}$ generations. Centered displacements accumulated
over those generations have order $\eps^{-1}$. For a smooth
compactly supported test function, the corresponding offspring
operator satisfies
\begin{equation}\label{eq:intro-generator-limit}
 \begin{aligned}
 \cK_\eps f(t,x)&=\E\left[f\left(t+\delta_\eps S,x+\eps\Xi\right)\right],\\
 \frac{a_\eps \cK_\eps f-f}{\delta_\eps}
 &\longrightarrow m\partial_t f+\frac12\Sigma:D_x^2f-\kappa f.
 \end{aligned}
\end{equation}
The convergence in \eqref{eq:intro-generator-limit} is uniform on
the observation region by \eqref{eq:global-generator} and
the near-critical scaling in \eqref{eq:scales}.
The hydrodynamic density consequently solves
\begin{equation}\label{eq:intro-hydrodynamic-equation}
 m\partial_t\rho=\frac12\Sigma:D_x^2\rho-\kappa\rho,
 \qquad \rho(0,x)=\lambda_0(x).
\end{equation}
The spatial diffusion matrix is $\Sigma/(2m)$ and the killing rate
is $\kappa/m$. Exact microscopic criticality $a_\eps=1$ and
near-critical sequences with $1-a_\eps=o(\eps^2)$ both give
$\kappa=0$ in this equation.

The density parameter $\alpha$ determines a separate transition.
The results of \cite{GWZ2026} identify three regimes: for
$\alpha>1$, branching noise averages out and the event measure
converges to $\rho(t,x)\dd t\dd x$; at $\alpha=1$, the limit is
a superprocess occupation measure; for $0<\alpha<1$, the event
measure converges locally to zero in probability, although its mean
has a nonzero limit when the initial density is nontrivial.
The last regime reflects the contribution of rare large clusters
to the first moment. In the high-density regime considered in this
paper, the effective population parameter is
\[
 n_\eps:=\eps^{-(\alpha-1)}\longrightarrow\infty,
 \qquad \alpha>1.
\]
The normalized equation-scale martingale bracket has order
$n_\eps^{-1}$. This gives the Gaussian scale $n_\eps^{-1/2}$
and suggests $n_\eps$ as the speed of macroscopic large deviations.
The proofs below establish both conclusions through the cluster
transform and the microscopic exponential martingale.

When the offspring mean approaches one, the distribution of a cluster
changes with $\eps$. The positive exponential moments needed for
an MDP must then be controlled uniformly over this changing family.
Our analytic cluster estimate provides that control under the finite
mean delay and finite spatial covariance assumptions above.
In the linear
model considered here, independent Poisson clusters permit exact
superposition at every positive real intensity of the prescribed past. This
property is used both to identify the cumulants and to prove the
spatial large-deviation lower bound.

\subsection{Our contributions}

Fix $0<T<\infty$, let $\NN_\eps$ be the normalized event measure \eqref{eq:event},
and write $\mu_\eps:=\E[\NN_\eps]$ for its exact mean and
$\mu:=\rho(t,x)\dd t\dd x$ for its hydrodynamic limit.
Theorems~\ref{thm:clt} and \ref{thm:mdp} concern
\begin{equation}\label{eq:intro-fluctuation-scales}
 Z_\eps:=\sqrt{n_\eps}(\NN_\eps-\mu_\eps),\qquad
 Y_\eps:=b_\eps^{-1}Z_\eps,\qquad
 b_\eps\to\infty,\quad b_\eps/\sqrt{n_\eps}\to0.
\end{equation}
The first family converges to a centered Gaussian random
distribution. The second satisfies an MDP at speed $b_\eps^2$
throughout the range in \eqref{eq:intro-fluctuation-scales}.
For $\B:=\Sigma:D_x^2/2$, the backward heat resolvent is
\[
 \cA f(t,x):=\frac1m\int_t^T
       e^{-(s-t)\kappa/m}\left(e^{(s-t)\B/m}f(s,\cdot)\right)(x)\dd s.
\]
The covariance is
\begin{equation}\label{eq:intro-gaussian-covariance}
 C(f,g):=\int_0^T\int_{\R^d}\rho(t,x)\cA f(t,x)\cA g(t,x)\dd x\dd t,
 \qquad f,g\in\cD((0,T)\times\R^d).
\end{equation}
The same resolvent determines the quadratic MDP rate function through
$L^2$ controls in the linear equation
$m\partial_tz=\B z-\kappa z+\sqrt\rho\,\zeta$, with zero initial
value. The Gaussian covariance is thus the quadratic variation of
branching noise propagated by the hydrodynamic equation.
The separate covariance term generated by the random Poisson past
has an additional factor $\delta_\eps$ and tends to zero.
Localized Fourier estimates give tightness and exponential tightness
in the weak distribution topology.

For finite nonzero initial mass
$L:=\int_{\R^d}\lambda_0(x)\dd x$, Theorem~\ref{thm:full-spatial-finite}
proves a full LDP at speed $n_\eps$ for the positive event measure
in the vague topology. Every finite-rate measure has a density
$q\in C([0,T];L^1(\R^d))$, initial trace $q_0=\lambda_0$, and a
control $v$ satisfying
\begin{equation}\label{eq:intro-spatial-action}
 m\partial_tq=\B q-\kappa q+vq,\qquad q_0=\lambda_0,
\end{equation}
such that the rate function for the LDP is
\begin{equation} 
\mathcal J_{\lambda_0}(q\dd t\dd x)
:=\frac12\int_0^T\int_{\R^d}(v(t,x))^2q(t,x)\dd x\dd t.
\end{equation}
The control changes the local reproduction rate. Its quadratic cost
is weighted by the density of events on which it acts. The
initial trace and finite-mass domain follow from the variational
weak equation and a mass--energy estimate. The upper bound uses
the microscopic exponential martingale, uniform generator
convergence, and total-mass truncation.

The quadratic control energy associated with \eqref{eq:intro-spatial-action} is related to
superprocess large deviations. Fleischmann, G\"artner, and
Kaj~\cite{FGK} and Schied~\cite{Schied} establish Schilder-type
and sample-path results for super-Brownian motion.
Schied~\cite{SchiedMDP} also studies moderate deviations and a
functional law of the iterated logarithm. Xiang~\cite{Xiang}
identifies a rate function for finite nonzero initial measures and treats
infinite initial measures in weighted path spaces~\cite{XiangInfinite}.
Here the Hawkes lower bound is proved directly. We average the forward
heat-and-time-translation semigroup and add a nonnegative correction
that preserves $q_0=\lambda_0$. A Gaussian derivative estimate controls
the energy of spatial truncation. Lemma~\ref{lem:finite-recovery}
then approximates every finite-action density in $L^1$ with convergence
of costs. A microscopic intensity change realizes each smooth control.
Finite total mass controls the global generator error, so the proof
requires only the original finite-moment assumptions.

Two extensions clarify the dependence on centering and reproduction
memory. Section~\ref{sec:bias} compares the exact mean $\mu_\eps$
with the hydrodynamic measure $\mu$. Their difference is
deterministic, but it must be small on the fluctuation scale before
one center can replace the other. Under $\E[S^2]<\infty$ and
$\E[|\Xi|^4]<\infty$, we estimate this difference and compute its
first two weak coefficients. Conditional spatial symmetry removes
the first coefficient. The resulting expansions identify Gaussian
limits with a deterministic shift at the corresponding thresholds for
the density parameter and give explicit sufficient conditions for
hydrodynamic centering in the MDP. These calculations require no
derivatives of $\lambda_0$.

For factorized reproduction laws with
$\Pp(S>s)\sim c_Hs^{-\eta}$, $0<\eta<1$, the appropriate spatial
scale becomes $\eps^\eta$ on the same time scale $\eps^2$.
With the matched scaling of the Poisson past in \cite{GWZ2026}, the effective
population is $\eps^{-\eta(\alpha-1)}$ and the near-critical
scale is $\eps^{2\eta}$. Section~\ref{sec:fractional} gives
Gaussian, moderate, and large deviations for the space--time measure. Under finite nonzero mass of the initial profile,
Theorem~\ref{thm:fractional-spatial-ldp} proves the full space--time
LDP at speed $n_{\eps,\eta}:=\eps^{-\eta(\alpha-1)}$; its rate function is
\begin{equation}\label{eq:intro-fractional-action}
 \mathcal J_{\eta,\lambda_0}(q\dd t\dd x)
 =\frac12\int_0^T\int_{\R^d}\frac{(j(t,x))^2}{q(t,x)}\dd x\dd t,
 \qquad q=\rho^\eta+\cA_\eta^*j.
\end{equation}
The residual includes the source generated by the prescribed past
$c_Ht^{-\eta}\lambda_0(x)\dd t\dd x$, which is integrable at zero.
The lower bound follows by smoothing with the stable-time--Brownian-space
semigroup and restoring this source exactly. A Gaussian derivative
estimate controls spatial truncation, and approximation by bounded
reproduction controls then covers every finite-action measure.
The upper bound uses the fractional generator and the microscopic
exponential martingale. These arguments apply throughout $0<\eta<1$.

Corollary~\ref{cor:fractional-action-tangent} identifies the
second-order expansion of the full fractional action with the Gaussian
control energy. The mean follows a Caputo fractional diffusion equation,
with coefficient $c_\eta=c_H\Gamma(1-\eta)$. The inverse-clock
interpretation is discussed in \cite{MS}. Heavy-tailed excitation marks, as
studied by Blanchet et al.~\cite{BLWZ}, give a different branching
mechanism from the heavy-tailed delay considered here.

The full measure-valued LDPs assume finite nonzero mass of the initial profile
in both delay regimes. For possibly infinite initial mass, the
finite-mean compact-kernel subclass has the spatial exponential
bounds and conditional approximation result stated in
Section~\ref{subsec:infinite-bounds}; the case of a zero initial profile is treated in
Proposition~\ref{prop:zero-obstruction}.

\subsection{Organization}
Section~\ref{sec:model} specifies the prescribed past, scaling,
operators, state spaces, and additional assumptions.
Section~\ref{sec:main-results} states the finite-mean results.
Sections~\ref{sec:analytic} and \ref{sec:clt-mdp} prove the renewal
and cluster estimates, the CLT, and the MDP.
Section~\ref{sec:spatialLDP} proves the finite-mass LDP
and the additional compact-kernel bounds.
Sections~\ref{sec:bias} and \ref{sec:fractional} treat deterministic
centering corrections and heavy-tailed delays, including the full
fractional space--time LDP.
Section~\ref{sec:conclusion} concludes.

\section{Model, assumptions, and scaling}\label{sec:model}

\subsection{The offspring law and critical scaling}
Fix $d\geq1$ and $0<T<\infty$. The parameter $\eps$ ranges over $(0,1]$, and all limits indexed by $\eps$ are taken as $\eps\downarrow0$.

\begin{assumption}[Finite-mean high-density model]\label{ass:finite-mean} 
Assume that $h$ is a nonnegative Boreal measurable function and that $\int_{\mathbb{R}^{d}}\int_{0}^{\infty}h(s,z)\dd s\dd z=1$.
Let $(S,\Xi)$ follow a probability law $h(s,z)\dd s\dd z$ on $(0,\infty)\times\R^d$. Moreover, assume that
\begin{equation}\label{eq:basic}
 m:=\E[S]\in(0,\infty),\qquad \E[\Xi]=0,\qquad
 \Sigma:=\E\left[\Xi\Xi^\top\right]\text{ is finite and positive definite}.
\end{equation}
The delay $S$ and displacement $\Xi$ may be dependent. Let $\lambda_0\geq0$ be bounded and uniformly continuous. Write
\begin{equation}\label{eq:scales}
 \delta_\eps:=\eps^2,\qquad
 n_\eps:=\eps^{-(\alpha-1)},\quad \alpha>1,\qquad
 a_\eps\in(0,1],\qquad
 \kappa_\eps:=\frac{1-a_\eps}{\delta_\eps}\longrightarrow\kappa\geq0.
\end{equation}
\end{assumption}
In the notation of \cite{GWZ2026}, $m=m_1$. We retain $m$ for the mean delay throughout this paper. The reproduction delay density and its tail are
\begin{equation}\label{eq:delay-marginal}
 H(s):=\int_{\R^d}h(s,z)\dd z,\qquad
 \overline H(s):=\int_s^\infty H(r)\dd r,\qquad
 \widehat H(z):=\int_0^\infty e^{-zs}H(s)\dd s\quad(z>0).
\end{equation}
\emph{Exact microscopic criticality} means $a_\eps\equiv1$. The condition $\kappa=0$ is weaker: it also allows $1-a_\eps=o(\eps^2)$. Both have the same limiting critical operator. A strictly positive $\kappa$ gives limiting killing rate $\kappa/m$. Both cases are covered by the same statement and proof. The parameter $\alpha>1$ controls the past intensity, while $a_\eps$ is the branching ratio.

\subsection{Prescribed past and the event observable}

We assume that the prescribed past is a Poisson measure $\Pi^\eps$ of intensity
\begin{equation}\label{eq:past}
 \eps^{-(\alpha+1-d)}\lambda_0(\eps y)\dd r\dd y,
 \qquad r<0.
\end{equation}
At positive times, 
i.e. on $(0,\infty)\times\R^d$, the space-time Hawkes process $N^\eps$ is a locally finite point process,
with no two events at the same time, with the predictable intensity
\begin{align}\label{eq:intensity}
 \lambda^\eps(t,x)
 &=a_\eps\int_{r<0}\int_{\R^d}h(t-r,x-y)\Pi^\eps(\dd r,\dd y)\nonumber\\
 &\qquad\qquad+a_\eps\int_{0<r<t}\int_{\R^d}h(t-r,x-y)N^\eps(\dd r,\dd y).
\end{align}
The atoms of $\Pi^\eps$ are prescribed parents. Each parent produces children at positive times, and these children generate positive-time descendants by independent Poisson reproduction. This gives the cluster construction of \cite{HawkesOakes1974} with the prescribed initial history. The process is constructed by increasing finite collections of prescribed parents and finite offspring generations. For each individual cluster, the Poisson Galton--Watson offspring law has mean $a_\eps\leq1$ and becomes extinct almost surely, so that the full cluster is finite almost surely. After superposing the infinitely many past parents, local finiteness on bounded space--time regions follows from the mean bound below. Monotone iteration of the mean-intensity equation, with the constant microscopic past intensity as a supersolution, gives the minimal bounded-mean solution and
\begin{equation}\label{eq:meanbound}
 \E\left[N^\eps(\dd t,\dd x)\right]
 \leq \eps^{-(\alpha+1-d)}\norm{\lambda_0}_\infty\dd t\dd x.
\end{equation}
Define the \textit{event measure} 
\begin{equation}\label{eq:event}
 \ip{\NN_\eps}{f}
 :=\frac{\delta_\eps^2}{n_\eps}
   \int_0^\infty\int_{\R^d}f(\delta_\eps t,\eps x)N^\eps(\dd t,\dd x),
 \qquad f\in C_c((0,\infty)\times\R^d), 
\end{equation}
and its exact deterministic mean by
\begin{equation}
 \mu_\eps:=\E\left[\NN_\eps\right].
\end{equation}
Thus $\delta_\eps^2/n_\eps=\eps^{\alpha+3}$, exactly as in \cite{GWZ2026}. The \emph{effective population parameter} $n_\eps$ controls the intensity of independent clusters.

\subsection{Microscopic operators and finite-horizon norms}

A pairing with a space--time density $q$ means pairing with the measure $q(t,x)\dd t\dd x$ on the stated time interval and $\R^d$. The microscopic event measure is defined on $(0,\infty)\times\R^d$ and is restricted to the observation horizon when needed.

For functions extended by zero beyond $T$, let
\begin{equation}\label{eq:K}
 \cK_\eps f(t,x):=\E\left[f\left(t+\delta_\eps S,x+\eps \Xi\right)\right],
 \qquad \cA_\eps:=\delta_\eps\sum_{j=0}^\infty a_\eps^j \cK_\eps^j,
\end{equation}
and define
\begin{equation}\label{eq:discrete-generator}
 \LL_\eps:=\frac{a_\eps \cK_\eps-\Id}{\delta_\eps}.
\end{equation}
For a measurable function $f:[0,T]\times\R^d\to\mathbb C$ and
$1\leq p\leq\infty$, the mixed norm is
\begin{equation}\label{eq:mixed-norm}
 \norm{f}_{L^\infty([0,T];L^p(\R^d))}
 :=\operatorname*{ess\,sup}_{0\leq t\leq T}
                  \norm{f(t,\cdot)}_{L^p(\R^d)}.
\end{equation}
The spatial $L^p(\R^d)$ norm is taken first, followed by the
essential supremum over $t\in[0,T]$. In particular, for $p=1$ the spatial norm is
$\int_{\R^d}|f(t,x)|\dd x$, and for $p=\infty$ it is
$\operatorname*{ess\,sup}_{x\in\R^d}|f(t,x)|$.
For $q\geq0$ and $p\in\{1,\infty\}$, we also use the weighted
uniform-in-time norm
\begin{equation}\label{eq:weighted-norm}
 \norm{f}_{q,p}:=\sup_{0\leq t\leq T}
                 e^{qt}\norm{f(t,\cdot)}_{L^p(\R^d)}.
\end{equation}
For continuous $L^p(\R^d)$-valued functions of time, the supremum and
essential supremum coincide, so that
$\norm{f}_{0,p}=\norm{f}_{L^\infty([0,T];L^p(\R^d))}$.
For a linear operator $\mathcal O$, the notation
$\norm{\mathcal O}_{q,p}$ denotes the induced operator norm
\begin{equation}\label{eq:weighted-operator-norm}
 \norm{\mathcal O}_{q,p}
 :=\sup_{0<\norm{f}_{q,p}<\infty}
       \frac{\norm{\mathcal O f}_{q,p}}{\norm{f}_{q,p}}.
\end{equation}

\subsection{Hydrodynamic density and backward resolvent}
Set
\begin{equation}\label{eq:operators}
 \B:=\frac12\Sigma:D_x^2,\qquad
 \LL_\kappa:=m\partial_t+\B-\kappa,\qquad
 \cQ_t:=e^{-\kappa t/m}\cP_t,\qquad \cP_t:=\exp(t\B/m).
\end{equation}
Thus, $\cP_t$ has covariance $t\Sigma/m$. The hydrodynamic density is $\rho_t=\cQ_t\lambda_0$, so
$m\partial_t\rho=\B\rho-\kappa\rho$, $\rho_0=\lambda_0$.
Write $\mu(\dd t,\dd x):=\rho(t,x)\dd t\dd x$ for the hydrodynamic reference measure. Section~\ref{sec:bias} compares it with the exact mean $\mu_\eps$.
The change of variables in \eqref{eq:event}, the mean bound
\eqref{eq:meanbound}, and the contraction of $\cQ_t$ give
\begin{equation}\label{eq:normalized-mean-density-bound}
 \mu_\eps(\dd t,\dd x)\leq\|\lambda_0\|_\infty\dd t\dd x,
 \qquad 0\leq\rho(t,x)\leq\|\lambda_0\|_\infty.
\end{equation}
On test functions vanishing after $T$, the backward resolvent is
\begin{equation}\label{eq:resolvent}
 \cA  f(t,x)=\frac1m\int_t^T \cQ_{s-t}f(s,\cdot)(x)\dd s,
 \qquad -\LL_\kappa \cA  f=f,\qquad \cA  f(T)=0.
\end{equation}
The adjoint $\mathcal{A}^{\ast}$ acts on space-time measures by $\langle\mathcal{A}^{\ast}\gamma,f\rangle=\langle \gamma,\mathcal{A}f\rangle$. 
In particular, for an integrable density $j$, $(\mathcal{A}^{\ast}j)(t)=m^{-1}\int_{0}^{t}\cQ_{t-s}j(s)\dd s$.
Our covariance form is
\begin{equation}\label{eq:cov}
 C(f,g)=\int_0^T\int_{\R^d}\rho(t,x) \cA  f(t,x)\cA  g(t,x)\dd x\dd t.
\end{equation}
Its finiteness follows from the $L^1$ and $L^\infty$ resolvent bounds. The factor $1/m$ in \eqref{eq:resolvent} is essential. Semigroup duality will repeatedly be used in the form
\begin{equation}\label{eq:mean-duality}
 m\int_{\R^d}\lambda_0(x)\cA  f(0,x)\dd x
 =\int_0^T\int_{\R^d}\rho(t,x)f(t,x)\dd x\dd t.
\end{equation}

\subsection{State spaces and deviation conventions}\label{subsec:spaces}

Let $E_T:=(0,T)\times\R^d$ and $\Omega_T:=[0,T)\times\R^d$.
For fluctuations we use $\cD(E_T)=C_c^\infty(E_T)$ and
$\cD'(E_T)$ with $\sigma(\cD',\cD)$, the weak distribution topology.
For the LDP we use
\begin{equation}\label{eq:spatial-state-space}
 \mathsf E:=\cM_{\mathrm{loc},+}(\Omega_T)
\end{equation}
with the vague topology, that is, convergence against $C_c(\Omega_T)$.
Standard random-measure
conventions are as in \cite{Kallenberg}. Test functions in
$C_c^\infty(\Omega_T)$ may be nonzero at $t=0$ and vanish near $t=T$;
this convention keeps the initial trace in weak equations. Fluctuation
test functions in $\cD(E_T)$ vanish near both endpoints. The microscopic measure
is defined on all positive times before any of these restrictions is made.

For $\cD'(E_T)$, probability laws are considered on the Borel
$\sigma$-field of the stated weak topology. A localized finite signed
measure is a measurable $H^{-s}$-valued random element for
$s>(d+1)/2$: its Fourier coordinates are measurable and their weighted
squares are bounded by the square of its local total variation times
a summable sequence. This applies to the centered microscopic event
measures before any limit is taken. The Gaussian limit is constructed
from white noise through \eqref{eq:Gaussian}, with covariance \eqref{eq:cov} and the same local
Sobolev realization. Lemma~\ref{lem:distribution-tightness} proves
tightness in this space. On its compact sets, countably many smooth Fourier
test functions determine the topology and the law, which permits the
projective arguments below.

The Gaussian initial condition is specified by the mild representation
\eqref{eq:Gaussian}. In a finite-action representation, $q(t,x)$ denotes
the density of the occupation measure.

Let $\mathcal{X}$ be a Polish space. A family of measurable random elements $X_\eps$ in $\mathcal{X}$
satisfies a \textit{large deviation principle} (LDP) at speed $s_\eps\to\infty$ and rate
function $I:\mathcal{X}\rightarrow[0,+\infty]$ that is lower-semicontinuous, 
when, for every open set $O$ and closed set $F$ in $\mathcal{X}$,
\begin{equation}\label{eq:LDP-convention}
 -\inf_{x\in O} I(x)\leq\liminf_{\eps\downarrow0}s_\eps^{-1}\log\Pp(X_\eps\in O),
 \qquad
 \limsup_{\eps\downarrow0}s_\eps^{-1}\log\Pp(X_\eps\in F)\leq-\inf_{x\in F} I(x).
\end{equation}
A rate function is said to be \emph{good} if every finite sublevel set is compact. Exponential
tightness means that for each $A<\infty$ there is a compact set $K_A$ with
$\limsup s_\eps^{-1}\log\Pp(X_\eps\notin K_A)\leq-A$.
An MDP is the LDP for the centered, moderately rescaled variables defined
below. We use these conventions throughout; see \cite{DemboZeitouni2010}.

\subsection{Additional assumptions and notation}
\begin{assumption}[Finite nonzero mass of the initial profile]\label{ass:finite-mass}
In addition to Assumption~\ref{ass:finite-mean},
\begin{equation}\label{eq:finiteinitial}
 L:=\int_{\R^d}\lambda_0(x)\dd x\in(0,\infty).
\end{equation}

\end{assumption}
Assumption~\ref{ass:finite-mass} is used in
Theorem~\ref{thm:full-spatial-finite}.
Proposition~\ref{prop:zero-obstruction} treats the case of a zero initial profile.

\begin{assumption}[A convenient compact-kernel subclass]\label{ass:quantitative}
The law $h$ is supported in $[s_-,s_+]\times B_R$, where
$B_R:=\{x\in\R^d:|x|<R\}$ is the Euclidean ball centered at the origin,
$0<R<\infty$, and $0<s_-<s_+<\infty$. The function $\lambda_0$ satisfies the boundedness and uniform
continuity conditions of Assumption~\ref{ass:finite-mean}. 
\end{assumption}

Assumption~\ref{ass:quantitative} is used
for the optional spatial bounds with initial profiles of possibly infinite mass;
the weaker moment assumption below suffices for the bias analysis.

\begin{assumption}[Moments for the weak bias expansion]\label{ass:bias-moments}
In addition to Assumption~\ref{ass:finite-mean}, assume
\begin{equation}\label{eq:bias-moments}
 \E\left[S^2\right]<\infty,\qquad \E\left[|\Xi|^4\right]<\infty.
\end{equation}
Conditional spatial symmetry, when invoked, means
$h(s,z)=h(s,-z)$ almost everywhere.
\end{assumption}
Assumption~\ref{ass:quantitative} implies Assumption~\ref{ass:bias-moments}.
The former is used for the optional compact-kernel spatial bounds;
the latter suffices for both bias propositions in Section~\ref{sec:bias}.
In particular, by Cauchy--Schwarz inequality,
\begin{equation}\label{eq:bias-mixed-moment}
 \E\left[S|\Xi|^2\right]\leq\left(\E\left[S^2\right]\right)^{1/2}\left(\E\left[|\Xi|^4\right]\right)^{1/2}<\infty.
\end{equation}

\begin{table}[htbp]
\centering\small
\caption{Principal notation. Each scale and observable has one meaning in the finite-mean theory.}\label{tab:notation}
\begin{tabularx}{\textwidth}{@{}p{.27\textwidth}X@{}}
\toprule
Symbol & Meaning\\
\midrule
$(S,\Xi)$; $h,H,\overline H$ & Offspring delay and displacement; joint density, delay density and tail.\\
$\Pi^\eps$, $N^\eps$, $\NN_\eps$ & Prescribed past, microscopic descendants, normalized event measure.\\
$\delta_\eps$, $n_\eps$ & Generator scale $\eps^2$; effective population $\eps^{-(\alpha-1)}$.\\
$\lambda_0$, $\rho$, $\mu_\eps$, $\mu$ & Initial density, hydrodynamic density, exact mean, limiting event measure.\\
$\B$, $\LL_\kappa$, $\cQ_t$, $\cA$ & Spatial operator, backward operator, killed heat semigroup, resolvent.\\
$Z$, $Y_\eps$, $\zeta$, $v$ & Gaussian field, moderate fluctuations, additive Gaussian control, multiplicative branching control.\\
$J_{\rm sp}$, $\mathcal J_{\lambda_0}$ & Local variational action and finite-mass space-time LDP rate function.\\
$\tau_\eps$, $\delta_{\eps,\eta}$, $n_{\eps,\eta}$ & Heavy-tail time scale, generator scale, effective population.\\
\bottomrule
\end{tabularx}
\end{table}

\section{Main results for finite-mean delays}\label{sec:main-results}
All results in this section use the model in Assumption~\ref{ass:finite-mean}.
The finite-mass condition is added where needed.

The CLT, MDP, and LDP concern the same space--time event
measure $\NN_\eps$ under the scaling
$(t,x)\mapsto(\eps^2t,\eps x)$. Theorem~\ref{thm:full-spatial-finite}
is the macroscopic large-deviation counterpart of
Theorems~\ref{thm:clt} and \ref{thm:mdp}.

\subsection{Hydrodynamic limit, CLT, and MDP}

The deterministic reference limit follows from
\cite[Theorem~4.1, parts~(1) and~(2)]{GWZ2026}, with $m=m_1$,
$\NN_\eps=\mathcal N_\eps^a$, and $\rho=\rho^\kappa$ in that paper's
notation. Its kernel assumptions, prescribed Poisson past, and
near-critical scaling coincide with Assumption~\ref{ass:finite-mean}
and \eqref{eq:past}--\eqref{eq:event}.

\begin{proposition}[Hydrodynamic limit]\label{prop:hydrodynamic-reference}
Under Assumption~\ref{ass:finite-mean}, $\mu_\eps\to\mu$ vaguely on
$\Omega_T$ and $\NN_\eps\to\mu$ in probability in $\mathsf E$.
\end{proposition}

\cite[Theorem~4.1]{GWZ2026} states vague convergence on $(0,\infty)\times\R^d$.
To include $t=0$ in $\Omega_T$, let $f\in C_c(\Omega_T)$ have spatial
support in a compact set $K$, and choose a smooth temporal cutoff
$0\leq\theta_\eta\leq1$ vanishing near zero and equal to one for
$t\geq\eta$, where $0<\eta<T$.
Then $\theta_\eta f\in C_c((0,T)\times\R^d)$, and
\eqref{eq:normalized-mean-density-bound} gives 
\begin{equation}\label{eq:hydrodynamic-initial-strip}
 \sup_{0<\eps\leq1}\E\left[\left|\NN_\eps((1-\theta_\eta)f)\right|\right]
 +\left|\mu((1-\theta_\eta)f)\right|
 \leq2\|\lambda_0\|_\infty\|f\|_\infty|K|\eta.
\end{equation}
Apply \cite[Theorem~4.1]{GWZ2026} to $\theta_\eta f$ and then let $\eta\downarrow0$.
Equation~\eqref{eq:hydrodynamic-initial-strip} and Markov's inequality
give the asserted mean and probability limits for every $f\in C_c(\Omega_T)$;
a countable convergence-determining family gives convergence in $\mathsf E$.

The next theorem identifies the Gaussian fluctuations of the event measure around its exact mean and expresses their covariance through the backward resolvent.

\begin{theorem}[Fluctuations]\label{thm:clt}
Under Assumption~\ref{ass:finite-mean},
\begin{equation}\label{eq:clt}
 Z_\eps:=\sqrt{n_\eps}(\NN_\eps-\mu_\eps)
 \ \Longrightarrow\ Z
 \quad\text{in }\cD'(E_T),
\end{equation}
where $Z$ is the centered Gaussian distribution with covariance \eqref{eq:cov}. It can be defined by a space--time white noise $W$ through
\begin{equation}\label{eq:Gaussian}
 Z(f)=\int_0^T\int_{\R^d}\sqrt{\rho(t,x)}\cA  f(t,x)W(\dd t,\dd x).
\end{equation}
Equivalently, it is the zero-initial-data distributional solution of
\begin{equation}\label{eq:SPDE}
 m\partial_t Z=\B Z-\kappa Z+\sqrt\rho\,\dot W.
\end{equation}
The zero initial condition in this formulation is defined by \eqref{eq:Gaussian}.
\end{theorem}

In the next theorem, the same covariance determines the quadratic cost of fluctuations between the Gaussian scale and the hydrodynamic scale.

\begin{theorem}[Moderate deviations]\label{thm:mdp}
Under the same assumptions, let $b_\eps\to\infty$ and $b_\eps/\sqrt{n_\eps}\to0$. Then
\begin{equation}\label{eq:mdpvariable}
 Y_\eps:=\frac{\sqrt{n_\eps}}{b_\eps}(\NN_\eps-\mu_\eps)
\end{equation}
satisfies an LDP in $\cD'(E_T)$ with speed $b_\eps^2$ and good rate function
\begin{align}\label{eq:MDPrate}
 I_G(z)
 &=\sup_{f\in\cD(E_T)}\left\{z(f)-\frac12 C(f,f)\right\}\nonumber\\
 &=\inf\left\{\frac12\int_0^T\int_{\mathbb{R}^{d}} \left(\zeta(t,x)\right)^2\dd x\dd t:
 z(f)=\int_0^T\int_{\mathbb{R}^{d}}\sqrt\rho\,\cA  f\,\zeta\dd x\dd t,\,\,\, \text{for every }f\right\}.
\end{align}
An empty infimum is $+\infty$. In mild notation the constraint is
\begin{equation}\label{eq:MDPcontrol}
 m\partial_t z=\B z-\kappa z+\sqrt\rho\,\zeta,\qquad z_0=0.
\end{equation}
\end{theorem}

The proofs of Theorems~\ref{thm:clt} and \ref{thm:mdp} are in
Section~\ref{sec:clt-mdp}. The control $\zeta\in L^2(E_T)$ enters
the additive forcing in \eqref{eq:MDPcontrol}. The results use $\mu_\eps$ as their
center. Conditions for replacing it by $\mu$, including endpoint bias
corrections, are proved in Section~\ref{sec:bias}.

\subsection{The space--time event-measure LDP at finite mass}

For $\nu\in\mathsf E$, define the local variational action
\begin{equation}\label{eq:spatialJ}
 J_{\rm sp}(\nu):=\sup_{\phi\in C_c^\infty([0,T)\times\R^d)}
 \left\{-\ip\nu{\LL_\kappa\phi+\phi^2/2}
       -m\int_{\mathbb{R}^{d}}\lambda_0(x)\phi(0,x)\dd x\right\}.
\end{equation}

Under Assumption~\ref{ass:finite-mass}, the following theorem
identifies the rate function and its finite domain in terms of
$J_{\rm sp}$. Proposition~\ref{prop:zero-obstruction} gives the rate
when the initial profile vanishes.

\begin{theorem}[Large deviations]
\label{thm:full-spatial-finite}
Suppose Assumption~\ref{ass:finite-mass} holds. Then $\NN_\eps$ satisfies a full LDP on
$\mathsf E$, at speed $n_\eps$, with good rate function
\begin{equation}\label{eq:full-spatial-rate}
 \mathcal J_{\lambda_0}(\nu):=
 \begin{cases}
 J_{\rm sp}(\nu),&\nu([0,T)\times\R^d)<\infty,\\
 +\infty,&\nu([0,T)\times\R^d)=\infty.
 \end{cases}
\end{equation}
Equivalently, its finite domain consists precisely of
$\nu(\dd t,\dd x)=q(t,x)\dd t\dd x$ for which
$q\in C([0,T];L^1(\R^d))$, $q\geq0$, $q(0)=\lambda_0$, and there is a
measurable $v$ satisfying
\begin{equation}\label{eq:full-spatial-controlled}
 m\partial_tq=\B q-\kappa q+vq,
 \qquad \int_0^T\int_{\mathbb{R}^{d}} (v(t,x))^2q(t,x)\dd x\dd t<\infty.
\end{equation}
The equation is distributional, including the initial trace in $L^1$;
$v$ is an $L^2(q\dd t\dd x)$ equivalence class, defined
$q\dd t\dd x$-almost everywhere. Moreover,
\begin{equation}\label{eq:full-spatial-energy}
 \mathcal J_{\lambda_0}(q\dd t\dd x)
     =\frac12\int_0^T\int_{\mathbb{R}^{d}} (v(t,x))^2q(t,x)\dd x\dd t.
\end{equation}
\end{theorem}

The proof is presented in Section~\ref{sec:spatialLDP}. The microscopic
upper bound uses uniform generator convergence and the exponential
identity; the lower bound uses the density-path representation of finite-action measures, the explicit
approximation by smooth controls in Lemma~\ref{lem:finite-recovery}, and
a microscopic intensity change.

\section{Renewal estimates and analytic cluster transforms}\label{sec:analytic}
\subsection{Linear renewal estimates}
We use the operators and weighted norms introduced in Section~\ref{sec:model}.

The renewal bounds in the following lemma control the discrete backward resolvents and identify their heat-resolvent limit in the two norms used for cluster transforms.

\begin{lemma}[Finite-horizon renewal bounds]\label{lem:renewal}
For $q>0$ and $p=1,\infty$,
\begin{equation}\label{eq:renewalbound}
 \norm{\cA_\eps}_{q,p}\leq
 \frac{\delta_\eps}{1-a_\eps\E[e^{-q\delta_\eps S}]}
 \longrightarrow\frac1{\kappa+mq}.
\end{equation}
On continuous functions vanishing after $T$ and at spatial infinity,
$\cA_\eps f\to \cA  f$ uniformly. For compactly supported continuous $f$, the convergence also holds in $L^\infty([0,T];L^1(\R^d))$. If two families converge in both the uniform and $L^\infty([0,T];L^1(\R^d))$ norms and are uniformly bounded in the uniform norm, then their products converge in both norms. The convergence extends to limits $f$ in the closure of compact continuous functions in both norms: if $f_\eps\to f$ in those norms, then $\cA_\eps f_\eps\to \cA f$ in both norms.
\end{lemma}

\begin{proof}
Let $(S_\ell,\Xi_\ell)$, $\ell\geq 1$, be independent and indentically distributed
pairs with the joint law in Assumption~\ref{ass:finite-mean}, and let
$T_j:=\sum_{\ell=1}^jS_\ell$ and $X_j:=\sum_{\ell=1}^j\Xi_\ell$,
with $T_0:=X_0:=0$. Translation invariance of the spatial $L^p$ norm gives
\begin{align*}
 e^{qt}\left\|\cK_\eps^jf(t)\right\|_p
 &\leq \E\left[e^{-q\delta_\eps T_j}\right]\|f\|_{q,p}
   =p_\eps^j\|f\|_{q,p},\qquad
 p_\eps=\E\left[e^{-q\delta_\eps S}\right],\\
 \|\cA_\eps f\|_{q,p}
 &\leq\delta_\eps\sum_{j=0}^\infty(a_\eps p_\eps)^j\|f\|_{q,p}
   =\frac{\delta_\eps}{1-a_\eps p_\eps}\|f\|_{q,p}.
\end{align*}
Since $(1-e^{-q\delta_\eps S})/\delta_\eps\leq qS$ and $\E[S]=m$,
the dominated convergence theorem and \eqref{eq:scales} give
\[
 \frac{1-a_\eps p_\eps}{\delta_\eps}
 =\kappa_\eps+a_\eps\E\left[
             \frac{1-e^{-q\delta_\eps S}}{\delta_\eps}\right]
 \longrightarrow\kappa+mq,
\]
which proves \eqref{eq:renewalbound}.

To identify the resolvent limit, let
\[
 U_\eps:=\delta_\eps\sum_{j=0}^\infty a_\eps^j
                    \operatorname{Law}\left(\delta_\eps T_j,\eps X_j\right),
 \qquad U(\dd u,\dd z):=\frac1m e^{-\kappa u/m}\cP_u(0,\dd z)\dd u.
\]
These measures are locally finite in $u$. For $q>0$ and $\xi\in\R^d$, let
$D_\eps:=1-e^{-q\delta_\eps S}$. The spatial second-order expansion in
\eqref{eq:basic} gives
$\E\left[e^{\mathrm{i}\eps\xi\cdot\Xi}\right]
=1-\delta_\eps\xi^\top\Sigma\xi/2+o(\delta_\eps)$.
The possible dependence of $S$ and $\Xi$ is controlled by
\begin{align*}
 \frac{|\E[D_\eps(e^{\mathrm{i}\eps\xi\cdot\Xi}-1)]|}{\delta_\eps}
 &\leq\frac{\eps|\xi|R\E[D_\eps]}{\delta_\eps}
 +\frac{\eps|\xi|}{\delta_\eps}
       \left(\E\left[D_\eps^2\right]\right)^{1/2}
       \left(\E\left[|\Xi|^2\one_{\{|\Xi|>R\}}\right]\right)^{1/2}\\
 &\leq mq|\xi|R\eps
       +\sqrt{mq}|\xi|\left(\E\left[|\Xi|^2\one_{\{|\Xi|>R\}}\right]\right)^{1/2}.
\end{align*}
Here $\E\left[D_\eps^2\right]\leq\E\left[D_\eps\right]\leq mq\delta_\eps$.
Taking $\eps\downarrow0$ and then $R\uparrow\infty$ proves
\begin{equation}\label{eq:joint-renewal-transform}
 \varphi_\eps(q,\xi):=\E\left[e^{-q\delta_\eps S+\mathrm{i}\eps\xi\cdot\Xi}\right]
 =1-\delta_\eps\left(mq+\xi^\top\Sigma\xi/2\right)+o(\delta_\eps),
\end{equation}
and hence
\begin{equation}\label{eq:renewal-measure-transform}
 \int_0^\infty\int_{\R^d}e^{-qu+\mathrm{i}\xi\cdot z}U_\eps(\dd u,\dd z)
 =\frac{\delta_\eps}{1-a_\eps\varphi_\eps(q,\xi)}
 \longrightarrow\frac1{\kappa+mq+\xi^\top\Sigma\xi/2}.
\end{equation}

We verify tightness of the finite measures $e^{-qu}U_\eps$ explicitly.
Let $z_\eps:=\E\left[\Xi e^{-q\delta_\eps S}\right]=-\E[\Xi D_\eps]$. The
Cauchy--Schwarz inequality gives
$|z_\eps|^2\leq mq\delta_\eps\E\left[|\Xi|^2\right]$.
For $j\geq1$, expand the square of
$X_j=\sum_{\ell=1}^j\Xi_\ell$ and use
$e^{-q\delta_\eps T_j}=\prod_{\ell=1}^j e^{-q\delta_\eps S_\ell}$:
\begin{align}\label{eq:renewal-independent-pair-expansion}
 \E\left[e^{-q\delta_\eps T_j}|X_j|^2\right]
 &=\sum_{\ell=1}^j
       \E\left[|\Xi_\ell|^2\prod_{k=1}^j e^{-q\delta_\eps S_k}\right]
       +\sum_{\substack{1\leq\ell,k\leq j\\\ell\ne k}}
       \E\left[(\Xi_\ell\cdot\Xi_k)
                    \prod_{i=1}^j e^{-q\delta_\eps S_i}\right]\nonumber\\
 &=j p_\eps^{j-1}\E\left[|\Xi|^2e^{-q\delta_\eps S}\right]
       +j(j-1)p_\eps^{j-2}|z_\eps|^2.
\end{align}
The pairs $(S_\ell,\Xi_\ell)$ are independent across $\ell$.
Each diagonal summand therefore factors into one weighted second
moment and $j-1$ factors $p_\eps$. Each off-diagonal summand factors
into two weighted first moments $z_\eps$ and $j-2$ factors $p_\eps$.
There are $j$ diagonal and $j(j-1)$ off-diagonal summands; for $j=1$
the second sum is empty. Dependence within each pair is retained in
the weighted moments.
Consequently, using $a_\eps\leq1$ and $1-p_\eps\sim mq\delta_\eps$,
\begin{align}\label{eq:renewal-spatial-moment}
 \int_0^\infty\int_{\R^d}e^{-qu}|z|^2U_\eps(\dd u,\dd z)
 &\leq\frac{\delta_\eps^2\E[|\Xi|^2]}{(1-p_\eps)^2}
       +\frac{2\delta_\eps^2|z_\eps|^2}{(1-p_\eps)^3}\leq C_q.
\end{align}
The temporal tail is bounded by
\begin{equation}\label{eq:renewal-temporal-tail}
 \int_L^\infty\int_{\R^d}e^{-qu}U_\eps(\dd u,\dd z)
 \leq e^{-qL/2}\frac{\delta_\eps}
                   {1-a_\eps\E[e^{-q\delta_\eps S/2}]}.
\end{equation}
For every sequence $\eps_n\downarrow0$, the finite measures
$e^{-qu}U_{\eps_n}(\dd u,\dd z)$ have uniformly bounded total masses
by \eqref{eq:renewal-measure-transform} with $\xi=0$ and are tight
by \eqref{eq:renewal-spatial-moment}--\eqref{eq:renewal-temporal-tail}.
They therefore have a weakly convergent subsequence on
$[0,\infty)\times\R^d$.
Equation~\eqref{eq:renewal-measure-transform}, also with $q$ replaced by
$q+r$ for $r\geq0$, identifies its Laplace--Fourier transform with that
of $e^{-qu}U$. Uniqueness of this transform identifies the limit.

By \eqref{eq:K} and \eqref{eq:resolvent},
\[
 \cA_\eps f(t,x)=\int_0^\infty\int_{\R^d}f(t+u,x+z)U_\eps(\dd u,\dd z),
 \quad \cA f(t,x)=\int_0^\infty\int_{\R^d}f(t+u,x+z)U(\dd u,\dd z).
\]
Let $\widetilde U_\eps(\dd u,\dd z):=e^{-qu}U_\eps(\dd u,\dd z)$,
$\widetilde U(\dd u,\dd z):=e^{-qu}U(\dd u,\dd z)$, and
\[
 F_{t,x}(u,z):=e^{qu}f(t+u,x+z),\qquad
 0\leq t\leq T,\quad x\in\R^d.
\]
Extending $f$ by zero for $t\geq T$ gives a uniformly continuous
function on $[0,\infty)\times\R^d$, and
$|F_{t,x}|\leq M_f:=e^{qT}\|f\|_\infty$ because $F_{t,x}=0$ for
$u>T$. On $K_{L,R}:=[0,L]\times\overline B_R$, these functions
have the common modulus
\[
 \left|F_{t,x}(u,z)-F_{t,x}(u',z')\right|
 \leq e^{qL}\omega_f\left(|u-u'|+|z-z'|\right)
       +q e^{qL}\|f\|_\infty\left|u-u'\right|,
\]
where $\omega_f$ is a modulus of uniform continuity of $f$.
For $\zeta>0$, choose a finite mesh of $K_{L,R}$ on which this
oscillation is at most $\zeta/3$. Partition the bounded vectors of
function values on the mesh into finitely many sets of diameter
at most $\zeta/3$ in the maximum norm. Choosing one member of the
family from each nonempty set gives bounded continuous functions
$F_1,\ldots,F_N$ on $[0,\infty)\times\R^d$ such that
\begin{equation}\label{eq:renewal-finite-net}
 \sup_{\substack{0\leq t\leq T\\x\in\R^d}}
       \min_{1\leq i\leq N}\|F_{t,x}-F_i\|_{C(K_{L,R})}\leq\zeta.
\end{equation}
By \eqref{eq:renewalbound}, the sum of the total masses of
$\widetilde U_\eps$ and $\widetilde U$ is bounded by a constant
$M_q$ for all sufficiently small $\eps$. Equations
\eqref{eq:renewal-spatial-moment} and \eqref{eq:renewal-temporal-tail},
and their limits, give
\begin{equation}\label{eq:renewal-tilted-tail}
 \left(\widetilde U_\eps+\widetilde U\right)(K_{L,R}^c)
 \leq C_q\left(R^{-2}+e^{-qL/2}\right).
\end{equation}
Let $\Delta_\eps:=\widetilde U_\eps-\widetilde U$.
Using \eqref{eq:renewal-finite-net} on $K_{L,R}$ and the bound
$|F_{t,x}-F_i|\leq2M_f$ on its complement yields
\begin{align}\label{eq:renewal-uniform-convolution}
 \|\cA_\eps f-\cA f\|_\infty
 &=\sup_{\substack{0\leq t\leq T\\x\in\R^d}}
       \left|\int_{[0,\infty)\times\R^d}F_{t,x}\dd\Delta_\eps\right|
       \nonumber\\
 &\leq\max_{1\leq i\leq N}
       \left|\int_{[0,\infty)\times\R^d}F_i\dd\Delta_\eps\right|
       +M_q\zeta+2M_fC_q\left(R^{-2}+e^{-qL/2}\right).
\end{align}
For fixed $L,R,\zeta$, the finite family in
\eqref{eq:renewal-finite-net} consists of bounded continuous functions.
Weak convergence of $\widetilde U_\eps$ to $\widetilde U$ therefore gives
\[
 \max_{1\leq i\leq N}
 \left|\int_{[0,\infty)\times\R^d}F_i\dd\Delta_\eps\right|
 \longrightarrow0.
\]
This is the first term on the right of
\eqref{eq:renewal-uniform-convolution}. Letting $L,R\to\infty$ and then
$\zeta\downarrow0$ proves $\|\cA_\eps f-\cA f\|_\infty\to0$.

For compactly supported $f$, translation of its support and
\eqref{eq:renewal-spatial-moment} give
\begin{equation}\label{eq:renewal-output-moment}
 \sup_{0\leq t\leq T}\int_{\R^d}\left(1+|x|^2\right)\cA_\eps|f|(t,x)\dd x
 \leq C_f e^{qT}\int_0^\infty\int_{\R^d}
                        e^{-qu}\left(1+|z|^2\right)U_\eps(\dd u,\dd z)
 \leq C_{f,T,q}.
\end{equation}
The analogous bound holds for $\cA |f|$. For every $R>0$,
\[
 \sup_{0\leq t\leq T}\left\|\cA_\eps f(t)-\cA f(t)\right\|_1
 \leq |B_R|\left\|\cA_\eps f-\cA f\right\|_\infty+C_{f,T,q}R^{-2}.
\]
Letting $\eps\downarrow0$ and then $R\uparrow\infty$ proves the $L^1$
convergence. Finally, for $p=1,\infty$,
\[
 \|u_\eps v_\eps-uv\|_{L^\infty([0,T];L^p(\R^d))}
 \leq\|v_\eps\|_\infty\|u_\eps-u\|_{L^\infty([0,T];L^p(\R^d))}
       +\|u\|_\infty\|v_\eps-v\|_{L^\infty([0,T];L^p(\R^d))},
\]
and
\[
 \left\|\cA_\eps f_\eps-\cA f\right\|_{L^\infty([0,T];L^p(\R^d))}
 \leq M_T\|f_\eps-f\|_{L^\infty([0,T];L^p(\R^d))}
       +\left\|\cA_\eps f-\cA f\right\|_{L^\infty([0,T];L^p(\R^d))}.
\]
The last estimate extends from compact continuous $f$ to their common closure in
two norms $\|\cdot\|_\infty$ and $\|\cdot\|_{L^\infty([0,T];L^1(\R^d))}$
by approximation and \eqref{eq:renewalbound}.
This completes the proof.
\end{proof}

\subsection{Positive and complex sources}

Let $\mathcal C_\eps(t,x)$ denote the random genealogical tree of the
full descendant cluster rooted at the macroscopic point $(t,x)$, with
the root included. For each $u\in\mathcal C_\eps(t,x)$, let $(t_u,x_u)$
be its macroscopic time and position. Define
\begin{equation}\label{eq:cluster-transform-definition}
 F_{\eps,f}(t,x):=\E_{t,x}\left[\exp\left\{\delta_\eps^2\sum_{u\in\mathcal C_\eps(t,x)}f(t_u,x_u)\right\}\right],
 \qquad v_{\eps,f}:=\frac{F_{\eps,f}-1}{\delta_\eps}.
\end{equation}
Here $f$ is extended by zero after $T$, so only events with $t_u\leq T$
contribute. We call $F_{\eps,f}$ the cluster transform and $v_{\eps,f}$
its rescaled version.
This definition is initially justified for sufficiently small positive sources; the following lemma establishes that justification. Conditioning on the Poisson offspring of the root gives
\begin{equation}\label{eq:poisson-offspring}
 F_{\eps,f}(t,x)=\exp\left\{\delta_\eps^2 f(t,x)
       +a_\eps\cK_\eps\left(F_{\eps,f}-1\right)(t,x)\right\}.
\end{equation}
Substituting $F_{\eps,f}=1+\delta_\eps v_{\eps,f}$ into
\eqref{eq:poisson-offspring} and using \eqref{eq:K} yields
\begin{equation}\label{eq:positivecluster}
 \log(1+\delta_\eps v)=\delta_\eps^2 f+a_\eps\delta_\eps \cK_\eps v,
 \qquad
 v=\cA_\eps\{f+G_\eps(v)\},
\end{equation}
where
\begin{equation}\label{eq:G}
 G_\eps(z):=\frac{\delta_\eps z-\log(1+\delta_\eps z)}{\delta_\eps^2}
 =\frac{z^2}{2}-\frac{\delta_\eps z^3}{3}+\mathcal O\left(\delta_\eps^2|z|^4\right).
\end{equation}
The quadratic term in \eqref{eq:G} is positive. For comparison, let
$g\geq0$ be continuous, compactly supported in $[0,T)\times\R^d$,
and sufficiently small in the uniform norm, and let
$u_{\eps,g}:=(1-F_{\eps,-g})/\delta_\eps=-v_{\eps,-g}\geq0$.
Substitution into \eqref{eq:positivecluster} gives
\begin{equation}\label{eq:nonnegative-laplace-discrete}
 u_{\eps,g}=\cA_\eps\left\{g-G_\eps(-u_{\eps,g})\right\}.
\end{equation}
The convergence proved in Lemma~\ref{lem:analytic} below and
\eqref{eq:G} identify its limit $u_g=-v_{-g}\geq0$ through
\begin{equation}\label{eq:nonnegative-laplace-cost}
 u_g=\cA\left(g-u_g^2/2\right),\qquad
 \LL_\kappa u_g-u_g^2/2+g=0,\qquad u_g(T)=0.
\end{equation}
Thus the nonnegative log-Laplace equation
\eqref{eq:nonnegative-laplace-cost} has the damping term $-u_g^2/2$.

The next lemma proves existence of positive exponential moments on a common source neighborhood and provides the analytic bounds needed to differentiate the cluster transform.

\begin{lemma}[Uniform analytic neighborhood]\label{lem:analytic}
Fix $T$ and a compact spatial support set $K$. There exists $r_T>0$, independent of sufficiently small $\eps$ and of both the location and volume of $K$, such that $F_{\eps,f}$ in \eqref{eq:cluster-transform-definition}
is finite for real bounded measurable $f$ with $\norm{f}_\infty<r_T$
and $\supp f\subset[0,T]\times K$.
For any finite family of bounded sources $f_1,\ldots,f_k$ with this support, let
\begin{equation}\label{eq:cluster-parameter-domain}
 f_\theta:=\sum_{\ell=1}^k\theta_\ell f_\ell,\qquad
 \mathcal U:=\left\{\theta\in\mathbb C^k:\|f_\theta\|_\infty<r_T\right\}.
\end{equation}
The maps $\theta\mapsto v_{\eps,f_\theta}$ and
$\theta\mapsto F_{\eps,f_\theta}=1+\delta_\eps v_{\eps,f_\theta}$
are holomorphic on $\mathcal U$, with locally uniform convergence of their
power series in the uniform norm; the same holds for $v_{\eps,f_\theta}$
in $L^\infty([0,T];L^1(\R^d))$.
For every source with $\|f\|_\infty<r_T$, with $C_T$ also independent
of the location and volume of $K$,
\begin{equation}\label{eq:vbounds}
 \norm{v_{\eps,f}}_\infty\leq C_T\norm{f}_\infty,
 \qquad
 \sup_{0\leq t\leq T}\norm{v_{\eps,f}(t)}_1\leq C_T\sup_{0\leq t\leq T}\norm{f(t)}_1.
\end{equation}
For continuous compactly supported sources $f$ that vanish at $T$ and
satisfy $\|f\|_\infty<r_T$, the rescaled transforms $v_{\eps,f}$,
which solve \eqref{eq:positivecluster},
converge in $\|\cdot\|_\infty$ and
$\|\cdot\|_{L^\infty([0,T];L^1(\R^d))}$ to the small bounded solution of
\begin{equation}\label{eq:positivePDE}
 v_f=\cA \left(f+v_f^2/2\right),\qquad
 m\partial_t v_f+\B v_f-\kappa v_f+v_f^2/2+f=0,
 \qquad v_f(T)=0.
\end{equation}
\end{lemma}

\begin{proof}
Let $M_T\geq1$ be a common bound for $\cA_\eps$ and $\cA$ in $L^\infty$ and
$L^\infty([0,T];L^1(\R^d))$, obtained by fixing $q>0$ in
\eqref{eq:renewalbound} and using \eqref{eq:resolvent}.
For $\delta_\eps R\leq1/2$, differentiation of \eqref{eq:G} gives
\[
 G_\eps'(z)=\frac{z}{1+\delta_\eps z},\qquad
 |G_\eps(z)|\leq |z|^2,\qquad
 |G_\eps(z)-G_\eps(w)|\leq2R|z-w|\quad(|z|,|w|\leq R).
\]
Choose
\begin{equation}\label{eq:cluster-source-radius}
 0<r_T\leq\min\{(8M_T^2)^{-1},(4M_T)^{-1}\},\qquad
 R_T:=2M_Tr_T,\qquad
 2M_TR_T\leq\frac12,\quad\delta_\eps R_T\leq\frac12.
\end{equation}
For a fixed $f$, let $R:=2M_T\|f\|_\infty<R_T$. For
$\Phi_\eps(v)=\cA_\eps(f+G_\eps(v))$, we have
\begin{align}\label{eq:cluster-contraction}
 \|\Phi_\eps(v)\|_\infty
 &\leq M_T\|f\|_\infty+M_TR^2\leq R,\nonumber\\
 \|\Phi_\eps(v)-\Phi_\eps(w)\|_{L^\infty([0,T];L^p(\R^d))}
 &\leq2M_TR\|v-w\|_{L^\infty([0,T];L^p(\R^d))}\nonumber\\
 &\leq\frac12\|v-w\|_{L^\infty([0,T];L^p(\R^d))},\qquad p=1,\infty.
\end{align}
The Banach fixed-point theorem gives a unique fixed point in the uniform
ball. Starting its iteration at zero and applying the $p=1$ estimate gives
\begin{align*}
 \|v_{\eps,f}\|_{L^\infty([0,T];L^1(\R^d))}
 &\leq M_T\|f\|_{L^\infty([0,T];L^1(\R^d))}
       +\frac12\|v_{\eps,f}\|_{L^\infty([0,T];L^1(\R^d))}.
\end{align*}
Absorbing the last term gives
\[
\|v_{\eps,f}\|_{L^\infty([0,T];L^1(\R^d))}
 \leq2M_T\|f\|_{L^\infty([0,T];L^1(\R^d))}.
\]
The source radius $r_T$ depends on $T$, $m$, and the bounded family
$\kappa_\eps$. The solution-ball radius $R=2M_T\|f\|_\infty$ and the
contraction constant $c_{T,f}:=2M_TR=4M_T^2\|f\|_\infty\leq1/2$
also depend on $\|f\|_\infty$. These choices are independent of the
location and volume of $K$.
To prove holomorphy, use $f_\theta$ and $\mathcal U$ from
\eqref{eq:cluster-parameter-domain} and define the Picard iterates
\begin{equation}\label{eq:cluster-analytic-picard}
 v_\eps^{[0]}(\theta)=0,\qquad
 v_\eps^{[n+1]}(\theta)
 =\cA_\eps\left(f_\theta+G_\eps(v_\eps^{[n]}(\theta))\right),
 \qquad n\geq0.
\end{equation}
The logarithm in \eqref{eq:G} is the branch with $\log1=0$.
On $|z|\leq R_T$, it has the absolutely convergent expansion
\[
 G_\eps(z)=\sum_{\ell=2}^\infty
       \frac{(-1)^\ell\delta_\eps^{\ell-2}}{\ell}z^\ell.
\]
Multiplication is continuous in
$X:=L^\infty([0,T]\times\R^d)\cap L^\infty([0,T];L^1(\R^d))$,
equipped with the sum of its two norms, because
$\|uv\|_{0,1}\leq\|u\|_\infty\|v\|_{0,1}$.
Thus \eqref{eq:cluster-analytic-picard} and the power series for $G_\eps$
prove inductively that each $v_\eps^{[n]}$ is holomorphic in
$\theta=(\theta_1,\ldots,\theta_k)\in\mathcal U$ as an $X$-valued map.
For a compact set $Q\subset\mathcal U$, \eqref{eq:cluster-contraction} gives
\[
 \sup_{\theta\in Q}
 \left\|v_\eps^{[n+1]}(\theta)-v_\eps^{[n]}(\theta)\right\|_{0,p}
 \leq 2^{-n}M_T\sup_{\theta\in Q}\|f_\theta\|_{0,p},
 \qquad p=1,\infty.
\]
The iterates therefore converge locally uniformly in both norms.
Passing to the limit in Cauchy's integral formula on each closed polydisc
contained in $\mathcal U$ proves holomorphy of the fixed point.

Let $|u|$ be the generation of $u$, with the root at generation zero.
For $j\geq0$, define the generation-truncated tree and its transform by
\begin{equation}\label{eq:cluster-generation-truncation}
 \begin{aligned}
 \mathcal C_\eps^{(j)}(t,x)
 &:=\left\{u\in\mathcal C_\eps(t,x):|u|\leq j\right\},\\
 F_{\eps,f}^{(j)}(t,x)
 &:=\E_{t,x}\left[\exp\left\{\delta_\eps^2
       \sum_{u\in\mathcal C_\eps^{(j)}(t,x)}f(t_u,x_u)\right\}\right].
 \end{aligned}
\end{equation}
Let $F_{\eps,f}^{(-1)}:=1$.
For a Poisson random measure $\mathcal P$ on a space $E$ with finite
intensity $\nu$, conditioning on its Poisson number of atoms gives, for bounded
complex $\psi$,
\begin{equation}\label{eq:offspring-poisson-product}
 \E\left[\prod_{y\in\mathcal P}\psi(y)\right]
 =e^{-\nu(E)}\sum_{n=0}^\infty\frac1{n!}
           \left(\int_E\psi(y)\nu(\dd y)\right)^n
 =\exp\left\{\int_E(\psi(y)-1)\nu(\dd y)\right\}.
\end{equation}
The empty product equals one. The root's offspring intensity is the
pushforward of $a_\eps h(s,z)\dd s\dd z$ under
$(s,z)\mapsto(t+\delta_\eps s,x+\eps z)$ and has total mass $a_\eps$.
Conditioning on these offspring, independence of their descendant trees
and \eqref{eq:offspring-poisson-product}, with
$\psi(t',x')=F_{\eps,f}^{(j-1)}(t',x')$, give
the recursion
\begin{equation}\label{eq:cluster-truncated-recursion}
 F_{\eps,f}^{(j)}
 =\exp\left\{\delta_\eps^2f+a_\eps\cK_\eps
                    \left(F_{\eps,f}^{(j-1)}-1\right)\right\},
 \qquad j\geq0.
\end{equation}
Induction in \eqref{eq:cluster-truncated-recursion} makes each
finite-generation transform finite for every bounded real source.
For $f\geq0$, abbreviate $F_j:=F_{\eps,f}^{(j)}$.
The fixed point is nonnegative. Exponentiating the first identity in
\eqref{eq:positivecluster} shows that
$\widehat F=1+\delta_\eps v_{\eps,f}$ satisfies
\begin{equation}\label{eq:cluster-fixed-transform}
 \widehat F(t,x)=\exp\left\{\delta_\eps^2 f(t,x)
       +a_\eps\cK_\eps\left(\widehat F-1\right)(t,x)\right\}.
\end{equation}
Induction in \eqref{eq:cluster-truncated-recursion}, using
\eqref{eq:cluster-fixed-transform}, gives
$1\leq F_j\leq F_{j+1}\leq\widehat F$.
Since $\mathcal C_\eps^{(j)}\uparrow\mathcal C_\eps$ and $f\geq0$,
the monotone convergence theorem gives
\begin{equation}\label{eq:nonnegative-cluster-limit}
 \begin{aligned}
 F_\infty(t,x)
 &:=\lim_{j\to\infty}F_j(t,x)\\
 &=\E_{t,x}\left[\lim_{j\to\infty}
       \exp\left\{\delta_\eps^2\sum_{u\in\mathcal C_\eps^{(j)}(t,x)}
                         f(t_u,x_u)\right\}\right]\\
 &=F_{\eps,f}(t,x)\leq\widehat F(t,x)<\infty.
 \end{aligned}
\end{equation}
Let $w_\infty:=(F_\infty-1)/\delta_\eps$. Then
$0\leq w_\infty\leq v_{\eps,f}$, so $\|w_\infty\|_\infty\leq R$.
Monotone convergence inside $\cK_\eps$ in
\eqref{eq:cluster-truncated-recursion} gives
\[
 F_\infty=\exp\left\{\delta_\eps^2f+a_\eps\cK_\eps(F_\infty-1)\right\},\qquad
 \left(\Id-a_\eps\cK_\eps\right)w_\infty
       =\delta_\eps\left\{f+G_\eps(w_\infty)\right\}.
\]
Applying the convergent resolvent series \eqref{eq:K}, justified by
\eqref{eq:renewalbound}, yields
$w_\infty=\cA_\eps\{f+G_\eps(w_\infty)\}$.
Uniqueness in \eqref{eq:cluster-contraction} gives
$F_\infty=\widehat F$.
For signed or complex sources in the class
\begin{equation}\label{eq:complex-source-class}
 f:[0,T]\times\R^d\longrightarrow\mathbb C,\qquad
 \supp f\subset[0,T]\times K,\qquad \|f\|_\infty<r_T,
\end{equation}
we have the domination
\begin{equation}\label{eq:cluster-complex-domination}
 \left|\exp\left\{\delta_\eps^2\sum_{u\in\mathcal C_\eps(t,x)}f(t_u,x_u)\right\}\right|
 \leq\exp\left\{\delta_\eps^2\sum_{u\in\mathcal C_\eps(t,x)}|f(t_u,x_u)|\right\}.
\end{equation}
The expectation of the right side of
\eqref{eq:cluster-complex-domination} is
\[
F_{\eps,|f|}(t,x)\leq1+\delta_\eps C_T\|f\|_\infty<\infty,
\]
by \eqref{eq:nonnegative-cluster-limit}
 and \eqref{eq:vbounds}.
The dominated convergence theorem applied to the transforms of
$\mathcal C_\eps^{(j)}(t,x)$ in
\eqref{eq:cluster-generation-truncation} allows $j\to\infty$ in
\eqref{eq:cluster-truncated-recursion}. It proves
\eqref{eq:poisson-offspring} for every source in
\eqref{eq:complex-source-class}. Moreover,
$|F_{\eps,f}-1|\leq F_{\eps,|f|}-1$, so the resulting rescaled
transform lies in the fixed-point ball of \eqref{eq:cluster-contraction}.
The bound
\[
 \left|\delta_\eps^2f+a_\eps\delta_\eps\cK_\eps v_{\eps,f}\right|
 \leq\delta_\eps^2r_T+\delta_\eps R_T\leq\frac38
\]
follows from \eqref{eq:cluster-source-radius} and $M_T\geq1$.
Thus taking the branch of the logarithm fixed above in
\eqref{eq:poisson-offspring} gives \eqref{eq:positivecluster};
uniqueness identifies the signed or complex transform with the analytic
fixed point constructed in \eqref{eq:cluster-analytic-picard}.
For the complex parameter vector $\theta\in\mathbb C^k$ in
\eqref{eq:cluster-parameter-domain}, consider
$f_\theta$ in \eqref{eq:complex-source-class}.
For every compact parameter set $Q\subset\mathcal U$, choose $r'<r_T$ with
$\sup_{\theta\in Q}\|f_\theta\|_\infty<r'$.
Replacing $|f|$ in \eqref{eq:cluster-complex-domination} by
$r'\one_{[0,T]\times K}$ gives an integrable common majorant.
The recursion \eqref{eq:cluster-truncated-recursion} shows that
$F_{\eps,f_\theta}^{(j)}$ in \eqref{eq:cluster-generation-truncation}
is analytic in $\theta$. Almost surely the full cluster is finite, so
the truncated exponential agrees with the full exponential for all
sufficiently large $j$, simultaneously for all $\theta$.
Dominated convergence theorem applied to the supremum of their difference
over each compact parameter set gives locally uniform convergence
of $F_{\eps,f_\theta}^{(j)}$ to $F_{\eps,f_\theta}$. The limit is
therefore analytic.

The same contraction argument with $\cA$ and $z^2/2$ gives the
small fixed point $v_f=\cA(f+v_f^2/2)$ in \eqref{eq:positivePDE}.
The integral form of
\eqref{eq:G} gives, for $|z|\leq R$,
\begin{equation}\label{eq:G-quadratic-remainder}
 \left|G_\eps(z)-z^2/2\right|\leq C\delta_\eps|z|^3.
\end{equation}
Subtracting \eqref{eq:positivePDE} from the fixed-point identity in
\eqref{eq:positivecluster} gives
\begin{align}\label{eq:cluster-resolvent-difference}
 v_{\eps,f}-v_f
 &=\cA_\eps\left\{G_\eps(v_{\eps,f})-G_\eps(v_f)\right\}
   +(\cA_\eps-\cA)\left(f+v_f^2/2\right)\nonumber\\
 &\quad+\cA_\eps\left\{G_\eps(v_f)-v_f^2/2\right\}.
\end{align}
Applying \eqref{eq:cluster-contraction} and
\eqref{eq:G-quadratic-remainder} to
\eqref{eq:cluster-resolvent-difference}, for $p=1,\infty$, yields
\begin{align*}
 \|v_{\eps,f}-v_f\|_{L^\infty([0,T];L^p(\R^d))}
 &\leq\frac12\|v_{\eps,f}-v_f\|_{L^\infty([0,T];L^p(\R^d))}
 \\
 &\qquad\qquad+\left\|(\cA_\eps-\cA )\left(f+v_f^2/2\right)\right\|_{L^\infty([0,T];L^p(\R^d))}\\
 &\qquad\qquad\qquad\qquad+C\delta_\eps M_TR^2\|v_f\|_{L^\infty([0,T];L^p(\R^d))}.
\end{align*}
The Picard iterates for $v_f$ belong to the common closure of compact
continuous functions in $\|\cdot\|_\infty$ and
$\|\cdot\|_{L^\infty([0,T];L^1(\R^d))}$.
Indeed, \eqref{eq:resolvent} preserves continuity and the terminal value
zero, and the limiting version of \eqref{eq:renewal-output-moment}
permits spatial cutoff approximation of $\cA g$ for compact continuous $g$.
The bounds in Lemma~\ref{lem:renewal} extend this property to the closure,
and its product estimate preserves it under squaring.
Consequently $f+v_f^2/2$ belongs to this closure, and
Lemma~\ref{lem:renewal}, applied in these two explicit norms,
shows that the second term tends to zero. The third term tends to zero by
\eqref{eq:G-quadratic-remainder} and $\delta_\eps\to0$.
Absorbing the first term proves
$v_{\eps,f}\to v_f$ in both the uniform norm and
$L^\infty([0,T];L^1(\R^d))$.
Apply $-\LL_\kappa$ to $v_f=\cA(f+v_f^2/2)$ in
\eqref{eq:positivePDE}. The inverse and terminal identities in
\eqref{eq:resolvent} give
\[
 -\LL_\kappa v_f
 =-\LL_\kappa\cA\left(f+v_f^2/2\right)=f+v_f^2/2,
 \qquad v_f(T)=\cA\left(f+v_f^2/2\right)(T)=0.
\]
The equation holds in the distributional sense, and classically
for smooth sources, which proves the remaining assertions in
\eqref{eq:positivePDE}.
\end{proof}

\subsection{Averaging over the prescribed past}
Define, for $r<0$,
\begin{equation}\label{eq:Jpast}
 \cJ_\eps w(r,x):=\int_0^\infty\int_{\R^d}
 h(q-r,z)w(\delta_\eps q,x+\eps z)\dd z\dd q,
 \quad
 \B_\eps w:=\int_{-\infty}^0\int_{\R^d}\lambda_0(x)\cJ_\eps w(r,x)\dd x\dd r.
\end{equation}
Let $\Lambda_\eps(f)$ denote the log moment generating functional of the unit-effective-population family ($n_\eps=1$, with the same $a_\eps$): its past density is $\eps^{d-2}\lambda_0(\eps y)$ and its mass per event is $\delta_\eps^2$. The independent cluster marks form a Poisson random measure by the marking theorem in \cite[Theorem~5.6]{LP}. Applying the Poisson exponential formula in \cite[Exercise~3.4]{LP} to these marks gives
\begin{equation}\label{eq:Lambdaexact}
 \Lambda_\eps(f)=\frac1{\delta_\eps}\int_{-\infty}^0\int_{\R^d}\lambda_0(x)
 \left[\exp\left\{a_\eps\delta_\eps \cJ_\eps v_{\eps,f}(r,x)\right\}-1\right]\dd x\dd r.
\end{equation}
The high-density model satisfies the \emph{exact} identity
\begin{equation}\label{eq:superposition}
 \log\E\left[e^{n_\eps\ip{\NN_\eps}{f}}\right]=n_\eps\Lambda_\eps(f).
\end{equation}
It holds for every positive real value of $n_\eps$.

Changing variables $S=q-r$ in \eqref{eq:Jpast} gives the crossing representation
\begin{equation}\label{eq:past-crossing-representation}
 \B_\eps w=\int_{\R^d}\lambda_0(x)
                 \E\left[\int_0^S w(\delta_\eps s,x+\eps\Xi)\dd s\right]\dd x.
\end{equation}
Fubini's theorem is justified by
$\E\left[\int_0^S\|w(\delta_\eps s)\|_1\dd s\right]\leq m\|w\|_{0,1}$.
Moreover, \eqref{eq:Jpast} gives
\[
|\cJ_\eps w(r,x)|\leq\norm{w}_\infty\overline H(-r)\leq\norm{w}_\infty.
\]
For $v=v_{\eps,f}$, the exponential Taylor remainder therefore satisfies
\[
 \left|\Lambda_\eps(f)-a_\eps \B_\eps v\right|
 \leq\frac{a_\eps^2\delta_\eps}{2}e^{\delta_\eps\|v\|_\infty}
          \|v\|_\infty \B_\eps|v|
 \leq\frac{m\delta_\eps}{2}e^{\delta_\eps\|v\|_\infty}
          \|\lambda_0\|_\infty\|v\|_\infty\|v\|_{0,1}.
\]
This proves the bounds
\begin{equation}\label{eq:Bbounds}
 |\B_\eps w|\leq m\norm{\lambda_0}_\infty\sup_{0\leq t\leq T}\norm{w(t)}_1,
 \qquad
 \Lambda_\eps(f)=a_\eps \B_\eps v_{\eps,f}
 +\mathcal O\left(\delta_\eps\norm{v_{\eps,f}}_\infty\sup_{0\leq t\leq T}\norm{v_{\eps,f}(t)}_1\right).
\end{equation}
If $w_\eps\to w$ in these two norms and $t\mapsto w(t,\cdot)$ is continuous at zero in $L^1_x$, then
\begin{equation}\label{eq:Blimit}
 \B_\eps w_\eps\longrightarrow m\int_{\R^d}\lambda_0(x)w(0,x)\dd x.
\end{equation}
To justify this limit, let
$D_R:=\{S\leq R,\ |\Xi|\leq R\}$ and
$\omega_w(r):=\sup_{|z|\leq r}\|w(0,\cdot+z)-w(0)\|_1$.
Equation~\eqref{eq:past-crossing-representation} gives
\begin{align*}
 \left|\B_\eps w_\eps-m\int_{\R^d}\lambda_0w(0)\dd x\right|
 &\leq m\|\lambda_0\|_\infty\|w_\eps-w\|_{0,1}\\
 &\quad+m\|\lambda_0\|_\infty
       \left[\sup_{0\leq t\leq\delta_\eps R}\|w(t)-w(0)\|_1
                                   +\omega_w(\eps R)\right]\\
 &\quad\quad+2\|\lambda_0\|_\infty\|w\|_{0,1}
                              \E\left[S\one_{D_R^c}\right].
\end{align*}
For fixed $R$, the first two lines tend to zero by the assumed norm
convergence, continuity at zero, and spatial translation continuity
in $L^1$. The last line tends to zero as $R\uparrow\infty$ by
$\E[S]=m<\infty$. These limits prove \eqref{eq:Blimit} for the
joint law in Assumption~\ref{ass:finite-mean}.

Differentiating the transforms of the offspring clusters and the prescribed Poisson past separates the limiting branching covariance from the smaller contribution of the random Poisson past, which leads to the following proposition.

\begin{proposition}[Cumulants and the contribution of the prescribed past]\label{prop:cumulants}
For every pair of smooth compactly supported real test functions,
\begin{align}\label{eq:cumulants}
 D\Lambda_\eps(0)f&=\ip{\mu_\eps}{f},\nonumber\\
 C_\eps(f,g):=D^2\Lambda_\eps(0)(f,g)
 &=a_\eps \B_\eps \cA_\eps[(\cA_\eps f)(\cA_\eps g)]\nonumber\\
 &\quad+a_\eps^2\delta_\eps\int_{-\infty}^0\int_{\R^d}\lambda_0(x)
 (\cJ_\eps \cA_\eps f)(r,x)(\cJ_\eps \cA_\eps g)(r,x)\dd x\dd r,\\
 C_\eps(f,g)&\longrightarrow C(f,g).\nonumber
\end{align}
Moreover, for test functions supported in a fixed compact set,
\begin{equation}\label{eq:third}
 \Lambda_\eps(zf)=z\ip{\mu_\eps}{f}
 +\frac12 z^2 C_\eps(f,f)+\mathcal O_K\left(|z|^3\norm{f}_\infty^3\right),
\end{equation}
uniformly whenever $|z|\norm{f}_\infty$ is sufficiently small.
\end{proposition}

\begin{proof}
By Lemma~\ref{lem:analytic}, the fixed-point equation
\eqref{eq:positivecluster} can be differentiated at zero. Since
$G_\eps(0)=G_\eps'(0)=0$ and $G_\eps''(0)=1$, its first two derivatives are
\begin{align}\label{eq:cluster-first-second-derivatives}
 Dv_{\eps,0}f&=\cA_\eps f,\nonumber\\
 D^2v_{\eps,0}(f,g)&=\cA_\eps[(\cA_\eps f)(\cA_\eps g)].
\end{align}
Differentiating \eqref{eq:Lambdaexact} under the integral is justified by
\eqref{eq:vbounds} and \eqref{eq:Bbounds}. Together with
\eqref{eq:superposition}, it gives the exact mean identity
\begin{equation}\label{eq:exact-mean-resolvent}
 \mu_\eps(f)=D\Lambda_\eps(0)f=a_\eps \B_\eps \cA_\eps f,
 \qquad f\in C_c^\infty(\Omega_T),
\end{equation}
which will also be used in Section~\ref{sec:bias}. The second derivative is
\begin{align*}
 D^2\Lambda_\eps(0)(f,g)
 &=a_\eps \B_\eps \cA_\eps[(\cA_\eps f)(\cA_\eps g)]\\
 &\quad+a_\eps^2\delta_\eps\int_{-\infty}^0\int_{\R^d}\lambda_0(x)
       (\cJ_\eps \cA_\eps f)(r,x)(\cJ_\eps \cA_\eps g)(r,x)\dd x\dd r.
\end{align*}
Differentiation of \eqref{eq:superposition} also identifies
$D^2\Lambda_\eps(0)(f,g)=n_\eps\Cov(\NN_\eps(f),\NN_\eps(g))$.
For the last integral, \eqref{eq:Jpast} gives
\begin{align*}
 &a_\eps^2\delta_\eps\left|\int_{-\infty}^0\int_{\R^d}\lambda_0(x)
       (\cJ_\eps \cA_\eps f)(r,x)(\cJ_\eps \cA_\eps g)(r,x)\dd x\dd r\right|\\
 &\quad\leq\delta_\eps\|\cA_\eps f\|_\infty \B_\eps|\cA_\eps g|
 \leq m\delta_\eps\|\lambda_0\|_\infty\|\cA_\eps f\|_\infty
                              \|\cA_\eps g\|_{L^\infty([0,T];L^1(\R^d))}.
\end{align*}
Lemma~\ref{lem:renewal} bounds the final product uniformly in $\eps$.
The same lemma gives convergence in both norms of
$w_\eps=(\cA_\eps f)(\cA_\eps g)$ to $w=(\cA f)(\cA g)$, and of
$\cA_\eps w_\eps$ to $\cA w$. Therefore \eqref{eq:Blimit} and
\eqref{eq:mean-duality} yield
\[
 a_\eps \B_\eps \cA_\eps w_\eps
 \longrightarrow m\int_{\R^d}\lambda_0(x)\cA w(0,x)\dd x
 =\int_0^T\int_{\R^d}\rho(t,x)\cA f(t,x)\cA g(t,x)\dd x\dd t.
\]
This proves all assertions in \eqref{eq:cumulants}.

For the third-order remainder, let $B_f:=\|f\|_\infty>0$,
$g:=f/B_f$, and $H_\eps(z):=\Lambda_\eps(zg)$. Lemma~\ref{lem:analytic}
and \eqref{eq:Bbounds} give $R>0$ and $M_K<\infty$ such that
$H_\eps$ is analytic on a neighborhood of $\{|z|\leq R\}$ and
$\sup_{|z|=R}|H_\eps(z)|\leq M_K$. These constants are independent of
$g$ and sufficiently small $\eps$; $R$ is also independent of $K$.
Cauchy's integral formula gives
\begin{equation}\label{eq:cauchy-coefficient-bound}
 \frac{H_\eps^{(j)}(0)}{j!}
 =\frac1{2\pi\mathrm{i}}\int_{|w|=R}\frac{H_\eps(w)}{w^{j+1}}\dd w,
 \qquad
 \left|\frac{H_\eps^{(j)}(0)}{j!}\right|\leq\frac{M_K}{R^j},
 \quad j\geq0.
\end{equation}
Using $H_\eps(0)=0$ and \eqref{eq:cumulants}, for $|z|B_f\leq R/2$,
\begin{align*}
 \left|\Lambda_\eps(zf)-z\mu_\eps(f)-\frac{z^2}{2}C_\eps(f,f)\right|
 \leq M_K\sum_{j=3}^\infty\left(\frac{|z|B_f}{R}\right)^j
 \leq\frac{2M_K}{R^3}|z|^3B_f^3.
\end{align*}
This is \eqref{eq:third}. For $B_f=0$ all terms vanish. Since
$M_K$ comes from spatial $L^1$ bounds, it depends on the support volume
and is independent of derivatives and oscillations of $f$.
\end{proof}

For the distribution-space tightness argument, the centered exponential estimates must depend on test function amplitudes and support volume without involving test function derivatives.

\begin{lemma}[Derivative-free centered exponential bounds]
\label{lem:centered-exponential}
Fix a compact spatial set $K$. There are $r>0$ and $0<C_K<\infty$ such that,
for every real bounded test function $f$ supported in $[0,T]\times K$, all sufficiently
small $\eps$, and every real $\theta$ with $|\theta|\|f\|_\infty\leq r$,
\begin{equation}\label{eq:centered-cgf-bound}
 0\leq\Lambda_\eps(\theta f)-\theta\mu_\eps(f)
       \leq C_K\theta^2\|f\|_\infty^2.
\end{equation}
Here $r$ and the small-$\eps$ cutoff can be chosen independently of $K$;
$C_K$ may depend on its volume but not on derivatives of $f$.
For $B_f=\|f\|_\infty>0$ and $x>0$, this implies
\begin{equation}\label{eq:bernstein-event}
 \Pp\left(\left|\left(\NN_\eps-\mu_\eps\right)(f)\right|>x\right)
 \leq 2\exp\left[-n_\eps\min\left\{
      \frac{x^2}{4C_KB_f^2},\frac{rx}{2B_f}\right\}\right].
\end{equation}
The inequality is trivial for $B_f=0$.
\end{lemma}

\begin{proof}
Let $X_\eps(f)$ be the unit-effective-population event observable, so
$\Lambda_\eps(\theta f)=\log\E\left[e^{\theta X_\eps(f)}\right]$ and
$\E[X_\eps(f)]=\mu_\eps(f)$ by \eqref{eq:cumulants}. Jensen's inequality gives
\[
 \Lambda_\eps(\theta f)-\theta\mu_\eps(f)
 =\log\E\left[e^{\theta(X_\eps(f)-\E X_\eps(f))}\right]\geq0.
\]
For $B_f=\|f\|_\infty>0$, apply the coefficient bound
\eqref{eq:cauchy-coefficient-bound} to $g=f/B_f$. For
$|\theta|B_f\leq R/2$,
\begin{align}\label{eq:centered-cgf-coefficient-sum}
 |\Lambda_\eps(\theta f)-\theta\mu_\eps(f)|
 &\leq M_K\sum_{j=2}^\infty
               \left(\frac{|\theta|B_f}{R}\right)^j
 \leq\frac{2M_K}{R^2}\theta^2B_f^2.
\end{align}
The construction of \eqref{eq:cluster-contraction} applies directly to
bounded measurable sources. Its dominating transform and the bounds in
\eqref{eq:Bbounds} therefore give
$\sup_{|z|=R}|\Lambda_\eps(zf/B_f)|\leq M_K$ for bounded measurable
$f$ as well. Cauchy's integral formula yields
\eqref{eq:cauchy-coefficient-bound} and its sum
\eqref{eq:centered-cgf-coefficient-sum} for every test function in
the lemma. Taking $r=R/2$ and $C_K=2M_K/R^2$ proves
\eqref{eq:centered-cgf-bound}. Specifically, choose $R=r_T/2$ in
\eqref{eq:cauchy-coefficient-bound} and $r=R/2=r_T/4$, where $r_T$
satisfies \eqref{eq:cluster-source-radius}; both radii are independent of $K$.

Let $D_\eps(f):=(\NN_\eps-\mu_\eps)(f)$. By \eqref{eq:superposition},
for $0<\theta\leq r/B_f$,
\begin{align*}
 \log\E\left[e^{\pm n_\eps\theta D_\eps(f)}\right]
 &=n_\eps\left[\Lambda_\eps(\pm\theta f)\mp\theta\mu_\eps(f)\right]
 \leq n_\eps C_K\theta^2B_f^2,\\
 \Pp\left(\pm D_\eps(f)>x\right)
 &\leq\exp\left\{-n_\eps\left(\theta x-C_K\theta^2B_f^2\right)\right\}.
\end{align*}
Choose $\theta=\min\left\{x/(2C_KB_f^2),r/B_f\right\}$. If the first value is
selected, $\theta x-C_K\theta^2B_f^2=x^2/(4C_KB_f^2)$.
Otherwise $x\geq2C_KrB_f$, and
\[
 \theta x-C_K\theta^2B_f^2=rx/B_f-C_Kr^2\geq rx/(2B_f).
\]
Adding the two one-sided bounds proves \eqref{eq:bernstein-event}.
The case $B_f=0$ has $D_\eps(f)=0$.
\end{proof}

\begin{remark}[Covariance of the prescribed past]\label{rem:initial}
In \eqref{eq:cumulants}, the prescribed-past contribution to
$n_\eps\Cov(\NN_\eps(f),\NN_\eps(g))$ is $\mathcal O(\delta_\eps)$ and
vanishes in the limit. The branching noise has variance of order
$n_\eps^{-1}$, while the initial equation-scale Poisson noise has variance
of order $\delta_\eps n_\eps^{-1}$. The limiting covariance is
\eqref{eq:cov}, with the zero initial condition specified by
\eqref{eq:Gaussian}.
\end{remark}

\section{Proofs of the CLT and MDP}\label{sec:clt-mdp}
This section proves Theorems~\ref{thm:clt} and \ref{thm:mdp} from the
cluster cumulants of Section~\ref{sec:analytic}.

\subsection{Finite-dimensional limits}\label{sec:finite:dim:limits}

The cumulant expansion gives Gaussian and moderate-deviation limits for
every finite family of test functions, including families with a singular
limiting covariance matrix.

\begin{lemma}[Finite-dimensional CLT and MDP]\label{lem:finite-dimensional-limits}
Under Assumption~\ref{ass:finite-mean}, let
$f_1,\ldots,f_k\in\cD(E_T)$ be real and let
$\Gamma:=(C(f_i,f_j))_{i,j=1}^k$. Then
\[
 (Z_\eps(f_1),\ldots,Z_\eps(f_k))
       \Longrightarrow\mathcal N_k(0,\Gamma).
\]
If $b_\eps\to\infty$ and $b_\eps/\sqrt{n_\eps}\to0$, the vectors
$X_\eps:=(Y_\eps(f_1),\ldots,Y_\eps(f_k))$ satisfy an LDP on $\R^k$
at speed $b_\eps^2$ with good rate function
\begin{equation}\label{eq:finite-dimensional-rate}
 I_\Gamma(x)=\sup_{\theta\in\R^k}
       \left\{\theta\cdot x-\theta^\top\Gamma\theta/2\right\}
 =\begin{cases}
   x^\top\Gamma^\dagger x/2,&x\in\operatorname{Ran}\Gamma,\\
   +\infty,&x\notin\operatorname{Ran}\Gamma,
  \end{cases}
\end{equation}
where $\Gamma^\dagger$ inverts the positive eigenvalues and is zero on
$\ker\Gamma$.
More precisely, for every open $O\subset\R^k$ and closed
$F\subset\R^k$,
\begin{align}\label{eq:finite-dimensional-mdp}
 \liminf_{\eps\downarrow0}b_\eps^{-2}\log\Pp(X_\eps\in O)
 &\geq-\inf_{x\in O}I_\Gamma(x),\nonumber\\
 \limsup_{\eps\downarrow0}b_\eps^{-2}\log\Pp(X_\eps\in F)
 &\leq-\inf_{x\in F}I_\Gamma(x).
\end{align}
\end{lemma}

\begin{proof}
Let
$f_\theta:=\sum_{i=1}^k\theta_if_i$ for $\theta\in\R^k$.
Let $\Gamma_\eps:=(C_\eps(f_i,f_j))_{i,j=1}^k$ and
$\Gamma:=(C(f_i,f_j))_{i,j=1}^k$.
Equations~\eqref{eq:superposition} and \eqref{eq:third} give
\begin{equation}\label{eq:finite-projection-cgf}
 \begin{aligned}
 \log\E\left[e^{\mathrm{i}\sum_{i=1}^k\theta_iZ_\eps(f_i)}\right]
 &=-\frac12\theta^\top\Gamma_\eps\theta+r_\eps(\theta),
 &|r_\eps(\theta)|&\leq C_K n_\eps^{-1/2}\|f_\theta\|_\infty^3,\\
 \lambda_\eps(\theta):=b_\eps^{-2}\log
                  \E \left[e^{b_\eps^2\sum_{i=1}^k\theta_iY_\eps(f_i)}\right]
 &=\frac12\theta^\top\Gamma_\eps\theta+s_\eps(\theta),
 &|s_\eps(\theta)|&\leq C_K\frac{b_\eps}{\sqrt{n_\eps}}
                                      \|f_\theta\|_\infty^3.
 \end{aligned}
\end{equation}
Here $[0,T]\times K$ contains the union of the test function supports,
and the remainder bounds in \eqref{eq:finite-projection-cgf} hold
for all sufficiently small $\eps$, uniformly on each compact
$\theta$-set. Proposition~\ref{prop:cumulants} gives
$\Gamma_\eps\to\Gamma$. Since $n_\eps\to\infty$ and
$b_\eps/\sqrt{n_\eps}\to0$, the characteristic functions converge to
$e^{-\theta^\top\Gamma\theta/2}$, proving the joint Gaussian limit.
The second line of \eqref{eq:finite-projection-cgf} gives the locally uniform limiting cumulant
$\lambda(\theta)=\theta^\top\Gamma\theta/2$.

The spectral decomposition of $\Gamma$ gives the Legendre transform
\eqref{eq:finite-dimensional-rate}: a nonzero component of $x$ in
$\ker\Gamma$ makes the supremum infinite, and maximizing on
$\operatorname{Ran}\Gamma$ gives $\theta=\Gamma^\dagger x$.
To prove the projected MDP \eqref{eq:finite-dimensional-mdp} with
rate function \eqref{eq:finite-dimensional-rate}, let $h_\eps:=b_\eps^2$.
If $x\in\operatorname{Ran}\Gamma$, let
$\theta:=\Gamma^\dagger x$ and change probability by
\begin{equation}\label{eq:finite-projection-tilt}
 \frac{\dd\Pp_{\eps,\theta}}{\dd\Pp}
 =e^{h_\eps[\theta\cdot X_\eps-\lambda_\eps(\theta)]},\qquad
 h_\eps^{-1}\log\E_{\eps,\theta}
                  \left[e^{h_\eps u\cdot(X_\eps-x)}\right]
 \longrightarrow\frac12u^\top\Gamma u.
\end{equation}
Indeed, the left side equals
$\lambda_\eps(\theta+u)-\lambda_\eps(\theta)-u\cdot x$,
and $\Gamma\theta=x$. For each coordinate, Chernoff's inequality
with $u=\pm a e_i$ gives an exponentially negative bound on
$|X_{\eps,i}-x_i|\geq r/\sqrt{k}$ by choosing $a>0$ with
$-ar/\sqrt{k}+a^2\Gamma_{ii}/2<0$. The union bound over $1\leq i\leq k$ gives
$\Pp_{\eps,\theta}(|X_\eps-x|<r)\to1$. Reversing
\eqref{eq:finite-projection-tilt} gives
\[
 \liminf_{\eps\downarrow0}h_\eps^{-1}\log
                        \Pp(|X_\eps-x|<r)
 \geq-\theta\cdot x+\lambda(\theta)-|\theta|r
 =-I_\Gamma(x)-|\theta|r.
\]
Taking arbitrarily small balls inside an open set proves its lower bound.

For the upper bound, any $a<I_\Gamma(x)$ admits a $\theta$ with
$\theta\cdot x-\lambda(\theta)>a$. A sufficiently small neighborhood
$U_x$ then satisfies
\[
 \limsup_{\eps\downarrow0}h_\eps^{-1}\log\Pp(X_\eps\in U_x)
 \leq-\inf_{y\in U_x}\theta\cdot y+\lambda(\theta)<-a.
\]
A finite cover proves the compact upper bound. Finally, for fixed $a>0$,
\[
 \limsup_{\eps\downarrow0}h_\eps^{-1}\log
                         \Pp(\|X_\eps\|_\infty>R)
 \leq-aR+\frac{a^2}{2}\max_{1\leq i\leq k}\Gamma_{ii}.
\]
This proves exponential tightness and extends the compact upper bound
to closed sets by mass outside a large cube. The sublevels of
$I_\Gamma$ are bounded closed ellipsoids in $\operatorname{Ran}\Gamma$,
so the projected MDP has a good rate function even when $\Gamma$ is singular.
\end{proof}

\subsection{Tightness in the weak distribution topology}

Uniform control of localized Fourier coordinates upgrades the finite-dimensional limits to tightness and exponential tightness in the weak distribution topology.

\begin{lemma}[Local distribution tightness]\label{lem:distribution-tightness}
The laws of $Z_\eps$ are tight in the weak distribution topology, and the laws of $Y_\eps$ are exponentially tight at speed $b_\eps^2$ in that topology. Both conclusions follow from compact sets in suitable local negative Sobolev spaces.
\end{lemma}

\begin{proof}
Let $D=d+1$, and choose cutoffs $\chi_j\in C_c^\infty(E_T)$, $j\geq1$,
equal to one on a compact exhaustion, with $0\leq\chi_j\leq1$ and
$\chi_{j+1}=1$ on $\supp\chi_j$.
For each $j$, fix a rectangular torus $\mathbb T_j$ whose fundamental
domain is
\[
 Q_j:=\prod_{\ell=0}^d\left(a_{j,\ell},a_{j,\ell}+L_{j,\ell}\right)\subset E_T,
 \qquad L_{j,\ell}>0,\qquad\supp\chi_j\Subset Q_j,
 \qquad V_j:=\prod_{\ell=0}^dL_{j,\ell}.
\]
Such a rectangle exists because $\supp\chi_j$ is compact in $E_T$.
Let $\Lambda_j:=\prod_{\ell=0}^d(L_{j,\ell}\mathbb Z)$ and
$\mathbb T_j:=\R^D/\Lambda_j$.
For $\varphi\in C^\infty(\mathbb T_j)$, write $\varphi_{\rm per}$ for
its smooth periodic lift to $\R^D$ and define
\[
 \langle\chi_jz,\varphi\rangle_{\mathbb T_j}
       :=z(\chi_j\varphi_{\rm per}).
\]
The product $\chi_j\varphi_{\rm per}$ is a test function in $\cD(E_T)$.
This formula first extends the compactly supported distribution $\chi_jz$
by zero to $\R^D$ and then identifies opposite faces of $Q_j$;
the support lies away from the faces.
All local Sobolev norms below refer to this distribution on $\mathbb T_j$.

We use the equivalent Fourier
norm $\|u\|_{H^{-s}}^2:=\sum_{k\in\mathbb Z^D}
(1+|k|^2)^{-s}|\widehat u(k)|^2$ throughout this proof.
For $y=(t,x)$ and $k=(k_0,\ldots,k_d)\in\mathbb Z^D$, let
\[
 \vartheta_{j,k}(y):=2\pi\sum_{\ell=0}^d
            \frac{k_\ell(y_\ell-a_{j,\ell})}{L_{j,\ell}},\qquad
 f_{j,k}^{(1)}:=\chi_j\cos\vartheta_{j,k},\quad
 f_{j,k}^{(2)}:=\chi_j\sin\vartheta_{j,k}.
\]
These real functions belong to $\cD(E_T)$, have common support contained
in $\supp\chi_j$, and satisfy $\left\|f_{j,k}^{(r)}\right\|_\infty\leq1$.
For real $z$, the Fourier coefficient, with the orthonormal normalization
on $Q_j$, is
\begin{equation}\label{eq:localized-fourier-coefficient}
 \widehat{\chi_jz}(k)
 :=V_j^{-1/2}z\left(\chi_je^{-\mathrm{i}\vartheta_{j,k}}\right)
 =V_j^{-1/2}\left[z\left(f_{j,k}^{(1)}\right)-\mathrm{i}z\left(f_{j,k}^{(2)}\right)\right].
\end{equation}
Taking the squared modulus in \eqref{eq:localized-fourier-coefficient}
and summing with weights $(1+|k|^2)^{-s}$ gives
\begin{equation}\label{eq:localized-fourier-norm}
 \left\|\chi_j z\right\|_{H^{-s}}^2
 \leq c_j\sum_{k\in\mathbb Z^D}\sum_{r=1}^2
             \left(1+|k|^2\right)^{-s}\left|z\left(f_{j,k}^{(r)}\right)\right|^2.
\end{equation}
With the chosen equivalent Fourier norm one may take $c_j=V_j^{-1}$,
and \eqref{eq:localized-fourier-norm} is an equality.
The usual Sobolev weights are
$(1+4\pi^2\sum_{\ell=0}^d k_\ell^2/L_{j,\ell}^2)^{-s}$;
they are comparable to $(1+|k|^2)^{-s}$ because all side lengths are
fixed and positive for each $j$.
By \eqref{eq:cumulants} and \eqref{eq:cauchy-coefficient-bound},
\[
 \sup_{0<\eps\leq\eps_0}\sup_{k\in\mathbb Z^D}\max_{1\leq r\leq2}
 \E\left[\left|Z_\eps\left(f_{j,k}^{(r)}\right)\right|^2\right]\leq C_j.
\]
For $s_0>D/2$, \eqref{eq:localized-fourier-norm} therefore gives
\[
 \sup_{0<\eps\leq\eps_0}\E\left[\left\|\chi_j Z_\eps\right\|_{H^{-s_0}}^2\right]
 \leq2c_jC_j\sum_{k\in\mathbb Z^D}\left(1+|k|^2\right)^{-s_0}=:B_j<\infty.
\]
Choose $R_j$ so that $B_j/R_j^2\leq\eta2^{-j}$. Markov's inequality
and the union bound yield
\begin{equation}\label{eq:distribution-tightness-probability}
 \Pp\left(\exists j\geq1:\left\|\chi_j Z_\eps\right\|_{H^{-s_0}}>R_j\right)
 \leq\sum_{j=1}^\infty B_j/R_j^2\leq\eta,\qquad 0<\eps\leq\eps_0.
\end{equation}
For a sequence $\mathbf R:=(R_j)_{j\geq1}$ of finite radii, define
\begin{equation}\label{eq:distribution-compact-sets}
 \mathsf K_{s_0}(\mathbf R)
 :=\bigcap_{j=1}^\infty
   \left\{z\in\cD'(E_T):\left\|\chi_jz\right\|_{H^{-s_0}(\mathbb T_j)}\leq R_j\right\}.
\end{equation}
We prove compactness in the local $H^{-s}$ topology for $s>s_0$.
Let $(z_n)_{n\geq1}\subset\mathsf K_{s_0}(\mathbf R)$, and let
$\widehat z_{n,j}(k)$ denote the Fourier coefficients of $\chi_jz_n$.
Then
\begin{equation}\label{eq:localized-fourier-tail}
 \begin{aligned}
 \left|\widehat z_{n,j}(k)\right|&\leq R_j\left(1+|k|^2\right)^{s_0/2},\\
 \sum_{k\in\mathbb Z^D,\,\,|k|>N}
       \left(1+|k|^2\right)^{-s}\left|\widehat z_{n,j}(k)\right|^2
 &\leq\left(1+N^2\right)^{-(s-s_0)}R_j^2.
 \end{aligned}
\end{equation}
Enumerate the countable set of pairs $(j,k)$. Successive subsequence
selection for each bounded coefficient, followed by taking the diagonal,
gives one subsequence $(z_{n_\ell})_{\ell\geq1}$ such that
$\widehat z_{n_\ell,j}(k)\to a_j(k)$ for every $j,k$.
For every $N$, convergence of the finite sum gives
\[
 \sum_{k\in\mathbb Z^D,\,\,|k|\leq N}
       \left(1+|k|^2\right)^{-s_0}|a_j(k)|^2\leq R_j^2.
\]
Letting $N\to\infty$ defines $u_j\in H^{-s_0}(\mathbb T_j)$ with
Fourier coefficients $a_j(k)$ and norm at most $R_j$.
Splitting at frequency $N$ and applying
\eqref{eq:localized-fourier-tail} to both $\chi_jz_{n_\ell}$ and $u_j$ gives
\begin{align*}
 \left\|\chi_jz_{n_\ell}-u_j\right\|_{H^{-s}}^2
 \leq\sum_{k\in\mathbb Z^D,\,\,|k|\leq N}
       \left(1+|k|^2\right)^{-s}\left|\widehat z_{n_\ell,j}(k)-a_j(k)\right|^2
       +4R_j^2\left(1+N^2\right)^{-(s-s_0)}\longrightarrow0,
\end{align*}
where the limit is taken first in $\ell$ and then in $N$.
The identities $\chi_jz_n=\chi_j(\chi_{j+1}z_n)$ pass to distributional
limits, so $u_j=\chi_ju_{j+1}$. For a test function $f$ supported where
$\chi_j=1$, define $z(f):=u_j(f)$. These compatibility identities
make the definition independent of $j$ and imply $\chi_jz=u_j$.
Thus $z\in\mathsf K_{s_0}(\mathbf R)$ and $z_{n_\ell}\to z$ for the metric
\begin{equation}\label{eq:local-sobolev-metric}
 d_s(z,w):=\sum_{j=1}^\infty2^{-j}
               \left(\left\|\chi_j(z-w)\right\|_{H^{-s}(\mathbb T_j)}\wedge1\right).
\end{equation}
Sequential compactness in this metric proves compactness.
Moreover, for $\supp f\subset\{\chi_j=1\}$,
\[
 |(z-w)(f)|\leq\|\chi_j(z-w)\|_{H^{-s}}\|f\|_{H^s}.
\]
Hence the inclusion into the weak distribution topology is continuous,
and \eqref{eq:distribution-compact-sets} is compact in that topology.
The probability bound \eqref{eq:distribution-tightness-probability}
 gives
$\Pp(Z_\eps\notin\mathsf K_{s_0}(\mathbf R))\leq\eta$, proving tightness.

For exponential tightness, let $A>0$. Applying
\eqref{eq:centered-cgf-bound} and \eqref{eq:superposition} with
$\theta=b_\eps/\sqrt{n_\eps}$ gives, simultaneously for every $j,k,r$,
\[
 \log\E\left[e^{\pm b_\eps^2Y_\eps(f_{j,k}^{(r)})}\right]
 \leq b_\eps^2C_j,
\]
once $\eps$ is sufficiently small. The cutoff is common to all indices
because the source radius in Lemma~\ref{lem:centered-exponential}
is independent of the spatial support. Chernoff's inequality gives
\begin{equation}\label{eq:fourier-coordinate-tail}
 \Pp\left(\left|Y_\eps\left(f_{j,k}^{(r)}\right)\right|>R_j(1+|k|)\right)
 \leq2\exp\left\{-b_\eps^2\left[R_j(1+|k|)-C_j\right]\right\}.
\end{equation}
Let $c_D:=\sum_{k\in\mathbb Z^D}\sum_{r=1}^2e^{-|k|}<\infty$.
For $b_\eps^2\geq1$ and $R_j\geq1$,
\[
 \sum_{k\in\mathbb Z^D}\sum_{r=1}^2
 \Pp\left(\left|Y_\eps\left(f_{j,k}^{(r)}\right)\right|>R_j(1+|k|)\right)
 \leq2c_D e^{-b_\eps^2(R_j-C_j)}.
\]
Choose $R_j\geq C_j+A+j+\log(4c_D)+1$. Summing over $j$ gives
\[
 \Pp\left(\exists j\geq1,k\in\mathbb Z^D,r\in\{1,2\}:
                 \left|Y_\eps\left(f_{j,k}^{(r)}\right)\right|>R_j(1+|k|)\right)
 \leq e^{-b_\eps^2A}.
\]
For $s_1>D/2+1$, define
\begin{equation}\label{eq:distribution-exponential-compact}
 \begin{aligned}
 \widetilde R_j^2
 &:=2c_jR_j^2\sum_{k\in\mathbb Z^D}
                     \left(1+|k|^2\right)^{-s_1}(1+|k|)^2<\infty,\\
 \mathsf K_A^{\rm dist}&:=\mathsf K_{s_1}(\widetilde{\mathbf R}).
 \end{aligned}
\end{equation}
On the complementary event, \eqref{eq:localized-fourier-norm} gives
$\|\chi_jY_\eps\|_{H^{-s_1}}\leq\widetilde R_j$ for every $j$.
The compactness proof for \eqref{eq:distribution-compact-sets}, with
$s_0=s_1$ and $s>s_1$, therefore gives a weakly compact set satisfying
\[
 \Pp\left(Y_\eps\notin\mathsf K_A^{\rm dist}\right)\leq e^{-b_\eps^2A}.
\]
This proves exponential tightness at speed $b_\eps^2$.
\end{proof}

\begin{proof}[Completion of Theorems~\ref{thm:clt} and \ref{thm:mdp}]
Let $\cTG f:=\sqrt\rho\,\cA f$. By \eqref{eq:normalized-mean-density-bound}
and \eqref{eq:resolvent},
\[
 \left\|\cTG f\right\|_{L^2(E_T)}^2
 \leq T\|\lambda_0\|_\infty\|\cA f\|_\infty
                                  \|\cA f\|_{L^\infty([0,T];L^1(\R^d))}<\infty.
\]
It\^o's isometry for the white-noise integral in \eqref{eq:Gaussian} gives
\[
\E[Z(f)Z(g)]=\left\langle\cTG f,\cTG g\right\rangle_{L^2(E_T)}=C(f,g).
\]
The finite-dimensional CLT in Lemma~\ref{lem:finite-dimensional-limits}
and the tightness in Lemma~\ref{lem:distribution-tightness} identify every subsequential
limit with this Gaussian field.
On each compact set \eqref{eq:distribution-compact-sets}, including
\eqref{eq:distribution-exponential-compact}, the map
$z\mapsto\left(z\left(f_{j,k}^{(r)}\right)\right)_{j,k,r}$ is continuous and injective:
its coordinates determine each $\chi_jz$ by Fourier expansion.
A continuous injection from a compact space into the Hausdorff
countable product of coordinate spaces is a homeomorphism onto its
image. Thus these coordinates determine both the relative topology
and the Borel probability laws on these compact sets.
This proves \eqref{eq:clt}.

We give the extension from finite-dimensional projections to the distribution space for the MDP explicitly.
For an open $O\subset\cD'(E_T)$ and $z\in O$, a cylinder
\[
 U=\left\{y:|y(f_i)-z(f_i)|<\eta,\ 1\leq i\leq\ell\right\}\subset O
\]
exists by the definition of the weak distribution topology. The
finite-dimensional lower bound of Lemma~\ref{lem:finite-dimensional-limits} gives
\[
 \liminf_{\eps\downarrow0}b_\eps^{-2}\log\Pp(Y_\eps\in O)
 \geq-I_{f_1,\ldots,f_\ell}\left(z(f_1),\ldots,z(f_\ell)\right)
 \geq-I_G(z),
\]
since its supremum ranges only over $\operatorname{span}\{f_1,\ldots,f_\ell\}$.
For a compact set $K$ and $a<\inf_{y\in K} I_G(y)$, each $z\in K$ admits a test
$f$ with $z(f)-C(f,f)/2>a$. Choose $\eta>0$ so that
$z(f)-\eta-C(f,f)/2>a$. The cumulant limit and Chernoff's inequality give
\[
 \limsup_{\eps\downarrow0}b_\eps^{-2}
 \log\Pp\left(Y_\eps(f)>z(f)-\eta\right)
 \leq-z(f)+\eta+C(f,f)/2<-a.
\]
A finite cover of $K$ proves its upper bound. If $F$ is closed and
$K_A$ is a compact set provided by exponential tightness, then
\[
 \limsup_{\eps\downarrow0}b_\eps^{-2}\log\Pp(Y_\eps\in F)
 \leq\max\left\{-\inf_{y\in F\cap K_A}I_G(y),-A\right\}
 \leq\max\left\{-\inf_{y\in F}I_G(y),-A\right\}.
\]
Letting $A\uparrow\infty$ proves the full upper bound with the first
expression in \eqref{eq:MDPrate}.

Suppose $I_G(z)<\infty$. Applying its defining supremum to $af$,
$a\in\R$, gives
\begin{equation}\label{eq:Gaussian-residual-functional-bound}
 |z(f)|^2\leq2I_G(z)\left\|\cTG f\right\|_2^2.
\end{equation}
Thus $z(f)=0$ whenever $\cTG f=0$, and the functional
$\cTG f\mapsto z(f)$ extends continuously to the Hilbert space
$\mathsf H=\overline{\cTG\cD(E_T)}^{L^2}$.
The Riesz representation theorem gives a unique $\zeta_0\in\mathsf H$
with $z(f)=\ip{\cTG f}{\zeta_0}_2$. Completing the square gives
\[
 z(f)-\frac12C(f,f)
 =\frac12\|\zeta_0\|_2^2
       -\frac12\left\|\cTG f-\zeta_0\right\|_2^2.
\]
Density of $\cTG\cD(E_T)$ in $\mathsf H$ proves
$I_G(z)=\|\zeta_0\|_2^2/2$. Every other representing $\zeta$ has
$\zeta-\zeta_0\perp\mathsf H$ and
$\|\zeta\|_2^2=\|\zeta_0\|_2^2+\|\zeta-\zeta_0\|_2^2$.
This proves the minimum-energy formula in \eqref{eq:MDPrate}.

Define $(\cS\zeta)(f):=\ip{\cTG f}{\zeta}_2$. For each fixed
$f\in\cD(E_T)$, the scalar evaluation map is
\[
 \zeta\longmapsto(\cS\zeta)(f)
 =\int_0^T\int_{\R^d}\sqrt{\rho(t,x)}\,\cA f(t,x)\zeta(t,x)\dd x\dd t.
\]
Since $\cTG f\in L^2(E_T)$, weak convergence $\zeta_\iota\rightharpoonup\zeta$
implies $(\cS\zeta_\iota)(f)\to(\cS\zeta)(f)$ by the definition of the
weak $L^2$ topology. Thus $\cS$ is continuous from weak $L^2(E_T)$ to
$\sigma(\cD',\cD)$, and
\[
 \{z:I_G(z)\leq M\}
 =\cS\left\{\zeta\in L^2(E_T):\|\zeta\|_2^2\leq2M\right\}.
\]
The Hilbert ball is weakly compact, so this equality proves goodness.
Finally, \eqref{eq:resolvent} gives $\cA (-\LL_\kappa\phi)=\phi$ for
$\phi\in\cD(E_T)$. Consequently,
\[
 -Z(\LL_\kappa\phi)=\int_0^T\int_{\R^d}\sqrt\rho\,\phi\,W(\dd t,\dd x),
 \qquad
 -z(\LL_\kappa\phi)=\int_0^T\int_{\R^d}\sqrt\rho\,\phi\,\zeta\dd x\dd t.
\]
Note that the field at time $t$ depends only on the
forcing during $[0,t]$. More explicitly, for
$\psi\in C_c^\infty(\R^d)$ the forward mild representations are
\begin{equation}\label{eq:Gaussian-forward-causal}
 \begin{aligned}
 Z_t(\psi)&=\frac1m\int_0^t\int_{\R^d}
       \sqrt{\rho(s,x)}\,\cQ_{t-s}\psi(x)W(\dd s,\dd x),\\
 z_t(\psi)&=\frac1m\int_0^t\int_{\R^d}
       \sqrt{\rho(s,x)}\,\cQ_{t-s}\psi(x)\zeta(s,x)\dd x\dd s.
 \end{aligned}
\end{equation}
The $L^2$ contraction of $\cQ_t$ and
\eqref{eq:normalized-mean-density-bound} give
\begin{align*}
 \E\left[\left|Z_t(\psi)\right|^2\right]
 &\leq\frac{t\|\lambda_0\|_\infty}{m^2}\|\psi\|_2^2,\\
 \left|z_t(\psi)\right|^2
 &\leq\frac{t\|\lambda_0\|_\infty}{m^2}\|\psi\|_2^2
                  \int_0^t\int_{\R^d}|\zeta(s,x)|^2\dd x\dd s.
\end{align*}
These estimates justify the stochastic Fubini theorem and Fubini's
theorem when \eqref{eq:Gaussian-forward-causal} is paired with a
space--time test function; \eqref{eq:resolvent} then yields
\eqref{eq:Gaussian} and \eqref{eq:MDPrate}, respectively.
They also give $Z_t(\psi)\to0$ in $L^2(\Pp)$ and $z_t(\psi)\to0$
as $t\downarrow0$. This proves \eqref{eq:SPDE} and
\eqref{eq:MDPcontrol} with the stated zero-initial mild convention.
The proof is complete.
\end{proof}

\section{Proofs of the LDP}\label{sec:spatialLDP}

We first prove Theorem~\ref{thm:full-spatial-finite} under
Assumption~\ref{ass:finite-mass}.
Section~\ref{subsec:infinite-bounds} gives compact-kernel bounds
for initial profiles of possibly infinite mass.
The proof uses the following order: the variational representation
and mass estimate, approximation by smooth controls, the controlled
law of large numbers and likelihood estimate, and finally
the LDP bounds.

\subsection{The variational action and its domain}

The variational action in \eqref{eq:spatialJ} represents the weak PDE residual by an $L^2(\nu)$ control; its quadratic expansion then yields the Gaussian action in the next proposition.

\begin{proposition}[Variational interpretation and its Gaussian tangent]\label{prop:action}
Whenever $J_{\rm sp}(\nu)<\infty$, there exists $v\in L^2(\nu)$ such that, for every $\phi\in C_c^\infty(\Omega_T)$,
\begin{equation}\label{eq:weakcontrolled}
 \begin{aligned}
 -\ip\nu{\LL_\kappa\phi}-m\int_{\R^d}\lambda_0(x)\phi(0,x)\dd x
 &=\int_{\Omega_T}v(t,x)\phi(t,x)\nu(\dd t,\dd x),\\
 J_{\rm sp}(\nu)&=\frac12\int_{\Omega_T}(v(t,x))^2\nu(\dd t,\dd x).
 \end{aligned}
\end{equation}
The representing $v$ is unique in the closure of the test space in $L^2(\nu)$. For a smooth positive density $q$ with the prescribed initial value, this yields
\begin{equation}\label{eq:formalaction}
 J_{\rm sp}(q\dd t\dd x)=\frac12\int_0^T\int_{\R^d}
 \frac{[m\partial_tq-\B q+\kappa q]^2}{q}\dd x\dd t,
\end{equation}
provided the weak identity has no additional boundary term and the integral in \eqref{eq:formalaction} is finite.
For a smooth compactly supported perturbation $z$ with $z_0=0$,
supported where $\rho$ is bounded below, define
\begin{equation}\label{eq:Gaussian-tangent-perturbation}
 q_s:=\rho+s z.
\end{equation}
Its second variation is
\begin{equation}\label{eq:Gaussian-tangent}
 \lim_{s\to0}s^{-2}J_{\rm sp}((\rho+s z)\dd t\dd x)
 =\frac12\int_0^T\int_{\R^d}\frac{[m\partial_tz-\B z+\kappa z]^2}{\rho}\dd x\dd t.
\end{equation}
\end{proposition}

\begin{proof}
Let
\[
 \ell_\nu(\phi):=-\nu(\LL_\kappa\phi)
                 -m\int_{\R^d}\lambda_0(x)\phi(0,x)\dd x,
 \qquad \phi\in C_c^\infty(\Omega_T).
\]
For every $a\in\R$, \eqref{eq:spatialJ} gives
$a\ell_\nu(\phi)-a^2\|\phi\|_{L^2(\nu)}^2/2\leq J_{\rm sp}(\nu)$.
Maximizing in $a$ yields
\begin{equation}\label{eq:spatial-residual-functional-bound}
 |\ell_\nu(\phi)|\leq\sqrt{2J_{\rm sp}(\nu)}\|\phi\|_{L^2(\nu)}.
\end{equation}
If $\|\phi\|_{L^2(\nu)}=0$, the same inequality for all $a$ forces
$\ell_\nu(\phi)=0$. Thus the Riesz representation theorem on the
closure of the test space gives $v\in L^2(\nu)$ with
$\ell_\nu(\phi)=\int_{\Omega_T}v\phi\dd\nu$. Completing the square gives
\[
 \ell_\nu(\phi)-\frac12\int_{\Omega_T}\phi^2\dd\nu
 =\frac12\int_{\Omega_T}v^2\dd\nu
       -\frac12\int_{\Omega_T}(v-\phi)^2\dd\nu.
\]
Approximation of $v$ in the closed test space proves
\eqref{eq:weakcontrolled}, including uniqueness there.
In fact $C_c^\infty(\Omega_T)$ is dense in $L^2(\nu)$.
To verify this, fix $u\in L^2(\nu)$ and a compact exhaustion
$K_n\uparrow\Omega_T$. For
$u_n=\one_{K_n}((-n)\vee(u\wedge n))$,
\begin{equation}\label{eq:L2-Radon-truncation}
 \|u-u_n\|_{L^2(\nu)}^2
 \leq\int_{\Omega_T\setminus K_n}|u|^2\dd\nu
       +\int_{\{|u|>n\}}|u|^2\dd\nu\longrightarrow0.
\end{equation}
Since $\nu(K_n)<\infty$, quantizing the bounded values of $u_n$
gives simple functions $s=\sum_{i=1}^N c_i\one_{A_i}$,
$A_i\subset K_n$, with arbitrarily small $L^2(\nu)$ error.
Radon regularity provides compact $F_i\subset A_i$ and relatively
compact open $O_i\supset A_i$ with $\nu(O_i\setminus F_i)$
arbitrarily small. Choose $g_i\in C_c(\Omega_T)$ with
$0\leq g_i\leq1$, $g_i=1$ on $F_i$, and $g_i=0$ outside $O_i$.
Such a function can be formed from distances to $F_i$ and to the
complement of $O_i$. Then
\begin{align}\label{eq:L2-Radon-continuous-approximation}
 \|\one_{A_i}-g_i\|_{L^2(\nu)}^2
 &\leq\nu(O_i\setminus F_i),\nonumber\\
 \left\|s-\sum_{i=1}^N c_i g_i\right\|_{L^2(\nu)}
 &\leq\sum_{i=1}^N|c_i|\left\{\nu(O_i\setminus F_i)\right\}^{1/2}.
\end{align}
Finally, for $g\in C_c(\Omega_T)$, extend $g$ continuously to negative
times by reflection at zero and by zero beyond its compact temporal
support. Convolution of this compactly supported continuous extension
with a smooth approximate identity gives
$g_\gamma\in C_c^\infty(\Omega_T)$ with
$\|g_\gamma-g\|_\infty\to0$. For all sufficiently small $\gamma$,
their supports lie in a fixed compact $K'\subset\Omega_T$, and hence
\begin{equation}\label{eq:L2-Radon-smooth-approximation}
 \|g_\gamma-g\|_{L^2(\nu)}
 \leq\left\{\nu(K')\right\}^{1/2}\|g_\gamma-g\|_\infty\longrightarrow0.
\end{equation}
Combining \eqref{eq:L2-Radon-truncation}--\eqref{eq:L2-Radon-smooth-approximation}
proves density, including at the boundary $t=0$, and therefore
uniqueness of the representing control $\nu$-almost everywhere.

For a smooth positive $q$ with $q(0)=\lambda_0$, choose $t_*<T$
and $R>0$ so that $\phi=0$ on $[t_*,T)\times\R^d$ and its spatial
support lies strictly inside $B_R$. Integration by parts on
$[0,t_*]\times B_R$ gives, with $n_R$ the outward unit normal,
\begin{align}\label{eq:spatial-action-integration-by-parts}
 -\int_0^{t_*}\int_{B_R}q\LL_\kappa\phi\dd x\dd t
 &=-m\left[\int_{B_R}q(t,x)\phi(t,x)\dd x\right]_{t=0}^{t=t_*}
       \nonumber\\
 &\qquad+\int_0^{t_*}\int_{B_R}
          (m\partial_tq-\B q+\kappa q)\phi\dd x\dd t\nonumber\\
 &\qquad\qquad+\frac12\int_0^{t_*}\int_{\partial B_R}
       (\phi\Sigma\nabla q-q\Sigma\nabla\phi)\cdot n_R
                                      \dd S\dd t.
\end{align}
The terminal and spatial boundary terms in
\eqref{eq:spatial-action-integration-by-parts} satisfy
\begin{equation}\label{eq:spatial-action-vanishing-boundaries}
 \begin{aligned}
 m\int_{B_R}q(t_*,x)\phi(t_*,x)\dd x&=0,\\
 \frac12\int_0^{t_*}\int_{\partial B_R}
       (\phi\Sigma\nabla q-q\Sigma\nabla\phi)\cdot n_R
                                      \dd S\dd t&=0,
 \end{aligned}
\end{equation}
because $\phi(t_*,\cdot)=0$ and both $\phi$ and $\nabla\phi$
vanish on $\partial B_R$. Thus
\eqref{eq:spatial-action-integration-by-parts} becomes
\begin{align}\label{eq:spatial-action-residual-pairing}
 -\int_0^T\int_{\R^d}q\LL_\kappa\phi\dd x\dd t
 &=m\int_{\R^d}q(0,x)\phi(0,x)\dd x\nonumber\\
 &\qquad\qquad+\int_0^T\int_{\R^d}
                (m\partial_tq-\B q+\kappa q)\phi\dd x\dd t.
\end{align}
The initial term in \eqref{eq:spatial-action-residual-pairing}
cancels with $m\int_{\R^d}\lambda_0\phi(0)\dd x$ in $\ell_\nu$. Hence
$v=(m\partial_tq-\B q+\kappa q)/q$, and
\eqref{eq:weakcontrolled} gives \eqref{eq:formalaction} when that energy
is finite. If the integral in \eqref{eq:formalaction} is infinite,
$J_{\rm sp}(q\dd t\dd x)<\infty$ would instead give an $L^2(q\dd t\dd x)$
control by \eqref{eq:weakcontrolled}. Its pairing with every compact
test function forces it to equal $(m\partial_tq-\B q+\kappa q)/q$
almost everywhere, a contradiction. Thus \eqref{eq:formalaction}
 also holds
with value $+\infty$.

For the perturbation $q_s$ in \eqref{eq:Gaussian-tangent-perturbation},
 let
$r_z:=m\partial_tz-\B z+\kappa z$ and choose a compact set containing
the supports of $z$ and $r_z$, on which $\rho\geq c>0$.
For $|s|\|z\|_\infty\leq c/2$, $\rho+sz\geq c/2$ on that set and
equals $\rho$ elsewhere. The PDE in \eqref{eq:operators} gives
\[
 m\partial_t(\rho+sz)-\B(\rho+sz)+\kappa(\rho+sz)=sr_z.
\]
The weak equation for $\rho$ in \eqref{eq:operators}, including its
trace $\lambda_0$, therefore gives residual $sr_z$ for $\rho+sz$.
Its control is $sr_z/(\rho+sz)$ on the compact support of $r_z$ and
zero elsewhere. Equation~\eqref{eq:weakcontrolled} yields
\[
 s^{-2}J_{\rm sp}((\rho+sz)\dd t\dd x)
 =\frac12\int_0^T\int_{\R^d}\frac{(r_z(t,x))^2}{\rho(t,x)+sz(t,x)}\dd x\dd t.
\]
Its integrand is bounded by $2r_z^2/c$ on a compact set and is zero
outside that set. Dominated convergence theorem proves \eqref{eq:Gaussian-tangent}.
For $\zeta=r_z/\sqrt\rho$, this energy is the control energy in
\eqref{eq:MDPcontrol} on the class of smooth compactly supported perturbations.
This completes the proof.
\end{proof}

\subsection{Global generator convergence}

For a bounded real test function $\phi$, set
\begin{equation}\label{eq:nonlinear-exponential-coefficient}
 c_\eps(\phi):=\frac{e^{\delta_\eps\phi}-1-\delta_\eps\phi}{\delta_\eps^2}.
\end{equation}

The uniform generator estimate in the following lemma bounds the error in the
microscopic upper bound by the total mass of the event measure.

\begin{lemma}[Uniform generator convergence]
\label{lem:global-generator}
Under Assumption~\ref{ass:finite-mean}, for every
$\phi\in C_c^\infty([0,T)\times\R^d)$,
\begin{equation}\label{eq:global-generator}
 \left\|\frac{\cK_\eps\phi-\phi}{\delta_\eps}
            -m\partial_t\phi-\B\phi\right\|_\infty\longrightarrow0.
\end{equation}
Consequently, with $\LL_\eps$ from \eqref{eq:discrete-generator},
\begin{equation}\label{eq:global-martingale-integrand}
 \left\|-\LL_\eps\phi-a_\eps \cK_\eps c_\eps(\phi)
                +\LL_\kappa\phi+\phi^2/2\right\|_\infty\longrightarrow0.
\end{equation}
The norms in \eqref{eq:global-generator} and
\eqref{eq:global-martingale-integrand} are taken over
$[0,\infty)\times\R^d$.
\end{lemma}
\begin{proof}
Extend $\phi$ by zero after $T$; it remains smooth because it vanishes
near $T$. Split the increment in \eqref{eq:K} into temporal and spatial
increments. For the temporal part, the fundamental theorem of calculus gives
\[
 R_{\eps,t}(t,x)=\E\left[S\int_0^1
       \partial_t\phi\left(t+u\delta_\eps S,x+\eps\Xi\right)\dd u\right]
       -m\partial_t\phi(t,x).
\]
If $\omega_t$ is a bounded modulus of uniform continuity of
$\partial_t\phi$, then
\begin{equation}\label{eq:global-generator-time-error}
 \|R_{\eps,t}\|_\infty
 \leq\E[S\omega_t(\delta_\eps S+\eps|\Xi|)]\longrightarrow0.
\end{equation}
Domination is by $2\|\partial_t\phi\|_\infty S$, integrable under
\eqref{eq:basic}. The spatial Taylor formula and $\E[\Xi]=0$ give
\begin{align*}
 \frac{\E[\phi(t,x+\eps\Xi)-\phi(t,x)]}{\eps^2}-\B\phi(t,x)
 &=\E\left[\int_0^1(1-u)\Xi^\top
      \left[D_x^2\phi(t,x+u\eps\Xi)-D_x^2\phi(t,x)\right]\Xi\dd u\right].
\end{align*}
With a bounded modulus $\omega_x$ for $D_x^2\phi$, this is bounded in
uniform norm by
\begin{equation}\label{eq:global-generator-space-error}
 \frac12\E\left[|\Xi|^2\omega_x(\eps|\Xi|)\right]\longrightarrow0.
\end{equation}
The integrable dominating variable is
$\|D_x^2\phi\|_\infty|\Xi|^2$. Since $\delta_\eps=\eps^2$,
adding \eqref{eq:global-generator-time-error} and
\eqref{eq:global-generator-space-error} proves \eqref{eq:global-generator}.

Let $r_\eps:=(\cK_\eps\phi-\phi)/\delta_\eps-m\partial_t\phi-\B\phi$.
By \eqref{eq:discrete-generator} and \eqref{eq:scales},
$\LL_\eps\phi=(\cK_\eps\phi-\phi)/\delta_\eps-\kappa_\eps \cK_\eps\phi$.
Also Taylor's integral remainder in
\eqref{eq:nonlinear-exponential-coefficient} gives
\[
 \left\|c_\eps(\phi)-\phi^2/2\right\|_\infty
 \leq\frac{\delta_\eps}{6}e^{\delta_\eps\|\phi\|_\infty}\|\phi\|_\infty^3.
\]
Uniform continuity and boundedness give $\cK_\eps\phi\to\phi$ and
$\cK_\eps(\phi^2)\to\phi^2$ uniformly on the full domain. Hence,
\begin{align*}
 &\left\|-\LL_\eps\phi-a_\eps \cK_\eps c_\eps(\phi)
                         +\LL_\kappa\phi+\phi^2/2\right\|_\infty\\
 &\quad\leq\|r_\eps\|_\infty
   +|\kappa_\eps-\kappa|\|\phi\|_\infty
   +\kappa_\eps\|\cK_\eps\phi-\phi\|_\infty
   +\left\|c_\eps(\phi)-\phi^2/2\right\|_\infty\\
 &\qquad\quad+\frac12\left\|\cK_\eps(\phi^2)-\phi^2\right\|_\infty
   +(1-a_\eps)\|c_\eps(\phi)\|_\infty\longrightarrow0.
\end{align*}
This proves \eqref{eq:global-martingale-integrand} with the stated global norm.
\end{proof}


\subsection{Microscopic exponential identity and the contribution of the prescribed past}\label{subsec:martingale-past}

For background on point-process exponentials and intensity changes, see
\cite{Bremaud}. For $\phi\in C_c^\infty(\Omega_T)$, define the point-process exponential
\begin{equation}\label{eq:expmart}
 \mathscr E_\eps(\phi):=\exp\left\{
 \delta_\eps\int_0^\infty\int_{\R^d}\phi(\delta_\eps t,\eps x)N^\eps(\dd t,\dd x)
 -\int_0^\infty\int_{\R^d}\left(e^{\delta_\eps\phi(\delta_\eps t,\eps x)}-1\right)
                \lambda^\eps(t,x)\dd t\dd x\right\}.
\end{equation}
Here $t$ and $x$ are microscopic coordinates. The
normalized past-parent compensator is
\begin{equation}\label{eq:pastcompensator-definition}
 \cI_\eps(g)=\frac{a_\eps\delta_\eps^2}{n_\eps}
    \int_{r<0}\int_{\R^d}\cJ_\eps g(r,\eps y)\Pi^\eps(\dd r,\dd y),
\end{equation}
where $\cJ_\eps$ is the crossing operator in \eqref{eq:Jpast}.
In particular, $\E[\cI_\eps(g)/\delta_\eps]=a_\eps \B_\eps g$.
The scaled functional of the prescribed past $\cI_\eps(\phi)/\delta_\eps$
appears in \eqref{eq:exactmartfunctional}.

The next lemma rewrites the exponential in \eqref{eq:expmart} as a functional of the event measure and identifies the contribution of the prescribed past at superexponential accuracy.

\begin{lemma}[Exponential identity and superexponential concentration of the past]\label{lem:martingale-past}
Under Assumption~\ref{ass:finite-mean}, the exponential
$\mathscr E_\eps(\phi)$ has conditional expectation at most one given
the past and satisfies
\begin{equation}\label{eq:exactmartfunctional}
 \frac1{n_\eps}\log\mathscr E_\eps(\phi)
 =\ip{\NN_\eps}{-\LL_\eps\phi-a_\eps \cK_\eps c_\eps(\phi)}
 -\frac{\cI_\eps(\phi)}{\delta_\eps}
 -\cI_\eps(c_\eps(\phi)).
\end{equation}
Moreover,
\begin{equation}\label{eq:pastsuperexp}
 \frac{\cI_\eps(\phi)}{\delta_\eps}
 \simeq m\int_{\R^d}\lambda_0(x)\phi(0,x)\dd x,\qquad
 \cI_\eps(c_\eps(\phi))\simeq0.
\end{equation}
Here $X_\eps\simeq x$ means that, for every $\eta>0$,
$n_\eps^{-1}\log\Pp(|X_\eps-x|>\eta)\to-\infty$.
\end{lemma}

\begin{proof}
Let $g_\eps(t,x):=\delta_\eps\phi(\delta_\eps t,\eps x)$, and let
$\mathscr E_{\eps,s}$ be the exponential in \eqref{eq:expmart} with
the two integrals restricted to $0<t\leq s$. Let
$A_{\eps,s}:=\log\mathscr E_{\eps,s}$. This process has finite variation
on bounded intervals, with continuous part
\[
 A_{\eps,s}^{\rm c}
 =-\int_0^s\int_{\R^d}\left(e^{g_\eps(r,x)}-1\right)\lambda^\eps(r,x)\dd x\dd r.
\]
At an event $(s,x)$ its jump is $g_\eps(s,x)$.
Let $K_\phi$ be a compact spatial set containing $\supp\phi$, in the
sense that $\phi(t,x)=0$ for $x\notin K_\phi$, and let
$K_{\eps,\phi}:=\eps^{-1}K_\phi$ and $H_\eps:=T/\delta_\eps$.
Only events with $x\in K_{\eps,\phi}$ contribute to these jumps.
Applying It\^o's formula with jumps to $a\mapsto e^a$
\cite[Chapter~II, Theorem~32]{Protter2004}, whose continuous quadratic
variation term is zero here, gives
\begin{align}\label{eq:point-process-exponential-jump}
 \mathscr E_{\eps,s}-1
 &=\int_0^s\mathscr E_{\eps,r-}\dd A_{\eps,r}^{\rm c}
   +\sum_{\substack{(r,x)\in N^\eps\\0<r\leq s,\ x\in K_{\eps,\phi}}}
          \mathscr E_{\eps,r-}\left(e^{g_\eps(r,x)}-1\right)\nonumber\\
 &=\int_0^s\int_{\R^d}\mathscr E_{\eps,r-}
       \left(e^{g_\eps(r,x)}-1\right)
       \left[N^\eps(\dd r,\dd x)-\lambda^\eps(r,x)\dd r\dd x\right].
\end{align}
Define the local event count and integrated intensity by
\[
 C_{\eps,s}:=N^\eps((0,s]\times K_{\eps,\phi}),\qquad
 B_{\eps,s}:=\int_0^s\int_{K_{\eps,\phi}}
                    \lambda^\eps(r,x)\dd x\dd r.
\]
The compensator property and \eqref{eq:meanbound} give
\begin{equation}\label{eq:local-exponential-finiteness}
 \E\left[C_{\eps,H_\eps}+B_{\eps,H_\eps}\right]
 =2\E\left[C_{\eps,H_\eps}\right]
 \leq2\eps^{-(\alpha+1-d)}\|\lambda_0\|_\infty
                    H_\eps|K_{\eps,\phi}|<\infty.
\end{equation}
Thus both quantities are finite almost surely. For $j\geq1$, let
\[
 \tau_j:=\inf\left\{s\geq0:C_{\eps,s}+B_{\eps,s}\geq j\right\},\qquad
 M_{g,\eps}:=\|g_\eps\|_\infty,\qquad d_{g,\eps}:=e^{M_{g,\eps}}-1,
\]
with $\inf\varnothing=\infty$. For $r\leq\tau_j$,
$C_{\eps,r-}+B_{\eps,r}\leq j$ because $B_{\eps,\cdot}$ is continuous,
and \eqref{eq:expmart} yields
\[
 0\leq\mathscr E_{\eps,r-}
 \leq\exp\left\{M_{g,\eps}C_{\eps,r-}+d_{g,\eps}B_{\eps,r}\right\}
 \leq e^{(M_{g,\eps}+d_{g,\eps})j},\qquad
 B_{\eps,H_\eps\wedge\tau_j}\leq j.
\]
Consequently the stopped integrand in
\eqref{eq:point-process-exponential-jump} satisfies
\[
 \int_0^{H_\eps\wedge\tau_j}\int_{\R^d}
  \left|\mathscr E_{\eps,r-}\left(e^{g_\eps(r,x)}-1\right)\right|^2
       \lambda^\eps(r,x)\dd x\dd r
 \leq d_{g,\eps}^2e^{2(M_{g,\eps}+d_{g,\eps})j}j.
\]
Its compensated integral is therefore a square-integrable martingale in
the filtration whose initial $\sigma$-field includes $\Pi^\eps$.
Taking conditional expectations in
\eqref{eq:point-process-exponential-jump} gives
$\E[\mathscr E_{\eps,H_\eps\wedge\tau_j}\mid\Pi^\eps]=1$.
By \eqref{eq:local-exponential-finiteness}, $\tau_j>H_\eps$ for all
sufficiently large $j$, almost surely. Since $g_\eps=0$ after $H_\eps$,
the conditional Fatou's lemma gives
\[
 \E\left[\mathscr E_\eps(\phi)\mid\Pi^\eps\right]
 =\E\left[\liminf_{j\to\infty}\mathscr E_{\eps,H_\eps\wedge\tau_j}
                                      \mid\Pi^\eps\right]
 \leq\liminf_{j\to\infty}
       \E\left[\mathscr E_{\eps,H_\eps\wedge\tau_j}\mid\Pi^\eps\right]=1.
\]

For a bounded compact test function $\psi$, Fubini's theorem applied to
\eqref{eq:intensity} gives the normalized compensator identity
\begin{equation}\label{eq:compensator-offspring-identity}
 \frac{\delta_\eps^2}{n_\eps}\int_0^\infty\int_{\R^d}
       \psi(\delta_\eps t,\eps x)\lambda^\eps(t,x)\dd x\dd t
 =\cI_\eps(\psi)+a_\eps\NN_\eps(\cK_\eps\psi).
\end{equation}
Absolute integrability follows from
\eqref{eq:normalized-mean-density-bound}, the spatial $L^1$ contraction
of $\cK_\eps$, and \eqref{eq:Bbounds}. Substitute
$e^{\delta_\eps\phi}-1=\delta_\eps\phi+\delta_\eps^2c_\eps(\phi)$
into \eqref{eq:expmart} and use
\eqref{eq:compensator-offspring-identity} for both test functions. Then
\begin{align*}
 n_\eps^{-1}\log\mathscr E_\eps(\phi)
 =\delta_\eps^{-1}\NN_\eps(\phi-a_\eps \cK_\eps\phi)
      -a_\eps\NN_\eps(\cK_\eps c_\eps(\phi))
 -\delta_\eps^{-1}\cI_\eps(\phi)
      -\cI_\eps(c_\eps(\phi)),
\end{align*}
which is \eqref{eq:exactmartfunctional} by \eqref{eq:discrete-generator}.

Let $\widetilde\Pi^\eps$ be the prescribed past in coordinates
$(r,x)=(r,\eps y)$; its intensity is
$(n_\eps/\delta_\eps)\lambda_0(x)\dd r\dd x$ on
$(-\infty,0)\times\R^d$ by \eqref{eq:past}. For
$X_\eps(\psi)=\cI_\eps(\psi)/\delta_\eps$, a jump of
$n_\eps X_\eps(\psi)$ is $a_\eps\delta_\eps \cJ_\eps\psi(r,x)$.
For a Poisson random measure $\Pi$ with intensity $\vartheta$ on
a space $\mathsf X$, and a bounded real $g\in L^1(\vartheta)$,
the Poisson exponential formula gives
\begin{equation}\label{eq:centered-Poisson-exponential-formula}
 \log\E\left[\exp\left\{\int_{\mathsf X}g\,\dd\Pi
                              -\int_{\mathsf X}g\,\dd\vartheta\right\}\right]
 =\int_{\mathsf X}\left(e^g-1-g\right)\dd\vartheta.
\end{equation}
Indeed, apply \cite[Theorem~3.9 and Exercise~3.4]{LP} separately
on $\{g<0\}$ and $\{g\geq0\}$, use independence of the two Poisson
restrictions, and subtract $\int_{\mathsf X}g\dd\vartheta$.
Boundedness and $L^1(\vartheta)$ integrability of $g$ give
$\int_{\mathsf X}|e^g-1|\dd\vartheta
 \leq e^{\|g\|_\infty}\int_{\mathsf X}|g|\dd\vartheta<\infty$.
Apply \eqref{eq:centered-Poisson-exponential-formula} with
$\mathsf X=(-\infty,0)\times\R^d$,
$\vartheta(\dd r,\dd x)=(n_\eps/\delta_\eps)\lambda_0(x)\dd r\dd x$,
and $g(r,x)=\theta a_\eps\delta_\eps\cJ_\eps\psi(r,x)$.
The bounds $\|g\|_\infty\leq|\theta|\delta_\eps\|\psi\|_\infty$ and
$\int_{\mathsf X}|g|\dd\vartheta
 \leq |\theta|n_\eps\B_\eps|\psi|<\infty$, from
\eqref{eq:Jpast} and \eqref{eq:Bbounds}, verify its hypotheses and give
\begin{align}\label{eq:past-centered-exponential-bound}
 &n_\eps^{-1}\log\E\left[e^{\theta n_\eps(X_\eps(\psi)-\E X_\eps(\psi))}\right]
 \nonumber\\
 &=\frac1{\delta_\eps}\int_{-\infty}^0\int_{\R^d}\lambda_0(x)
 \left[e^{\theta a_\eps\delta_\eps \cJ_\eps\psi}-1
                       -\theta a_\eps\delta_\eps \cJ_\eps\psi\right]\dd x\dd r
 \nonumber\\
 &\leq\frac{\theta^2\delta_\eps}{2}
       e^{|\theta|\delta_\eps\|\psi\|_\infty}
       \|\psi\|_\infty \B_\eps|\psi|
 \leq C_\psi\theta^2\delta_\eps
       e^{|\theta|\delta_\eps\|\psi\|_\infty},
\end{align}
where we used $|\cJ_\eps\psi|\leq\|\psi\|_\infty\overline H(-r)
\leq\|\psi\|_\infty$ and \eqref{eq:Bbounds}. The identity follows
first for finite parent regions and then by truncation, since the
integral of $|\cJ_\eps\psi|$ against $\lambda_0\dd r\dd x$ is finite.
By \eqref{eq:Blimit},
\[
\E[X_\eps(\phi)]=a_\eps \B_\eps\phi\to m\int_{\R^d}\lambda_0(x)\phi(0,x)\dd x.
\]
For any fixed $\theta>0$, Chernoff's inequality and
\eqref{eq:past-centered-exponential-bound} give
\[
 \limsup_{\eps\downarrow0}n_\eps^{-1}\log\Pp\left(
     \left|X_\eps(\phi)-m\int_{\R^d}\lambda_0(x)\phi(0,x)\dd x\right|>\eta\right)
 \leq-\theta\eta/2.
\]
Letting $\theta\uparrow\infty$ proves the first part of
\eqref{eq:pastsuperexp}. For $\psi_\eps=\delta_\eps c_\eps(\phi)$,
$\|\psi_\eps\|_\infty+\|\psi_\eps\|_{L^\infty([0,T];L^1(\R^d))}=\mathcal O_\phi(\delta_\eps)$.
Thus, $\E[X_\eps(\psi_\eps)]=\mathcal O_\phi(\delta_\eps)$ and the last bound
in \eqref{eq:past-centered-exponential-bound} is $\mathcal O_\phi(\theta^2\delta_\eps^3)$.
Since $X_\eps(\psi_\eps)=\cI_\eps(c_\eps(\phi))$, the same
two-sided Chernoff argument proves the second part of
\eqref{eq:pastsuperexp}.
This completes the proof.
\end{proof}

\subsection{Density-path representation and coercivity}\label{sec:spatial-completion}

Finite total mass turns the weak controlled equation into a continuous $L^1$ density path with a prescribed initial trace, and its mass estimate provides compact rate function sublevels in the following proposition.

\begin{proposition}[Finite-mass lifting, trace, and coercivity]
\label{prop:finite-mass-lift}
Retain Assumption~\ref{ass:finite-mean} and suppose
$L:=\int_{\R^d}\lambda_0(x)\dd x<\infty$; here $L=0$ is also permitted. If $\nu$ has finite
total mass and $J_{\rm sp}(\nu)<\infty$, it has the unique density-path
representation and control in \eqref{eq:full-spatial-controlled}.
Conversely that representation implies \eqref{eq:full-spatial-energy}.
Writing $M(t):=\int_{\R^d}q(t,x)\dd x$,  one has $\sqrt M\in H^1(0,T)$ and
\begin{equation}\label{eq:mass-energy-bound}
 \frac12\int_0^T\left(2m\left(\sqrt M\right)' +\kappa\sqrt M\right)^2\dd t
       \leq J_{\rm sp}(\nu).
\end{equation}
In particular,
\begin{equation}\label{eq:mass-coercivity}
 \sup_{0\leq t\leq T}\sqrt{M(t)}
 \leq\sqrt L+\sqrt{\frac{T J_{\rm sp}(\nu)}{2m^2}},
 \qquad
 \nu\left([0,T)\times\R^d\right)
 \leq T\left(\sqrt L+\sqrt{\frac{T J_{\rm sp}(\nu)}{2m^2}}\right)^2.
\end{equation}
The function in \eqref{eq:full-spatial-rate} is lower semicontinuous and
has compact sublevels in $\mathsf E$.
\end{proposition}
\begin{proof}
Let $v$ be the control from Proposition~\ref{prop:action} and set
$j=v\nu$. The Cauchy--Schwarz inequality and \eqref{eq:weakcontrolled} give
\begin{equation}\label{eq:residual-TV}
 |j|(\Omega_T)=\int_{\Omega_T}|v|\dd\nu
 \leq\left(\int_{\Omega_T}v^2\dd\nu\right)^{1/2}\left\{\nu(\Omega_T)\right\}^{1/2}
 =\left\{2J_{\rm sp}(\nu)\nu(\Omega_T)\right\}^{1/2}<\infty.
\end{equation}
Extend both measures by zero to $t<0$. The weak identity
\eqref{eq:weakcontrolled} becomes, on $(-\infty,T)\times\R^d$,
\begin{equation}\label{eq:measure-PDE-trace}
 m\partial_t\nu-\B^*\nu+\kappa\nu
 =j+m\delta_0(\dd t)\lambda_0(x)\dd x.
\end{equation}
For $\psi\in C_c^\infty(\R^d)$, let
$a_\psi(\dd t):=\int_{\R^d}\psi(x)\nu(\dd t,\dd x)$ and define
\[
 \sigma_\psi(\dd t):=\frac1m\int_{\R^d}(\B\psi-\kappa\psi)(x)
                                     \nu(\dd t,\dd x)
       +\frac1m\int_{\R^d}\psi(x)j(\dd t,\dd x)
       +\delta_0(\dd t)\int_{\R^d}\lambda_0(x)\psi(x)\dd x.
\]
Equation~\eqref{eq:measure-PDE-trace} gives
$\partial_ta_\psi=\sigma_\psi$ in distributions. This is a finite
signed measure, since
\[
 |\sigma_\psi|((-\infty,T))
 \leq m^{-1}\|\B\psi-\kappa\psi\|_\infty\nu(\Omega_T)
       +m^{-1}\|\psi\|_\infty|j|(\Omega_T)+L\|\psi\|_\infty.
\]
Let $b_\psi(t):=\sigma_\psi((-\infty,t])$. Integration against a
compact smooth temporal test function shows
$\partial_t(b_\psi\dd t)=\sigma_\psi$. Thus
$a_\psi-b_\psi\dd t$ has zero distributional derivative and is a
constant multiple of Lebesgue measure. Both terms vanish for negative
times, so the constant is zero. In particular $a_\psi$ is absolutely
continuous in time. Choosing nonnegative $\psi$ equal to one on
each compact spatial set proves
\[
 \nu(\dd t,\dd x)=\dd t\,\nu_t(\dd x),\qquad
 \int_0^T\nu_t(\R^d)\dd t=\nu(\Omega_T)<\infty.
\]
The same disintegration holds for $j=v\nu$, and
$\int_0^T|j_t|(\R^d)\dd t=|j|(\Omega_T)$.
Therefore the density of $a_\psi$ is absolutely continuous on $[0,T]$
and satisfies
\begin{equation}\label{eq:lift-spatial-test-equation}
 m\frac{\dd}{\dd t}\nu_t(\psi)
 =\nu_t(\B\psi-\kappa\psi)+j_t(\psi),\qquad
 \nu_0(\psi)=\int_{\R^d}\lambda_0\psi\dd x.
\end{equation}

Let the candidate finite-measure path be
\begin{equation}\label{eq:finite-mild-lift}
 q_t\dd x=\cQ_t^*(\lambda_0\dd x)
       +\frac1m\int_0^t\cQ_{t-s}^*j_s\dd s,
       \qquad 0\leq t\leq T.
\end{equation}
For $s<t$, $\cQ_{t-s}^*j_s$ has an $L^1$ heat-kernel density and
$\|\cQ_{t-s}^*j_s\|_1\leq|j_s|(\R^d)$. It follows that
\[
 \|q_t\|_1\leq L+\frac1m\int_0^t|j_s|(\R^d)\dd s,
 \qquad
 \int_0^T\|q_t\|_1\dd t\leq TL+\frac Tm|j|(\Omega_T).
\]
At zero, \eqref{eq:finite-mild-lift} gives
\[
 \|q_t-\lambda_0\|_1
 \leq\|\cQ_t\lambda_0-\lambda_0\|_1
       +\frac1m\int_0^t|j_s|(\R^d)\dd s\longrightarrow0,
\]
by strong continuity of the heat semigroup on $L^1$ and
\eqref{eq:residual-TV}.
For positive times, the convolution kernel of $\cQ_r$ in
\eqref{eq:operators} is
\begin{equation}\label{eq:killed-heat-kernel-density}
 p_r^\kappa(x)
 =\frac{m^{d/2}e^{-\kappa r/m}}{(2\pi r)^{d/2}\sqrt{\det\Sigma}}
       \exp\left\{-\frac{m x^\top\Sigma^{-1}x}{2r}\right\},
 \qquad r>0.
\end{equation}
Differentiating \eqref{eq:killed-heat-kernel-density} in $r$ gives
\begin{align}\label{eq:heat-time-derivative-bound}
 \partial_rp_r^\kappa(x)
 &=p_r^\kappa(x)\left(-\frac\kappa m-\frac d{2r}
                     +\frac{m x^\top\Sigma^{-1}x}{2r^2}\right),\nonumber\\
 \left\|\partial_rp_r^\kappa\right\|_1
 &\leq\frac\kappa m+\frac dr,\qquad r>0.
\end{align}
Here the Gaussian moment is
$\int_{\R^d}x^\top\Sigma^{-1}x\,p_r^\kappa(x)\dd x
=e^{-\kappa r/m}dr/m$.
Fix $0<t\leq T$, $0<\eta<t$, and $|h|<\eta/2$ with
$t+h\in[0,T]$. For $0\leq s\leq t-\eta$, all times between
$t-s$ and $t+h-s$ are at least $\eta/2$. Thus
\eqref{eq:heat-time-derivative-bound} implies
\[
 \left\|\left(\cQ_{t+h-s}^*-\cQ_{t-s}^*\right)j_s\right\|_1
 \leq |h|\left(\frac\kappa m+\frac{2d}{\eta}\right)|j_s|(\R^d).
\]
Splitting both integrals in \eqref{eq:finite-mild-lift} at $t-\eta$
and using the $L^1$ contraction of $\cQ_r$ on their remaining parts yields
\begin{align}\label{eq:finite-mild-time-continuity}
 \|q_{t+h}-q_t\|_1
 &\leq\|(\cQ_{t+h}-\cQ_t)\lambda_0\|_1\nonumber\\
 &\quad+\frac{|h|}{m}\left(\frac\kappa m+\frac{2d}{\eta}\right)
                      \int_0^{t-\eta}|j_s|(\R^d)\dd s
                      +\frac2m\int_{t-\eta}^{(t+|h|)\wedge T}|j_s|(\R^d)\dd s.
\end{align}
For fixed $\eta$, the first two terms tend to zero as $h\to0$.
The last term converges to $2m^{-1}\int_{t-\eta}^t|j_s|(\R^d)\dd s$,
which tends to zero as $\eta\downarrow0$ by \eqref{eq:residual-TV}.
This proves $q\in C([0,T];L^1(\R^d))$, including the terminal value.

The candidate solves \eqref{eq:lift-spatial-test-equation} by differentiating
the heat semigroup against $\psi$. To identify it with $\nu$, set
$\zeta=\nu-q\dd t\dd x$, a finite signed measure. It has zero initial
trace and $\zeta(\LL_\kappa\phi)=0$ for compact test functions.
For $f\in C_c^\infty([0,T')\times\R^d)$, $T'<T$, let $\phi=\cA f$ and
choose $\chi_R(x)=\chi(x/R)$, where $\chi=1$ on $B_1$ and vanishes
outside $B_2$. The product rule gives
\[
 \LL_\kappa(\chi_R\phi)
 =-\chi_R f+\phi\B\chi_R+(\Sigma\nabla\chi_R)\cdot\nabla\phi.
\]
The derivatives of $\phi$ are bounded by \eqref{eq:resolvent}, so that
\[
 |\zeta(f)|\leq C_f\left(R^{-1}+R^{-2}\right)|\zeta|(\Omega_T)\longrightarrow0,
\]
once $R$ contains the support of $f$. Thus $\zeta=0$ and
$\nu_t=q_t\dd x$ almost everywhere. Positivity of $\nu$ gives
$q_t\geq0$ almost everywhere in time; $L^1$ continuity extends this
to every $t$. Also $j_t=v(t,\cdot)q_t\dd x$ almost everywhere, so
\eqref{eq:finite-mild-lift} proves \eqref{eq:full-spatial-controlled}.
Conversely, the distributional equation \eqref{eq:full-spatial-controlled} with the $L^1$ trace
$q(0)=\lambda_0$, equivalently \eqref{eq:measure-PDE-trace} with
$\nu=q\dd t\dd x$ and $j=vq\dd t\dd x$,
 gives
\eqref{eq:weakcontrolled}; density of the test space and square
completion as in Proposition~\ref{prop:action} prove
\eqref{eq:full-spatial-energy} and uniqueness of the control.

Integrating \eqref{eq:finite-mild-lift} over space yields
\[
 M(t)=Le^{-\kappa t/m}
       +\frac1m\int_0^t e^{-\kappa(t-s)/m}j_s(\R^d)\dd s.
\]
It is absolutely continuous and satisfies
\begin{equation}\label{eq:spatial-mass-balance}
 m\dot M(t)+\kappa M(t)=\int_{\R^d}v(t,x)q(t,x)\dd x,
 \qquad M(0)=L.
\end{equation}
Let $e(t):=\int_{\R^d}(v(t,x))^2q(t,x)\dd x$. The Cauchy--Schwarz
inequality gives $|m\dot M+\kappa M|^2\leq M e$ almost everywhere.
For $u_\eta=\sqrt{M+\eta}$,
\[
 2m\dot u_\eta+\kappa u_\eta
 =\frac{m\dot M+\kappa M}{\sqrt{M+\eta}}
       +\frac{\kappa\eta}{\sqrt{M+\eta}},\qquad
 \left\|2m\dot u_\eta+\kappa u_\eta\right\|_2
 \leq\left\|\sqrt e\right\|_2+\kappa\sqrt{\eta T}.
\]
The integrating-factor formula bounds $u_\eta$ in $H^1$ uniformly.
Since $u_\eta\to\sqrt M$ uniformly, weak lower semicontinuity gives
\[
 \frac12\int_0^T\left[2m\left(\sqrt M\right)'(t)+\kappa\sqrt{M(t)}\right]^2\dd t
 \leq\frac12\int_0^T e(t)\dd t=J_{\rm sp}(\nu),
\]
which is \eqref{eq:mass-energy-bound}. For
$g=2m(\sqrt M)'+\kappa\sqrt M$, the same integrating factor gives
\[
 \sqrt{M(t)}=e^{-\kappa t/(2m)}\sqrt L
     +\frac1{2m}\int_0^t e^{-\kappa(t-s)/(2m)}g(s)\dd s
 \leq\sqrt L+\frac{\sqrt T}{2m}\|g\|_2.
\]
This proves both bounds in \eqref{eq:mass-coercivity}.

Finally let
\begin{equation}\label{eq:compact-total-mass}
 \mathsf K_R:=\left\{\nu\in\mathsf E:\nu(\Omega_T)\leq R\right\}.
\end{equation}
On this locally compact second-countable space, finite positive
measures identify with positive linear functionals on $C_0(\Omega_T)$.
The Banach--Alaoglu theorem \cite[Theorem~3.15]{Rudin1991}, applied
to $C_0(\Omega_T)^*$, makes its closed norm ball of radius $R$
compact for $\sigma(C_0(\Omega_T)^*,C_0(\Omega_T))$, the weak-star
topology. The positive cone is weak-star closed, since each inequality
$\nu(f)\geq0$ for $f\in C_0(\Omega_T)$, $f\geq0$, defines a
closed set. Thus the positive part of the ball is also compact.
On this ball, convergence on $C_c$ and on $C_0$ agree by uniform
approximation and the total-mass bound. Thus $\mathsf K_R$ is vaguely
compact. Total mass is lower semicontinuous because
$\nu(\Omega_T)=\sup_{0\leq f\leq1,\ f\in C_c(\Omega_T)}\nu(f)$.
The action $J_{\rm sp}$ is a supremum of continuous compact-test
functionals by \eqref{eq:spatialJ}. Each sublevel of
\eqref{eq:full-spatial-rate} lies in a fixed $\mathsf K_R$ by
\eqref{eq:mass-coercivity}; its vague limit has finite mass and no
larger $J_{\rm sp}$. These facts prove lower semicontinuity and goodness of the rate function.
The proof is complete.
\end{proof}

\subsection{Approximation of every finite-action density}\label{subsec:finite-recovery}
The finite-mass lower bound follows by approximating every finite-action
density with bounded smooth reproduction controls. The regularization
below preserves the prescribed initial density, including where that
density vanishes.

\begin{lemma}[Approximation by smooth controls at finite mass]\label{lem:finite-recovery}
Assume $0<L<\infty$. If $\mathcal J_{\lambda_0}(q\dd t\dd x)<\infty$,
there are real $v_k\in C_c^\infty(E_T)$ for which the solutions of
\eqref{eq:controlledPDE} satisfy
\begin{equation}\label{eq:finite-recovery-conclusion}
 \|q^{v_k}-q\|_{L^1(\Omega_T)}\longrightarrow0,\qquad
 \frac12\int_0^T\int_{\R^d}v_k^2q^{v_k}\dd x\dd t
       \longrightarrow\mathcal J_{\lambda_0}(q\dd t\dd x).
\end{equation}
\end{lemma}
\begin{proof}
Let $j=vq$ be the residual from Proposition~\ref{prop:finite-mass-lift},
so $q=\rho+\cA^*j$ and
$\int_{\Omega_T}j^2/q\dd t\dd x=2\mathcal J_{\lambda_0}(q\dd t\dd x)$.
In all expressions $j^{2}/q$, the value is zero if $j=q=0$ and $+\infty$ if $q=0$ and $j\neq 0$.
Both $q$ and $j$ belong to $L^1(\Omega_T)$, by finite mass and the
Cauchy--Schwarz inequality. In this proof all norms without a domain
are over $\Omega_T$. Let
\begin{equation}\label{eq:spatial-Gaussian-kernel}
 g_r^\Sigma(x)=\frac{1}{(2\pi r)^{d/2}\sqrt{\det\Sigma}}
       \exp\left\{-\frac{x^\top\Sigma^{-1}x}{2r}\right\},
 \qquad r>0,
\end{equation}
be the Gaussian density of covariance $r\Sigma$, and define the
forward semigroup
\begin{equation}\label{eq:finite-forward-smoothing-semigroup}
 (\mathsf Q_r f)(t,x)
 :=e^{-\kappa r}\one_{\{t>mr\}}
       \int_{\R^d}g_r^\Sigma(x-y)f(t-mr,y)\dd y.
\end{equation}
Its potential is $\cA^*=\int_0^\infty\mathsf Q_r\dd r$.
For $0<h<T/m$, average this semigroup with
\[
 w_h(r)=h^{-1}e^{-(r-h)/h}\one_{\{r\geq h\}},\qquad
 \mathsf S_h f=\int_h^\infty w_h(r)\mathsf Q_r f\dd r.
\]
The operators $\mathsf S_h$ are positive $L^1$ contractions, commute
with $\cA^*$, and converge strongly to the identity as $h\downarrow0$.
The hydrodynamic formula in \eqref{eq:operators} gives
$\mathsf Q_r\rho(t,x)=\rho(t,x)\one_{\{t>mr\}}$.
Consequently,
\begin{equation}\label{eq:finite-source-correction}
 d_h(t,x):=\rho(t,x)-\mathsf S_h\rho(t,x)
 =\rho(t,x)\int_{\max\{h,t/m\}}^\infty w_h(r)\dd r\geq0.
\end{equation}
Let $q_h:=\mathsf S_hq+d_h$ and $j_h:=\mathsf S_hj$.
Then $q_h=\rho+\cA^*j_h$, $q_h\geq0$, and $q_h\to q$ in $L^1$.
The Cauchy--Schwarz inequality for the kernel of $\mathsf S_h$ gives
\begin{equation}\label{eq:finite-smoothing-cost}
 \int_{\Omega_T}\frac{j_h^2}{q_h}\dd t\dd x
 \leq\int_{\Omega_T}\mathsf S_h(j^2/q)\dd t\dd x
 \leq\int_{\Omega_T}\frac{j^2}{q}\dd t\dd x.
\end{equation}

To truncate space, choose smooth $0\leq\chi_R\leq1$, equal to one
on $B_R$, zero outside $B_{2R}$, and flat at both ends of its
transition, so that
\[
 \left|\nabla\chi_R\right|^2\leq CR^{-2}\chi_R(1-\chi_R),\qquad
 \left|D^2\chi_R\right|^2\leq CR^{-4}\chi_R(1-\chi_R).
\]
Such a cutoff is obtained by rescaling the smooth transition built
from $e^{-1/s}\one_{\{s>0\}}$.
Write $a_h:=\mathsf S_hq$, $b_h:=\mathsf S_h\rho$, and set
\begin{align}\label{eq:finite-truncated-pair}
 q_{h,R}&:=\rho+\chi_R(a_h-b_h)
        =\chi_Ra_h+(1-\chi_R)b_h+d_h,\nonumber\\
 j_{h,R}&:=\chi_Rj_h-e_{h,R},\nonumber\\
 e_{h,R}&:=(\B\chi_R)(a_h-b_h)
                 +(\Sigma\nabla\chi_R)\cdot\nabla(a_h-b_h).
\end{align}
The product rule for $m\partial_t-\B+\kappa$ gives
$q_{h,R}=\rho+\cA^*j_{h,R}$ with unchanged initial value.
In particular $q_{h,R}\geq\chi_Rq_h$ and $q_{h,R}\geq d_h$.
Differentiating the explicit Gaussian kernel
\eqref{eq:spatial-Gaussian-kernel} gives
\[
 \Sigma^{1/2}\nabla_xg_r^\Sigma(x-y)
 =-\frac{\Sigma^{-1/2}(x-y)}r g_r^\Sigma(x-y).
\]
For $f\geq0$, applying the Cauchy--Schwarz inequality to
\eqref{eq:finite-forward-smoothing-semigroup} therefore yields
\begin{align*}
 \frac{|\Sigma^{1/2}\nabla\mathsf Q_r f(t,x)|^2}{\mathsf Q_r f(t,x)}
 &\leq e^{-\kappa r}\one_{\{t>mr\}}
   \int_{\R^d}\frac{|\Sigma^{-1/2}(x-y)|^2}{r^2}
       g_r^\Sigma(x-y)f(t-mr,y)\dd y.
\end{align*}
The Gaussian moment
$\int_{\R^d}|\Sigma^{-1/2}z|^2g_r^\Sigma(z)\dd z=dr$
and Tonelli's theorem give
\begin{equation}\label{eq:finite-semigroup-fisher}
 \int_{\Omega_T}\frac{|\Sigma^{1/2}\nabla\mathsf Q_r f|^2}
                              {\mathsf Q_r f}\dd t\dd x
 \leq\frac dr\|f\|_1,
\end{equation}
since the time translation and killing decrease the remaining mass.
Applying the Cauchy--Schwarz inequality also to the $r$ integral gives
\begin{equation}\label{eq:finite-smoothing-fisher}
\int_{\Omega_T}\frac{|\Sigma^{1/2}\nabla\mathsf S_h f|^2}
                              {\mathsf S_h f}\dd t\dd x
 \leq\int_h^\infty w_h(r)\frac dr\|f\|_1\dd r
 \leq\frac dh\|f\|_1.
\end{equation}
Using the cutoff bounds in \eqref{eq:finite-truncated-pair} and
then \eqref{eq:finite-smoothing-fisher}, we obtain
\begin{align}\label{eq:finite-cutoff-cost}
 \int_{\Omega_T}\frac{e_{h,R}^2}{q_{h,R}}\dd t\dd x
 &\leq CR^{-4}\int_{\Omega_T}(a_h+b_h)\dd t\dd x\nonumber\\
 &\quad+CR^{-2}\int_{\Omega_T}
       \left(\frac{|\Sigma^{1/2}\nabla a_h|^2}{a_h}
                    +\frac{|\Sigma^{1/2}\nabla b_h|^2}{b_h}\right)\dd t\dd x
 \nonumber\\
 &\leq C\left(R^{-4}+h^{-1}R^{-2}\right)(\|q\|_1+\|\rho\|_1).
\end{align}
The convention is zero at simultaneous zeros. For every $u>0$,
\[
 \int_{\Omega_T}\frac{j_{h,R}^2}{q_{h,R}}\dd t\dd x
 \leq(1+u)\int_{\Omega_T}\frac{j_h^2}{q_h}\dd t\dd x
       +\left(1+u^{-1}\right)\int_{\Omega_T}\frac{e_{h,R}^2}{q_{h,R}}\dd t\dd x.
\]
Since $\|q_{h,R}-q_h\|_1\to0$, first letting $R\to\infty$ and
then $u\downarrow0$ gives approximation of $q_h$ without increasing cost.

For the pair in \eqref{eq:finite-truncated-pair}, define the truncated
control by
\begin{equation}\label{eq:finite-truncated-control}
 v_{h,R}(t,x):=
 \begin{cases}
 j_{h,R}(t,x)/q_{h,R}(t,x),&q_{h,R}(t,x)>0,\\
 0,&q_{h,R}(t,x)=0.
 \end{cases}
\end{equation}
We next prove that \eqref{eq:finite-truncated-control} is bounded.
Changing variables $r=(t-s)/m$ writes
$\mathsf S_h$ as a space--time integral kernel
\[
 \frac1m w_h((t-s)/m)e^{-\kappa(t-s)/m}
       g_{(t-s)/m}^\Sigma(x-y)\one_{\{0<s<t-mh\}}.
\]
This kernel and its first spatial derivatives are bounded for fixed
$h$, so $\|D_x^\alpha\mathsf S_h f\|_\infty\leq C_h\|f\|_1$
for $|\alpha|\leq1$. Its convolution with an $L^1$ function is
continuous, by translation continuity in $L^1$.
For $t\leq mh$, both $j_{h,R}$ and $q_{h,R}-\rho$ vanish.
For $t\geq mh$, equation~\eqref{eq:finite-source-correction} gives
$d_h(t,x)=e^{1-t/(mh)}\rho(t,x)$. Since $L>0$, the positive heat
kernel in \eqref{eq:operators} makes $\rho$ continuous and strictly
positive on $[mh,T]\times\overline B_{2R}$. Thus,
\[
 \inf_{\substack{mh\leq t\leq T\\|x|\leq2R}}d_h(t,x)
 \geq e^{1-T/(mh)}
       \inf_{\substack{mh\leq t\leq T\\|x|\leq2R}}\rho(t,x)>0.
\]
Thus the control \eqref{eq:finite-truncated-control} is bounded,
continuous for $t>0$, and spatially compactly supported.

For bounded $w$, the mild form of \eqref{eq:controlledPDE} is
\begin{equation}\label{eq:controlled-mild-equation}
 q^w=\rho+\cA^*(wq^w).
\end{equation}
Let $M\geq\|w\|_\infty$ and
$\|f\|_{1,-z}=\int_0^T\int_{\R^d}e^{-zt}|f(t,x)|\dd x\dd t$.
The forward form of \eqref{eq:resolvent} and Fubini's theorem give
\begin{align}\label{eq:controlled-weighted-contraction}
 \left\|\cA^*(wf)\right\|_{1,-z}
 \leq\frac1m\int_0^T\int_s^T
     e^{-zt-\kappa(t-s)/m}\left\|w(s,\cdot)f(s,\cdot)\right\|_1\dd t\dd s
 \leq\frac{M}{\kappa+mz}\|f\|_{1,-z}.
\end{align}
Choose $z>0$ with $\kappa+mz>2M$, and let
$\vartheta=M/(\kappa+mz)<1$.
The Banach fixed-point theorem gives a unique $q^w\in L^1(\Omega_T)$,
since $e^{-zT}\|f\|_1\leq\|f\|_{1,-z}\leq\|f\|_1$.
To verify positivity, let
\[
 \rho_M(t)=e^{-Mt/m}\rho(t),\qquad
 (\cA_M^*f)(t)=\frac1m\int_0^t
                  e^{-M(t-s)/m}\cQ_{t-s}f(s)\dd s.
\]
The controlled equation is equivalently
\[
q^w=\rho_M+\cA_M^*\left((w+M)q^w\right).
\]
Its contraction constant is at most $2M/(\kappa+M+mz)<1$, and hence
\[
 q^w=\sum_{\ell=0}^\infty
             \bigl[\cA_M^*((w+M)\,\cdot)\bigr]^\ell\rho_M\geq0
\]
in $L^1$, because $\rho_M\geq0$ and $w+M\geq0$.

Now fix $h,R$, let $v=v_{h,R}$, and take $M=\|v\|_\infty$ in
\eqref{eq:controlled-weighted-contraction}.
Choose $\beta_k\in C_c^\infty(0,T)$ with $0\leq\beta_k\leq1$ and
$\beta_k(t)\to1$ for each $t\in(0,T)$.
Extend $\beta_kv$ by zero outside $E_T$ and let
\[
 v_{h,R,k}=\omega_{\gamma_k}*(\beta_kv),\qquad
 0<\gamma_k<\frac12\operatorname{dist}(\supp\beta_k,\{0,T\}),
 \qquad\gamma_k\downarrow0,
\]
where $\omega_{\gamma_k}$ is a nonnegative smooth mollifier of integral
one supported in the ball of radius $\gamma_k$ in $\R^{d+1}$.
The cutoffs may be chosen equal to one on an increasing compact
exhaustion of $(0,T)$. Continuity of $v$ on $E_T$ then gives
$v_{h,R,k}\to v$ pointwise, with
$v_{h,R,k}\in C_c^\infty(E_T)$ and $\|v_{h,R,k}\|_\infty\leq M$.
Let $q_k:=q^{v_{h,R,k}}$ and $\bar q:=q_{h,R}=q^v$; the last equality
follows from uniqueness and \eqref{eq:finite-truncated-pair}.
In particular, \eqref{eq:controlled-mild-equation} gives
\begin{equation}\label{eq:recovery-controlled-equations}
 q_k=\rho+\cA^*(v_{h,R,k}q_k),\qquad
 \bar q=\rho+\cA^*(v\bar q).
\end{equation}
Subtracting the two equations in \eqref{eq:recovery-controlled-equations}
gives
\[
 q_k-\bar q=
\cA^*\left\{v_{h,R,k}\left(q_k-\bar q\right)
                         +\left(v_{h,R,k}-v\right)\bar q\right\}.
\]
Applying \eqref{eq:controlled-weighted-contraction} yields
\begin{equation}\label{eq:controlled-density-stability}
 \|q_k-\bar q\|_1
 \leq\frac{e^{zT}}{(1-\vartheta)(\kappa+mz)}
          \left\|(v_{h,R,k}-v)\bar q\right\|_{1,-z}\longrightarrow0.
\end{equation}
Indeed, $|(v_{h,R,k}-v)\bar q|\leq2M\bar q\in L^1$,
so the dominated convergence theorem applies. Finally,
\begin{align}\label{eq:controlled-energy-stability}
 \left|\int_{\Omega_T}v_{h,R,k}^2q_k\dd t\dd x
       -\int_{\Omega_T}v^2\bar q\dd t\dd x\right|
 \leq M^2\left\|q_k-\bar q\right\|_1
 +\int_{\Omega_T}\left|v_{h,R,k}^2-v^2\right|\bar q\dd t\dd x
 \longrightarrow0.
\end{align}
The last integral converges to zero by the dominated convergence
theorem. Recall that $v=j_{h,R}/q_{h,R}$ where $q_{h,R}>0$ and $v=0$
where $q_{h,R}=0$. The finite-cost bound \eqref{eq:finite-cutoff-cost}
implies $j_{h,R}=0$ almost everywhere on $\{q_{h,R}=0\}$.
The quotient in the energy uses the convention
\[
 \frac{j_{h,R}^2}{q_{h,R}}
 :=\begin{cases}
      j_{h,R}^2/q_{h,R},&q_{h,R}>0,\\
      0,&q_{h,R}=0\text{ and }j_{h,R}=0,\\
      +\infty,&q_{h,R}=0\text{ and }j_{h,R}\ne0.
    \end{cases}
\]
Thus $v^2\bar q=j_{h,R}^2/q_{h,R}$ almost everywhere, with both sides
equal to zero where $q_{h,R}=j_{h,R}=0$.
Equations \eqref{eq:finite-smoothing-cost}--\eqref{eq:finite-cutoff-cost}
and a diagonal choice as $h\downarrow0$, $R,k\to\infty$ prove the
upper bound on the limiting costs in \eqref{eq:finite-recovery-conclusion}.
Lower semicontinuity from Proposition~\ref{prop:finite-mass-lift}
gives equality and completes the proof.
\end{proof}

\subsection{Proof of the full finite-mass LDP}

\begin{proof}[Proof of Theorem~\ref{thm:full-spatial-finite}]
For the upper bound, choose $0<c<r_T$ from Lemma~\ref{lem:analytic}.
Let $0\leq\chi_k\leq1$ be compact spatial cutoffs increasing to one,
and let $f_k(t,x):=c\one_{[0,T]}(t)\chi_k(x)$.
The support-independent bound \eqref{eq:vbounds} gives
$0\leq v_{\eps,f_k}\leq V:=C_Tc$ for all $k$ and sufficiently
small $\eps$. The inequality $e^u-1\leq ue^u$ for $u\geq0$,
\eqref{eq:Lambdaexact}, and \eqref{eq:past-crossing-representation}
give, uniformly in $k$,
\begin{align*}
 \Lambda_\eps(f_k)
 &\leq a_\eps e^{\delta_\eps V}\B_\eps v_{\eps,f_k}
 \leq e^{\delta_\eps V}mLV\leq C.
\end{align*}
Here $C$ depends only on the fixed model and $T$, since
$\delta_\eps\leq1$ for small $\eps$.
The monotone convergence theorem and \eqref{eq:superposition} yield
\begin{equation}\label{eq:spatial-mass-exponential-moment}
 \log\E\left[e^{n_\eps c\NN_\eps(\Omega_T)}\right]\leq n_\eps C.
\end{equation}
In particular, the event measure has finite total mass almost surely.
Chernoff's inequality applied to
\eqref{eq:spatial-mass-exponential-moment} gives
\begin{equation}\label{eq:total-mass-exp-tail}
 \Pp\left(\NN_\eps\notin\mathsf K_R\right)
 =\Pp\left(\NN_\eps(\Omega_T)>R\right)\leq\exp\left\{-n_\eps(cR-C)\right\}.
\end{equation}
The negligible deterministic endpoints have zero event mass.
The compact sets \eqref{eq:compact-total-mass} thus prove exponential
tightness in the vague topology.

For real-valued $\phi\in C_c^\infty(\Omega_T)$, define
\begin{equation}\label{eq:spatial-exponential-test-functional}
 \begin{aligned}
 \mathcal F_\phi(\nu)&:=-\nu\left(\LL_\kappa\phi+\phi^2/2\right)
             -m\int_{\R^d}\lambda_0(x)\phi(0,x)\dd x,\\
 J_{\rm sp}(\nu)&:=\sup_{\phi\in C_c^\infty(\Omega_T;\R)}
                         \mathcal F_\phi(\nu).
 \end{aligned}
\end{equation}
Thus the supremum is over real smooth compactly supported test functions on
$[0,T)\times\R^d$, as in \eqref{eq:spatialJ}; they may be nonzero
at $t=0$ and vanish near $t=T$.
Let $r_{\eps,\phi}$ be the supremum norm in
\eqref{eq:global-martingale-integrand}. On $\{\NN_\eps\in\mathsf K_R\}$,
\eqref{eq:exactmartfunctional} gives
\begin{align*}
 \left|n_\eps^{-1}\log\mathscr E_\eps(\phi)-\mathcal F_\phi(\NN_\eps)\right|
 \leq Rr_{\eps,\phi}
 +\left|\cI_\eps(\phi)/\delta_\eps-m\int_{\R^d}\lambda_0(x)\phi(0,x)\dd x\right|
 +|\cI_\eps(c_\eps(\phi))|.
\end{align*}
The last two terms exceed $2\eta$ only on an event
$B_{\eps,\phi,\eta}$ with superexponentially small probability by
\eqref{eq:pastsuperexp}. Since $\E\left[\mathscr E_\eps(\phi)\right]\leq1$, Markov's
inequality gives
\begin{equation}\label{eq:spatial-local-upper-estimate}
 \Pp\left(\NN_\eps\in\mathsf K_R,\ \mathcal F_\phi(\NN_\eps)>a\right)
 \leq e^{-n_\eps(a-Rr_{\eps,\phi}-2\eta)}
                         +\Pp(B_{\eps,\phi,\eta}).
\end{equation}
Let $K\subset\mathsf K_R$ be compact and $a<\inf_{\nu\in K}J_{\rm sp}(\nu)$.
For each $\nu\in K$ choose $\phi$ with $\mathcal F_\phi(\nu)>a$.
The continuous functions $\mathcal F_\phi$ give an open cover;
choose a finite subcover with test functions $\phi_1,\ldots,\phi_N$.
Apply \eqref{eq:spatial-local-upper-estimate} to each test function, then let
$\eps\downarrow0$ and $\eta\downarrow0$. The finite union gives
\[
 \limsup_{\eps\downarrow0}n_\eps^{-1}\log\Pp(\NN_\eps\in K)\leq-a.
\]
Taking $a\uparrow\inf_{\nu\in K}J_{\rm sp}(\nu)$ proves the compact upper bound.

For a closed $F\subset\mathsf E$, $F\cap\mathsf K_R$ is compact and
\eqref{eq:total-mass-exp-tail} yields
\begin{align*}
 \limsup_{\eps\downarrow0}n_\eps^{-1}\log\Pp(\NN_\eps\in F)
 &\leq\max\left\{-\inf_{\nu\in F\cap\mathsf K_R}J_{\rm sp}(\nu),-cR+C\right\}\\
 &\leq\max\left\{-\inf_{\nu\in F}\mathcal J_{\lambda_0}(\nu),-cR+C\right\}.
\end{align*}
On $\mathsf K_R$ the two actions agree by \eqref{eq:full-spatial-rate}.
Letting $R\uparrow\infty$ proves the full upper bound, including closed
sets of infinite-total-mass measures.

For the lower bound, let $O$ be vague open and
$\nu=q\dd t\dd x\in O$ have finite $\mathcal J_{\lambda_0}$.
Choose the approximating sequence in Lemma~\ref{lem:finite-recovery}.
For every sufficiently large fixed $k$, $q^{v_k}\dd t\dd x\in O$.
Lemma~\ref{lem:smooth-controlled-lower}, proved below under either
finite mass of the initial profile or the compact-kernel assumptions, gives
\[
 \liminf_{\eps\downarrow0}n_\eps^{-1}\log\Pp(\NN_\eps\in O)
 \geq-\frac12\int_0^T\int_{\R^d}\left(v_k(t,x)\right)^2q^{v_k}(t,x)\dd x\dd t.
\]
Let $k\to\infty$ and use \eqref{eq:finite-recovery-conclusion}.
Optimizing over finite-rate $\nu\in O$ proves the full lower bound.
The rate function and its goodness follow from Proposition~\ref{prop:finite-mass-lift}.
If its energy is zero, $v^2q=0$ almost everywhere, so $j=vq=0$.
Equation~\eqref{eq:finite-mild-lift} gives $q_t=\cQ_t\lambda_0=\rho_t$,
which proves the unique zero-cost assertion.
\end{proof}

\subsection{Compact-kernel bounds beyond finite mass}\label{subsec:infinite-bounds}

In this subsection, Assumptions~\ref{ass:finite-mean} and
\ref{ass:quantitative} hold. The bounded uniformly continuous density
$\lambda_0\geq0$ may have finite or infinite total mass.

For a real control $v\in C_c^\infty([0,T)\times\R^d)$, let $q^v$ be the bounded mild solution of
\begin{equation}\label{eq:controlledPDE}
 m\partial_t q^v=\B q^v-\kappa q^v+v q^v,
 \qquad q^v_0=\lambda_0,
\end{equation}
and set
\begin{equation}\label{eq:controlcost}
 \mathscr S(v):=\frac12\int_0^T\int_{\R^d}q^v(t,x)(v(t,x))^2\dd x\dd t.
\end{equation}
This cost is finite despite a possibly infinite initial total mass, because $v$ is compactly supported and $q^v$ is bounded.

For the compact-kernel subclass, localized exponential estimates
give an upper bound for all closed sets and a lower bound at the densities generated by smooth controls
in the following theorem.

\begin{theorem}[Spatial exponential bounds]\label{thm:spatialbounds}
Under Assumptions~\ref{ass:finite-mean} and \ref{ass:quantitative}, the laws of $\NN_\eps$ are exponentially tight at speed $n_\eps$ in the vague topology. For every closed set $F$,
\begin{equation}\label{eq:spatialupper}
\limsup_{\eps\to0}\frac1{n_\eps}\log\Pp(\NN_\eps\in F)
 \leq-\inf_{\nu\in F}J_{\rm sp}(\nu).
\end{equation}
For every open set $O$ and every smooth control $v$ with $q^v\dd t\dd x\in O$,
\begin{equation}\label{eq:spatiallower}
\liminf_{\eps\to0}\frac1{n_\eps}\log\Pp(\NN_\eps\in O)
 \geq-\mathscr S(v).
\end{equation}
Moreover, $J_{\rm sp}(q^v\dd t\dd x)=\mathscr S(v)$.
\end{theorem}

\begin{proof}[Proof of the upper bound in Theorem~\ref{thm:spatialbounds}]
Let $K_j$, $j\geq1$, exhaust $\Omega_T$ by compact sets with
$K_j\subset\operatorname{int}_{\Omega_T}K_{j+1}$, and choose
$f_j\in C_c^\infty(\Omega_T)$ with $0\leq f_j\leq1$ and $f_j=1$ on
$K_j$. The common radius in Lemma~\ref{lem:analytic} and
\eqref{eq:Bbounds} give a fixed $c>0$, a common small-$\eps$ cutoff,
and finite constants $C_j$ such that
$\Lambda_\eps(cf_j)\leq C_j$ for every $j$.
By \eqref{eq:superposition} and Chernoff's inequality,
\[
 \Pp\left(\NN_\eps(f_j)>R_j\right)\leq e^{-n_\eps(cR_j-C_j)}.
\]
Choose $cR_j\geq C_j+A+j+\log2$ and $n_\eps\geq1$. Then
\begin{equation}\label{eq:local-mass-exponential-exhaustion}
 \Pp\left(\NN_\eps\notin\mathsf K_A^{\rm loc}\right)
 \leq\sum_{j=1}^\infty e^{-n_\eps(A+j+\log2)}\leq e^{-n_\eps A},
 \quad
 \mathsf K_A^{\rm loc}=\bigcap_{j=1}^\infty\left\{\nu:\nu(f_j)\leq R_j\right\}.
\end{equation}
To prove vague compactness of $\mathsf K_A^{\rm loc}$, let
$(\nu_n)_{n\geq1}$ be any sequence in this set. For each $j$,
\begin{equation}\label{eq:local-compactness-mass-bound}
 \sup_{n\geq1}\nu_n(K_j)
 \leq\sup_{n\geq1}\nu_n(f_j)\leq R_j,
 \qquad
 |\nu_n(g)|\leq R_j\|g\|_\infty
 \quad\bigl(g\in C_c(\Omega_T),\ \supp g\subset K_j\bigr).
\end{equation}
The space $C_{K_j}(\Omega_T)$ of continuous functions supported in
$K_j$ is separable in the uniform norm. Choose a countable dense
family $\{g_{j,\ell}:\ell\geq1\}$ in each such space. The bounds
\eqref{eq:local-compactness-mass-bound} and diagonal extraction give
one subsequence $(\nu_{n_k})_{k\geq1}$ for which
$\nu_{n_k}(g_{j,\ell})$ converges for every $j,\ell$.
For $g\in C_{K_j}(\Omega_T)$ and any $k,k',\ell$,
\begin{equation}\label{eq:local-compactness-Cauchy}
 \left|\nu_{n_k}(g)-\nu_{n_{k'}}(g)\right|
 \leq2R_j\left\|g-g_{j,\ell}\right\|_\infty
       +\left|\nu_{n_k}(g_{j,\ell})-\nu_{n_{k'}}(g_{j,\ell})\right|.
\end{equation}
First approximate $g$ uniformly and then let $k,k'\to\infty$.
Equation~\eqref{eq:local-compactness-Cauchy} proves convergence for
every $g\in C_c(\Omega_T)$, since its support is contained in some
$K_j$. The limit $\ell(g):=\lim_{k\to\infty}\nu_{n_k}(g)$ is
linear, nonnegative for $g\geq0$, and satisfies
$|\ell(g)|\leq R_j\|g\|_\infty$ on $C_{K_j}(\Omega_T)$.
The Riesz representation theorem therefore gives a locally finite
positive Radon measure $\nu$ with $\ell(g)=\nu(g)$. In particular,
\begin{equation}\label{eq:local-compactness-cutoff-limit}
 \nu_{n_k}\longrightarrow\nu\quad\hbox{vaguely},\qquad
 \nu(f_j)=\lim_{k\to\infty}\nu_{n_k}(f_j)\leq R_j
 \quad(j\geq1).
\end{equation}
Thus $\nu\in\mathsf K_A^{\rm loc}$. On this locally bounded set,
the countable pairings $g_{j,\ell}$ determine the vague topology
by \eqref{eq:local-compactness-Cauchy}; for an enumeration $(g_i)_{i\geq1}$,
a compatible metric is
$\sum_{i=1}^\infty2^{-i}(1\wedge|\nu(g_i)-\nu'(g_i)|)$.
Sequential compactness is therefore compactness.
Together with \eqref{eq:local-mass-exponential-exhaustion}, this
proves exponential tightness. The source radius and small-$\eps$
cutoff used there are common to all $j$ by Lemma~\ref{lem:analytic}.

For a fixed test function $\phi$, let
\[
\eta_{\eps,\phi}:=-\LL_\eps\phi-a_\eps \cK_\eps c_\eps(\phi)
                 +\LL_\kappa\phi+\phi^2/2.
\]
Lemma~\ref{lem:global-generator} gives $\|\eta_{\eps,\phi}\|_\infty\to0$.
Assumption~\ref{ass:quantitative} makes its support lie in a fixed
compact set $K_\phi$ for all small $\eps$, because offspring shifts
have size at most $\eps R$ and $\delta_\eps s_+$.
For a vague compact set $K$ of measures, choose a cutoff $\psi=1$ on
$K_\phi$. The continuous pairing with $\psi$ gives
\[
 \sup_{\nu\in K}\left|\nu(\eta_{\eps,\phi})\right|
 \leq\left\|\eta_{\eps,\phi}\right\|_\infty\sup_{\nu\in K}\nu(\psi)\to0.
\]
Together with \eqref{eq:pastsuperexp} and \eqref{eq:exactmartfunctional},
this yields the local upper estimate
\eqref{eq:spatial-local-upper-estimate}, with the total-mass term
$Rr_{\eps,\phi}$ replaced by this local pairing bound.
For $a<\inf_{\nu\in K}J_{\rm sp}(\nu)$, the affine functionals in
\eqref{eq:spatial-exponential-test-functional} give a finite open
cover of $K$ by sets on which a selected functional exceeds $a$.
The exponential Markov inequality and the finite union give
\[
\limsup_{\eps\downarrow0}n_\eps^{-1}\log\Pp(\NN_\eps\in K)
 \leq-\inf_{\nu\in K}J_{\rm sp}(\nu).
\]
For a closed $F$, apply this to $F\cap\mathsf K_A^{\rm loc}$ and
use \eqref{eq:local-mass-exponential-exhaustion}:
\[
\limsup_{\eps\downarrow0}n_\eps^{-1}\log\Pp(\NN_\eps\in F)
 \leq\max\left\{-\inf_{\nu\in F}J_{\rm sp}(\nu),-A\right\}.
\]
Letting $A\uparrow\infty$ proves \eqref{eq:spatialupper}.
\end{proof}

The smooth-control lower bound uses a change of intensity whose hydrodynamic limit is the controlled equation \eqref{eq:controlledPDE}.

\begin{lemma}[Controlled law of large numbers]\label{lem:controlledLLN}
Assume either Assumption~\ref{ass:finite-mass}, or
Assumptions~\ref{ass:finite-mean} and \ref{ass:quantitative}.
Keep the law of the past unchanged and replace the positive-time intensity by
\begin{equation}\label{eq:tiltintensity}
 \lambda^{\eps,v}(t,x)=
 e^{\delta_\eps v(\delta_\eps t,\eps x)}
\lambda^\eps\left(t,x;N^{\eps,v},\Pi^\eps\right).
\end{equation}
For a bounded smooth compactly supported $v$, the corresponding event measure converges in probability, vaguely, to $q^v\dd t\dd x$, with $q^v$ from \eqref{eq:controlledPDE}. Its normalized mean is bounded uniformly on each finite macroscopic time interval.
\end{lemma}

\begin{proof}
Let $M:=\|v\|_\infty$ and $b_\eps:=a_\eps e^{\delta_\eps M}$.
For a parent at microscopic coordinates $(r,y)$, its positive-time
offspring have intensity density
\[
 a_\eps e^{\delta_\eps v(\delta_\eps t,\eps x)}
        h(t-r,x-y)\one_{\{t>\max\{r,0\}\}}
 \leq b_\eps h(t-r,x-y)\one_{\{t>\max\{r,0\}\}}.
\]
Thus replacing $\lambda_0$ by $\|\lambda_0\|_\infty$ and the factor
$a_\eps e^{\delta_\eps v(\delta_\eps t,\eps x)}$ by $b_\eps$
dominates the mean, by monotone iteration of the nonnegative mean equation.
Let $H_\delta(t):=\delta^{-1}H(t/\delta)$. Its normalized temporal
mean $R_\eps$ satisfies
$R_\eps=b_\eps\overline H(\cdot/\delta_\eps)+b_\eps H_{\delta_\eps}*R_\eps$.
Since
$(\Id-b_\eps H_{\delta_\eps}*)1=1-b_\eps+b_\eps\overline H(\cdot/\delta_\eps)$,
iteration gives
\begin{equation}\label{eq:tilted-mean-renewal}
 R_\eps(t)=1+(b_\eps-1)\sum_{j=0}^\infty b_\eps^j
             \Pp\left(\delta_\eps T_j\leq t\right).
\end{equation}
Choose $q_*>0$ with $mq_*+\kappa>M$. Then
\begin{align*}
 \sum_{j=0}^\infty b_\eps^j\Pp\left(\delta_\eps T_j\leq t\right)
 &\leq e^{q_*T}\sum_{j=0}^\infty
                   \left[b_\eps\widehat H(q_*\delta_\eps)\right]^j
 =\frac{e^{q_*T}}{1-b_\eps\widehat H(q_*\delta_\eps)},\\
 1-b_\eps\widehat H(q_*\delta_\eps)
 &=\delta_\eps(mq_*+\kappa-M+o(1)),\qquad |b_\eps-1|\leq C_M\delta_\eps.
\end{align*}
Thus \eqref{eq:tilted-mean-renewal} gives the uniform bound
\begin{equation}\label{eq:controlled-mean-density-bound}
 \sup_{0\leq t\leq T}R_\eps(t)\leq C_{T,M},\qquad
 \E\left[\NN_\eps^v(\dd t,\dd x)\right]
       \leq C_{T,M}\|\lambda_0\|_\infty\dd t\dd x.
\end{equation}
For the compact-kernel model, fix a bounded microscopic spatial
set $K$ and the microscopic horizon $H_\eps:=T/\delta_\eps$.
If a prescribed parent has a positive-time descendant by $H_\eps$,
its first child has time $t_1>0$ and its own time lies in $[-s_+,0)$.
Along a line of descent with $k$ positive-time events,
\begin{equation}\label{eq:controlled-finite-generation-depth}
 t_k\geq t_1+(k-1)s_-,\qquad
 k\leq J_\eps:=1+\lfloor H_\eps/s_-\rfloor.
\end{equation}
If the last event lies in $K$, the displacement bound in
Assumption~\ref{ass:quantitative} places its prescribed ancestor in
$K^{+J_\eps R}:=\left\{y:\operatorname{dist}(y,K)\leq J_\eps R\right\}$.
The number of possible prescribed ancestors is almost surely finite,
because \eqref{eq:past} gives
\begin{equation}\label{eq:controlled-local-parent-count}
 \E\left[\Pi^\eps\left([-s_+,0)\times K^{+J_\eps R}\right)\right]
 \leq\eps^{-(\alpha+1-d)}\|\lambda_0\|_\infty
              s_+\left|K^{+J_\eps R}\right|<\infty.
\end{equation}
For a parent at $(r,y)$, give its children the Poisson intensity
\begin{equation}\label{eq:controlled-offspring-construction}
 a_\eps h(s,z)
 e^{\delta_\eps v(\delta_\eps(r+s),\eps(y+z))}
 \one_{\{r+s>0\}}\dd s\dd z,
 \qquad (s,z)\in(0,\infty)\times\R^d.
\end{equation}
Its total mass is at most $b_\eps$, so every parent has finitely
many children almost surely. Equations
\eqref{eq:controlled-finite-generation-depth} and
\eqref{eq:controlled-local-parent-count} show that only finitely many
ancestors and generations can contribute to $(0,H_\eps]\times K$.
Using the same offspring Poisson measures for nested observation
regions gives consistent local constructions. Superposition of
\eqref{eq:controlled-offspring-construction} over past parents and
positive-time events yields exactly \eqref{eq:tiltintensity}.
Under Assumption~\ref{ass:finite-mass}, integration
of the mean iteration gives
\begin{equation}\label{eq:finite-controlled-total-mean}
 \E^v\left[\NN_\eps^v(\Omega_T)\right]\leq C_{T,M}LT.
\end{equation}
The bound holds for every finite-generation construction; increasing
the generation cutoff gives a process with finite expected total count
on the observation horizon, and hence no explosion. In both cases, the
count and compensator are finite on every compact observation region.

For a compact deterministic test function $g$, define the equation-scale martingale
\begin{equation}\label{eq:controlled-equation-martingale}
 \mathcal M_\eps^v(g)=\frac{\delta_\eps}{n_\eps}
 \int_0^\infty\int_{\R^d}g(\delta_\eps t,\eps x)
    \left[N^{\eps,v}(\dd t,\dd x)-\lambda^{\eps,v}(t,x)\dd t\dd x\right].
\end{equation}
The compensated-integral isometry and
\eqref{eq:controlled-mean-density-bound} give
\begin{equation}\label{eq:controlled-martingale-second-moment}
 \E^v\left[|\mathcal M_\eps^v(g)|^2\right]
 =n_\eps^{-1}\E^v\left[\NN_\eps^v(g^2)\right]
 \leq\frac{C_{T,M}\|\lambda_0\|_\infty}{n_\eps}
                       \int_0^T\int_{\R^d}g^2\dd x\dd t.
\end{equation}
Let $\mathcal R_\eps^v$ be the original-intensity compensator evaluated
along the tilted process:
\[
 \mathcal R_\eps^v(g)=\frac{\delta_\eps^2}{n_\eps}
 \int_0^\infty\int_{\R^d}g(\delta_\eps t,\eps x)
       \lambda^\eps(t,x;N^{\eps,v},\Pi^\eps)\dd x\dd t.
\]
Fubini's theorem as in \eqref{eq:compensator-offspring-identity} and
\eqref{eq:tiltintensity} give the exact identities
\begin{align}\label{eq:controlled-compensator-identities}
 \mathcal R_\eps^v(g)&=\cI_\eps(g)+a_\eps\NN_\eps^v(\cK_\eps g),\nonumber\\
 \NN_\eps^v(g)&=\mathcal R_\eps^v(e^{\delta_\eps v}g)
                                      +\delta_\eps\mathcal M_\eps^v(g).
\end{align}
Let $v_\eps=(1-e^{-\delta_\eps v})/\delta_\eps$; then
$\|v_\eps-v\|_\infty\leq C_M\delta_\eps$.
Applying the second identity in \eqref{eq:controlled-compensator-identities}
to $e^{-\delta_\eps v}g$ and subtracting the first identity at $g$ yields
\begin{equation}\label{eq:controlled-weak-microscopic-equation}
 \NN_\eps^v(-\LL_\eps g-v_\eps g)
 =\frac{\cI_\eps(g)}{\delta_\eps}
       +\mathcal M_\eps^v(g)-\delta_\eps\mathcal M_\eps^v(v_\eps g).
\end{equation}
The two martingale terms tend to zero in $L^2$ by
\eqref{eq:controlled-martingale-second-moment}. The past law is unchanged,
so \eqref{eq:pastsuperexp} applies to its term. Let $e_\eps:=\LL_\eps g+v_\eps g-\LL_\kappa g-vg$.
The global generator convergence gives $\|e_\eps\|_\infty\to0$.
When the initial profile has finite mass, \eqref{eq:finite-controlled-total-mean}
yields
\[
 \E^v\left[|\NN_\eps^v(e_\eps)|\right]\leq C_{T,M}LT\|e_\eps\|_\infty\to0.
\]
Under the compact-kernel assumptions, $e_\eps$ has support in a
common compact set $K_g$, so
\eqref{eq:controlled-mean-density-bound} instead gives the upper bound
$C_{T,M}\|\lambda_0\|_\infty|K_g|\|e_\eps\|_\infty\to0$.

The mean bound \eqref{eq:controlled-mean-density-bound}
 and a countable compact exhaustion give tightness
of the positive measures. Any limiting random measure $\nu$ satisfies
\begin{equation}\label{eq:controlled-limit-weak-equation}
 -\nu(\LL_\kappa g+vg)
 =m\int_{\R^d}\lambda_0(x)g(0,x)\dd x,
 \qquad g\in C_c^\infty(\Omega_T).
\end{equation}
A countable dense test family and continuity of the distributional
pairing make the identity hold simultaneously almost surely.
Moreover,
\begin{equation}\label{eq:controlled-limit-mean-density-bound}
\E\left[\nu(f)\right]\leq C_{T,M}\|\lambda_0\|_\infty
\int_0^T\int_{\R^d}f(t,x)\dd x\dd t,
\end{equation}
for
nonnegative compact $f$, by applying weak convergence to truncated
nonnegative pairings and then the monotone convergence theorem.

We prove uniqueness in the class satisfying
\eqref{eq:controlled-limit-mean-density-bound}. On $C_0(\Omega_T)$,
with zero extension for $t\geq T$, define
\[
 (\mathsf T_r w)(t,x)
 :=e^{-\kappa r}\cP_{mr}[w(t+mr,\cdot)](x),\qquad
 \cV w:=vw.
\]
The operators $\mathsf T_r$ form a strongly continuous contraction
semigroup with generator $\LL_\kappa$ and potential $\cA$ from
\eqref{eq:resolvent}; also $\|\cV\|\leq M$.
The Dyson--Phillips expansion for a bounded perturbation
\cite[Theorem~III.1.10]{EngelNagel2000}, integrated over genealogical
time, gives the following Duhamel series for the backward solution
$-(\LL_\kappa+v)\phi=f$, $\phi(T)=0$:
\begin{equation}\label{eq:controlled-backward-duality}
 \phi=\sum_{j=0}^\infty(\cA\cV)^j\cA f
 =\cA f+\cA(v\cA f)+\cA(v\cA(v\cA f))+\cdots.
\end{equation}
Here $\mathsf T_r=0$ for $r\geq T/m$, so the potential only involves
a finite interval. The same is true for each perturbed term, since
multiplication by $v$ does not change the time coordinate.
Convergence can also be checked directly: the time-ordered simplex
in the iterated
resolvent and the $L^p$ contraction of $\cQ_t$ give
\[
 \left\|(\cA \cV )^j\cA f\right\|_{0,p}
 \leq\frac{(T/m)^{j+1}M^j}{(j+1)!}\|f\|_{0,p},
 \qquad j\in\mathbb N_0,\quad p\in\{1,\infty\}.
\]
Thus \eqref{eq:controlled-backward-duality} converges and solves the
equation. Differentiating this series, with the product rule for the
compact smooth $v$, gives bounded derivatives. Choose $R_0$ containing
the spatial supports of $f$ and $v$. By \eqref{eq:resolvent},
$\phi=\cA (f+v\phi)$, with bounded forcing supported in $B_{R_0}$.
The Gaussian heat kernel and its first derivatives satisfy, for
$|\beta|\leq1$ and $|x|\geq2R_0+1$,
\[
 \sup_{0\leq t\leq T}\left|D_x^\beta\phi(t,x)\right|
 \leq C\int_0^T s^{-(d+|\beta|)/2}
             e^{-(|x|-R_0)^2/(Cs)}\dd s
 \leq C_{T,f,v}e^{-c_T|x|^2}.
\]
This bounds the cutoff errors by integrable functions independent of $R$.
Apply \eqref{eq:controlled-limit-weak-equation} to $\chi_R\phi$.
The product rule and the bounded-mean estimate \eqref{eq:controlled-limit-mean-density-bound}
 give
\begin{align*}
 \E\left[\left|\nu(f)-m\int_{\R^d}\lambda_0(x)\phi(0,x)\dd x\right|\right]
 &\leq C\int_0^T\int_{\{|x|\geq R\}}
             \left(R^{-1}|\nabla\phi|+R^{-2}|\phi|\right)\dd x\dd t\\
 &\quad+m\|\lambda_0\|_\infty
                \int_{\{|x|\geq R\}}|\phi(0,x)|\dd x\longrightarrow0.
\end{align*}
The same semigroup duality applied to \eqref{eq:controlledPDE} gives
\[
\int_0^T\int_{\R^d}q^v(t,x)f(t,x)\dd x\dd t=m\int_{\R^d}\lambda_0(x)\phi(0,x)\dd x.
\]
Thus every limiting measure equals $q^v\dd t\dd x$ almost surely.
Tightness and this deterministic identification prove convergence
in probability. Equation~\eqref{eq:controlled-mean-density-bound}
proves the asserted uniform mean bound.
\end{proof}

The likelihood of the preceding intensity change gives the same
smooth-control lower bound under both sets of assumptions in
Lemma~\ref{lem:controlledLLN}.

\begin{lemma}[Lower bound along smooth controls]\label{lem:smooth-controlled-lower}
Under either set of assumptions in Lemma~\ref{lem:controlledLLN},
for every bounded smooth compactly supported $v$ and every vague
open $O$ containing $q^v\dd t\dd x$, the bound
\begin{align*}
\liminf_{\varepsilon\downarrow 0}n_{\varepsilon}^{-1}\log\mathbb{P}(\mathcal{N}_{\varepsilon}\in O)\geq-\mathscr S(v)
\end{align*}
holds and
$J_{\rm sp}(q^v\dd t\dd x)=\mathscr S(v)$.
This lemma will also give the lower bound in
Theorem~\ref{thm:spatialbounds}.
\end{lemma}

\begin{proof}
Let $M:=\|v\|_\infty$, $t_\eps:=T/\delta_\eps$, and let a compact
spatial set $K$ contain the spatial support of $v$. Write
$K_\eps:=\eps^{-1}K$ and $u_\eps(t,x):=\delta_\eps v(\delta_\eps t,\eps x)$.
On the canonical space with past $\Pi^\eps$ and positive-time event
measure $N$, let
\[
 \lambda_\eps^0(t,x):=\lambda^\eps(t,x;N,\Pi^\eps),\qquad
 \lambda_\eps^v(t,x):=e^{u_\eps(t,x)}\lambda_\eps^0(t,x),\qquad
 \Pp^0:=\Pp.
\]
Under $\Pp^v$, $N$ has the law of $N^{\eps,v}$. We first justify
the likelihood identity on the observation horizon:
\begin{equation}\label{eq:likelihood}
 \begin{aligned}
 \ell_\eps^v:=\log\frac{\dd\Pp^v}{\dd\Pp}
 &=\delta_\eps\int_0^\infty\int_{\R^d}
                 v(\delta_\eps t,\eps x)N^{\eps,v}(\dd t,\dd x)\\
 &\quad-\int_0^\infty\int_{\R^d}
       \left[e^{\delta_\eps v(\delta_\eps t,\eps x)}-1\right]
                      \lambda_\eps^0(t,x)\dd x\dd t.
 \end{aligned}
\end{equation}
All likelihood terms are supported in $(0,t_\eps]\times K_\eps$.
For $0\leq t\leq t_\eps$ and $w\in\{0,v\}$, define
\[
 C_\eps(t)=N((0,t]\times K_\eps),\qquad
 A_\eps^w(t)=\int_0^t\int_{K_\eps}\lambda_\eps^w(s,x)\dd x\dd s,
\]
and let
\[
 \sigma_{\eps,j}:=\inf\left\{t\in[0,t_\eps]:
             C_\eps(t)+A_\eps^0(t)+A_\eps^v(t)\geq j\right\},
\]
with the convention that $\inf\varnothing=\infty$.
The likelihood process on this region is
\begin{equation}\label{eq:localized-likelihood-process}
 Z_\eps(t)=\exp\left\{
  \int_0^t\int_{K_\eps}u_\eps(s,x)N(\dd s,\dd x)
  -\int_0^t\int_{K_\eps}\left(e^{u_\eps(s,x)}-1\right)
                                      \lambda_\eps^0(s,x)\dd x\dd s\right\}.
\end{equation}
For $\tau=t\wedge\sigma_{\eps,j}$, this is the product of the
event-hazard ratios and the no-event factors:
\[
 Z_\eps(\tau)
 =\left(\prod_{\substack{(s,x)\text{ an atom of }N\\0<s\leq\tau,\ x\in K_\eps}}
       e^{u_\eps(s,x)}\right)
       \exp\left\{-\int_0^\tau\int_{K_\eps}
                \left(\lambda_\eps^v-\lambda_\eps^0\right)(s,x)\dd x\dd s\right\}.
\]
In particular, $Z_\eps(0)=1$, and the jump product rule gives
\[
 Z_\eps(t)=1+\int_0^t\int_{K_\eps}Z_\eps(s-)\left(e^{u_\eps(s,x)}-1\right)
           \left[N(\dd s,\dd x)-\lambda_\eps^0(s,x)\dd x\dd s\right].
\]
Since $|u_\eps|\leq\delta_\eps M$ and the count jumps by at most one,
\[
 0<Z_\eps\left(t\wedge\sigma_{\eps,j}\right)
 \leq\exp\left\{\delta_\eps M(j+1)+\left(e^{\delta_\eps M}+1\right)j\right\}.
\]
Thus, the stopped local martingale is a bounded martingale, with
$\E[Z_\eps(t_\eps\wedge\sigma_{\eps,j})\mid\Pi^\eps]=1$.
Changing measure by this density preserves the past law and changes
the predictable intensity to $\lambda_\eps^v$ up to $\sigma_{\eps,j}$,
as follows by the product rule for the compensated event integrals.
The conditional event-hazard construction therefore agrees with
$\Pp^v$ up to this stopping time.

To remove the stopping, use
$e^{-\delta_\eps M}\lambda_\eps^0\leq\lambda_\eps^v
\leq e^{\delta_\eps M}\lambda_\eps^0$ and compensation under
$\Pp^w$. Equations~\eqref{eq:event} and
\eqref{eq:controlled-mean-density-bound}, also with $v=0$, give
\begin{align}\label{eq:likelihood-stopping-bound}
 \E^w\left[C_\eps(t_\eps)\right]
 &=\E^w\left[A_\eps^w(t_\eps)\right]
 \leq\frac{n_\eps}{\delta_\eps^2}
       C_{T,M}\|\lambda_0\|_\infty T|K|,\nonumber\\
 \Pp^w\left(\sigma_{\eps,j}\leq t_\eps\right)
 &\leq\frac{\E^w[C_\eps(t_\eps)+A_\eps^0(t_\eps)+A_\eps^v(t_\eps)]}{j}
 \nonumber\\
 &\leq\frac{(2+e^{\delta_\eps M})n_\eps}
                  {j\delta_\eps^2}
       C_{T,M}\|\lambda_0\|_\infty T|K|
 \xrightarrow[j\to\infty]{}0,\qquad w\in\{0,v\}.
\end{align}
Here $\eps$ is fixed and sufficiently small while $j\to\infty$.
For every event $B$ measurable on the observation horizon, the stopped
density identity on $G_{\eps,j}:=\{\sigma_{\eps,j}>t_\eps\}$ is
\[
 \Pp^v(B\cap G_{\eps,j})
   =\E\left[Z_\eps(t_\eps)\one_{B\cap G_{\eps,j}}\right].
\]
Since $G_{\eps,j}$ increases to an event of probability one under both
laws by \eqref{eq:likelihood-stopping-bound}, the monotone convergence theorem gives
$\Pp^v(B)=\E[Z_\eps(t_\eps)\one_B]$.
The density $Z_\eps(t_\eps)$ is strictly positive and finite, so the
two laws are equivalent. Taking its logarithm under $\Pp^v$ proves
\eqref{eq:likelihood}.

Under $\Pp^v$, compensation in \eqref{eq:likelihood} gives
\begin{equation}\label{eq:controlled-likelihood-decomposition}
 n_\eps^{-1}\ell_\eps^v
 =\mathcal M_\eps^v(v)+\mathcal R_\eps^v(\mathfrak h_\eps(v)),\qquad
 \mathfrak h_\eps(z)=\frac{e^{\delta_\eps z}\delta_\eps z
                          -e^{\delta_\eps z}+1}{\delta_\eps^2}.
\end{equation}
For $|z|\leq M$, Taylor's integral remainder gives
$\left|\mathfrak h_\eps(z)-z^2/2\right|\leq C_M\delta_\eps|z|^3$.
By \eqref{eq:controlled-martingale-second-moment},
\[
 \E^v\left[|\mathcal M_\eps^v(v)|^2\right]
 \leq n_\eps^{-1}C_{T,M}\|\lambda_0\|_\infty
                               \int_0^T\int_{\R^d}(v(t,x))^2\dd x\dd t\to0.
\]
The compensator identity used in
\eqref{eq:controlled-weak-microscopic-equation} also gives
\[
 \mathcal R_\eps^v(g)
 =\NN_\eps^v\left(e^{-\delta_\eps v}g\right)
                  -\delta_\eps\mathcal M_\eps^v\left(e^{-\delta_\eps v}g\right).
\]
For $g=\mathfrak h_\eps(v)$, its support is fixed and
$e^{-\delta_\eps v}\mathfrak h_\eps(v)\to v^2/2$ uniformly.
More explicitly, let $F_\eps=e^{-\delta_\eps v}\mathfrak h_\eps(v)$.
Then $\supp F_\eps\subset\supp v$ and
$\|F_\eps-v^2/2\|_\infty\leq C_M\delta_\eps$.
The uniform mean bound \eqref{eq:controlled-mean-density-bound} and
the martingale estimate \eqref{eq:controlled-martingale-second-moment}
give
\begin{align}\label{eq:controlled-likelihood-test-error}
 \E^v\left[\left|\NN_\eps^v\left(F_\eps-v^2/2\right)\right|\right]
 &\leq C_{T,M}\|\lambda_0\|_\infty T|K|
                         \left\|F_\eps-v^2/2\right\|_\infty\longrightarrow0,
 \nonumber\\
 \E^v\left[\left|\delta_\eps\mathcal M_\eps^v(F_\eps)\right|^2\right]
 &\leq\frac{\delta_\eps^2}{n_\eps}C_{T,M}\|\lambda_0\|_\infty
              \int_0^T\int_{\R^d}|F_\eps(t,x)|^2\dd x\dd t
              \longrightarrow0.
\end{align}
Consequently, \eqref{eq:controlled-likelihood-decomposition},
\eqref{eq:controlled-likelihood-test-error}, and
Lemma~\ref{lem:controlledLLN}, applied to the fixed test function $v^2/2$, yield
\begin{equation}\label{eq:controlled-likelihood-limit}
 n_\eps^{-1}\ell_\eps^v\longrightarrow
 \frac12\int_0^T\int_{\R^d}q^v(t,x)(v(t,x))^2\dd x\dd t
 =\mathscr S(v)\quad\text{in }\Pp^v\text{-probability}.
\end{equation}
If $O$ is open and contains $q^v\dd t\dd x$, let
$B_\eps:=\left\{\NN_\eps^v\in O,\ \ell_\eps^v\leq n_\eps(\mathscr S(v)+\eta)\right\}$.
By \eqref{eq:controlled-likelihood-limit}, $\Pp^v(B_\eps)\to1$, and
\[
 \Pp(\NN_\eps\in O)
 =\E^v\left[e^{-\ell_\eps^v}\one_{\{\NN_\eps^v\in O\}}\right]
 \geq e^{-n_\eps(\mathscr S(v)+\eta)}\Pp^v(B_\eps).
\]
Taking logarithms, then $\eps\downarrow0$ and $\eta\downarrow0$,
proves \eqref{eq:spatiallower}.

Finally, \eqref{eq:controlledPDE} with its initial value gives, for
$\phi\in C_c^\infty(\Omega_T)$,
\[
 -\int_0^T\int_{\R^d}q^v\LL_\kappa\phi\dd x\dd t
       -m\int_{\R^d}\lambda_0(x)\phi(0,x)\dd x
 =\int_0^T\int_{\R^d}q^vv\phi\dd x\dd t.
\]
The terminal term vanishes because $\phi$ vanishes near $T$;
its spatial compact support and the bounded mild solution justify the
distributional integration by parts. Thus the affine functional in
\eqref{eq:spatialJ} is
\[
 \int_0^T\int_{\R^d}q^v\left(v\phi-\phi^2/2\right)\dd x\dd t
 =\mathscr S(v)-\frac12\int_0^T\int_{\R^d}q^v(\phi-v)^2\dd x\dd t.
\]
It is at most $\mathscr S(v)$, with equality at the admissible test
$\phi=v$. This proves $J_{\rm sp}(q^v\dd t\dd x)=\mathscr S(v)$.
The proof is complete.
\end{proof}

\subsection{Approximation conditions and a zero initial profile}\label{subsec:recovery}

\begin{definition}[Approximation by smooth controls]\label{def:recovery}
The model satisfies the spatial approximation property if for every $\nu$ with $J_{\rm sp}(\nu)<\infty$ there are real smooth compactly supported controls $v_j$ such that
\begin{equation}\label{eq:recovery}
 q^{v_j}\dd t\dd x\longrightarrow\nu\quad\text{vaguely},\qquad
 \limsup_{j\to\infty}\mathscr S(v_j)\leq J_{\rm sp}(\nu).
\end{equation}
\end{definition}

If smooth controls approximate every finite-action measure with convergence of their costs as in \eqref{eq:recovery}, the preceding bounds extend to a full LDP in the following corollary.

\begin{corollary}[Conditional completion of the LDP]\label{cor:conditionalLDP}
Under the assumptions of Theorem~\ref{thm:spatialbounds}, if Definition~\ref{def:recovery} holds, then $\NN_\eps$ satisfies a full LDP at speed $n_\eps$ with good rate function $J_{\rm sp}$ in the vague topology.
\end{corollary}

\begin{proof}
Let $O\subset\mathsf E$ be open and $\nu\in O$ have finite
$J_{\rm sp}(\nu)$. Definition~\ref{def:recovery} gives
$q^{v_j}\dd t\dd x\in O$ for all large $j$. For every $\eta>0$, it
also gives $\mathscr S(v_j)\leq J_{\rm sp}(\nu)+\eta$ for all large $j$.
Fix such a $j$ before applying \eqref{eq:spatiallower}. Then
\[
\liminf_{\eps\downarrow0}n_\eps^{-1}\log\Pp(\NN_\eps\in O)
 \geq-\mathscr S(v_j)\geq-J_{\rm sp}(\nu)-\eta.
\]
Take $\eta\downarrow0$ and optimize over $\nu\in O$. The upper bound
is \eqref{eq:spatialupper}, so the full LDP follows.

The rate function is lower semicontinuous by \eqref{eq:spatialJ}. To prove
goodness for the rate function, choose a compact $K_A$ from exponential tightness and apply
the just-proved lower bound to its open complement:
\[
 -\inf_{\nu\in\mathsf E\setminus K_A}J_{\rm sp}(\nu)
 \leq\liminf_{\eps\downarrow0}n_\eps^{-1}\log\Pp(\NN_\eps\notin K_A)
 \leq\limsup_{\eps\downarrow0}n_\eps^{-1}\log\Pp(\NN_\eps\notin K_A)
 \leq-A.
\]
Thus $\{J_{\rm sp}\leq M\}\subset K_A$ for $A>M$. This sublevel is
closed by lower semicontinuity and hence compact.
This completes the proof.
\end{proof}

\begin{remark}[The smooth-control approximation assumption]\label{rem:gap}
Under Assumptions~\ref{ass:finite-mean} and \ref{ass:quantitative},
Theorem~\ref{thm:spatialbounds}
gives the closed-set upper bound and the open-set lower bound along
the measures $q^v\dd t\dd x$. Corollary~\ref{cor:conditionalLDP} uses
the approximation property in Definition~\ref{def:recovery} to extend
this lower bound to every finite-action measure.
For finite nonzero initial mass, Theorem~\ref{thm:full-spatial-finite}
establishes the full LDP under Assumption~\ref{ass:finite-mass}
by the approximation lemma and the controlled likelihood argument.
Proposition~\ref{prop:zero-obstruction}
gives a counterexample to approximation by smooth controls when the initial profile vanishes.
\end{remark}

A zero initial profile separates the actual microscopic rate from the variational action and gives an explicit obstruction to 
approximation by smooth controls in the following proposition.

\begin{proposition}[A zero initial profile obstructs the candidate full LDP]\label{prop:zero-obstruction}
The hypotheses of Theorem~\ref{thm:spatialbounds} alone do not imply Definition~\ref{def:recovery}. If $\lambda_0\equiv0$, then $\NN_\eps\equiv0$ almost surely and its actual good rate function, at any diverging speed, is
\begin{equation}\label{eq:zero-true-rate}
 I_0(\nu)=0\quad\text{if }\nu=0,\qquad
 I_0(\nu)=+\infty\quad\text{otherwise}.
\end{equation}
Nevertheless, $J_{\rm sp}$ is finite at some nonzero measures with zero initial trace.
\end{proposition}

\begin{proof}
If $\lambda_0=0$, the prescribed Poisson past has zero intensity and has
no points almost surely. The descendant construction therefore gives
$\NN_\eps=0$ almost surely. For a set containing zero its probability
is one, and for a set excluding zero it is zero. These facts directly
give the LDP with \eqref{eq:zero-true-rate} at every diverging speed.

Let $g$ be a strictly positive bounded smooth integrable density,
$L_g=\int_{\R^d}g(x)\dd x>0$, and $r_t=\cQ_tg$. Set
$q(t,x):=t^2r_t(x)$. The heat equation in \eqref{eq:operators} gives
\[
 m\partial_tq-\B q+\kappa q
 =2mt\,r_t=(2m/t)q\quad(t>0),\qquad q(0,x)=0.
\]
For $v(t,x)=2m/t$, preservation of spatial mass by $\cP_t$ gives
\begin{align*}
 \int_0^T\int_{\R^d}v^2q\dd x\dd t
 &=4m^2\int_0^T\int_{\R^d}r_t(x)\dd x\dd t
 =4m^2L_g\int_0^T e^{-\kappa t/m}\dd t<\infty.
\end{align*}
Its weak residual equals $\int_0^T\int_{\R^d}vq\phi\dd x\dd t$,
since the initial trace is zero. Completing a square bounds the
variational action by half this energy. Conversely, let
$\phi_j(t,x)=\theta_j(t)\chi_j(x)2m/t$, with smooth cutoffs vanishing
near $0,T$ and at spatial infinity and tending to one on $(0,T)\times\R^d$.
Then $\phi_j\in C_c^\infty(\Omega_T)$ and
$\int_0^T\int_{\R^d}q|\phi_j-v|^2\dd x\dd t\to0$ by dominated
convergence theorem. Square completion therefore gives
\begin{equation}\label{eq:zero-counterexample-action}
 J_{\rm sp}(q\dd t\dd x)=2m^2L_g\int_0^T e^{-\kappa t/m}\dd t<\infty.
\end{equation}

For every bounded smooth multiplicative control $w$, its zero-initial
mild solution satisfies
\[
 q_t^w=\frac1m\int_0^t\cQ_{t-s}(w(s)q_s^w)\dd s,\qquad
 \|q_t^w\|_\infty\leq\frac{\|w\|_\infty}{m}
                                  \int_0^t\|q_s^w\|_\infty\dd s.
\]
Iterating this inequality $j$ times gives
$\|q_t^w\|_\infty\leq C_T(\|w\|_\infty T/m)^j/j!\to0$;
thus $q^w=0$. These controlled measures cannot converge vaguely to
the nonzero $q\dd t\dd x$ above.
Choose $f\geq0$ compact with $x=\int_0^T\int_{\R^d}qf\dd x\dd t>0$.
The vague open set $O=\{\nu:|\nu(f)-x|<x/2\}$ excludes zero, so
\[
\liminf_{\eps\downarrow0}n_\eps^{-1}\log\Pp(\NN_\eps\in O)=-\infty
 <-J_{\rm sp}(q\dd t\dd x).
\]
This contradicts a lower bound with the finite candidate action and
proves the obstruction to Definition~\ref{def:recovery}.
\end{proof}

\section{Centering at the PDE and deterministic bias}\label{sec:bias}

Retain $\mu$ and $\mu_\eps$ from Section~\ref{sec:model}, and let $d_\eps=\mu_\eps-\mu$.
Replacing the exact mean by the hydrodynamic measure adds a deterministic shift whose size must be compared with the CLT or MDP normalization.

\begin{proposition}[Exact conditions for changing the center]\label{prop:centering}
The CLT remains valid with $\mu$ in place of $\mu_\eps$ provided
\begin{equation}\label{eq:CLTbias}
\sqrt{n_\eps}\,d_\eps\longrightarrow0\quad\text{in }\cD'(E_T).
\end{equation}
The MDP remains valid with that replacement provided
\begin{equation}\label{eq:MDPbias}
 \frac{\sqrt{n_\eps}}{b_\eps}d_\eps\longrightarrow0
 \quad\text{in }\cD'(E_T).
\end{equation}
If the scaled deterministic bias instead tends to a nonzero distribution $d$, the Gaussian limit is shifted by $d$, or the MDP rate becomes $z\mapsto I_G(z-d)$, respectively.
\end{proposition}

\begin{proof}
Let $h_\eps=\sqrt{n_\eps}d_\eps$ for the CLT and
$h_\eps=(\sqrt{n_\eps}/b_\eps)d_\eps$ for the MDP. The two centered
families satisfy the exact identities
\[
 \sqrt{n_\eps}(\NN_\eps-\mu)=Z_\eps+h_\eps,
 \qquad
 \frac{\sqrt{n_\eps}}{b_\eps}(\NN_\eps-\mu)=Y_\eps+h_\eps,
\]
with the corresponding choice of $h_\eps$ in each identity.
If $h_\eps\to d$ weakly in $\cD'(E_T)$, then for every finite test
family $f_1,\ldots,f_k$,
$h_\eps(f_i)\to d(f_i)$, $1\leq i\leq k$.
Consequently the finite-dimensional CLT shifts by this deterministic
vector. The finite-dimensional MDP cumulants shift by
\[
 b_\eps^{-2}\log\E\left[\exp\left\{b_\eps^2
                    \sum_{i=1}^k\theta_i\left[Y_\eps(f_i)+h_\eps(f_i)\right]\right\}\right]
 \longrightarrow\frac12C(f_\theta,f_\theta)+d(f_\theta),
 \qquad f_\theta=\sum_{i=1}^k\theta_i f_i.
\]
The Legendre transform is the original projected rate
$x-(d(f_1),\ldots,d(f_k))$.

To retain tightness, for a sequence $\eps_j\downarrow0$ let
$H=\{h_{\eps_j}:j\geq1\}\cup\{d\}$, which is compact by convergence.
If $K_A$ is a compact tightness or exponential-tightness set for the
original family, then $K_A+H$ is compact because addition is continuous
in the weak distribution topology. Moreover,
\[
 \Pp\left(Y_{\eps_j}+h_{\eps_j}\notin K_A+H\right)
 \leq\Pp\left(Y_{\eps_j}\notin K_A\right),
\]
and the analogous bound holds for $Z_{\eps_j}$. The finite-cylinder
and extension by finite covers of compact sets in the proof of Theorems~\ref{thm:clt} and
\ref{thm:mdp} now gives the shifted CLT and the rate function
\begin{equation}\label{eq:shifted-centering-limits}
 \begin{aligned}
 \sqrt{n_\eps}(\NN_\eps-\mu)&\Longrightarrow Z+d,\\
 I_{G,d}(z)&:=\sup_{f\in\cD(E_T)}
       \left\{z(f)-d(f)-C(f,f)/2\right\}=I_G(z-d).
 \end{aligned}
\end{equation}
This reasoning applies to every sequence $\eps_j\downarrow0$ and thus
to the limits in \eqref{eq:shifted-centering-limits}.
 Conditions~\eqref{eq:CLTbias} and
\eqref{eq:MDPbias} are precisely the cases $d=0$.
\end{proof}

Proposition~\ref{prop:hydrodynamic-reference} gives convergence of
the exact means. Changing the center at the CLT or MDP scale requires
the rates in \eqref{eq:CLTbias} or \eqref{eq:MDPbias}, respectively.

The additional moments in Assumption~\ref{ass:bias-moments} quantify this deterministic shift and give usable ranges for PDE centering.

\begin{proposition}[A usable bias estimate]\label{prop:biasrate}
Under Assumption~\ref{ass:bias-moments}, for every fixed test
$f\in\cD(E_T)$,
\begin{equation}\label{eq:biasrate}
 |\ip{d_\eps}{f}|\leq C_f\{\eps+|\kappa_\eps-\kappa|\}.
\end{equation}
Under conditional spatial symmetry in Assumption~\ref{ass:bias-moments}, namely
\begin{equation}\label{eq:bias-conditional-symmetry}
 \operatorname{Law}(\Xi\mid S=s)=\operatorname{Law}(-\Xi\mid S=s)
 \quad\text{for }H(s)\dd s\text{-almost every }s,
\end{equation}
equivalently $h(s,z)=h(s,-z)$ almost everywhere, the $\eps$ term in
\eqref{eq:biasrate} can be replaced by $\eps^2$.
Consequently, with $\kappa_\eps=\kappa$, PDE-centered CLTs hold for
$1<\alpha<3$, or for $1<\alpha<5$ under conditional symmetry.
More precisely, let $\beta=1$ in general and $\beta=2$ under
\eqref{eq:bias-conditional-symmetry}. If
$\eps^{-\beta}d_\eps\to d^{(\beta)}$ in $\cD'(E_T)$, then at
$\alpha=1+2\beta$,
\begin{equation}\label{eq:bias-endpoint-scaling}
 \sqrt{n_\eps}\,d_\eps
 =\eps^{-\beta}d_\eps\longrightarrow d^{(\beta)},\qquad
 \sqrt{n_\eps}(\NN_\eps-\mu)\Longrightarrow Z+d^{(\beta)}.
\end{equation}
Thus a coefficient $d^{(\beta)}\ne0$ produces a nonzero Gaussian shift.
Proposition~\ref{prop:biascoefficients} computes these coefficients under
the corresponding expansions of $\kappa_\eps$.
\end{proposition}

\begin{proof}
Let $u:=\cA f$ and write $\|w\|_{0,p}:=\sup_{0\leq t\leq T}\|w(t)\|_p$.
For each spatial multi-index $\beta$, \eqref{eq:resolvent} gives
\begin{equation}\label{eq:backward-derivative-bounds}
 \left\|D_x^\beta u\right\|_{0,p}\leq(T/m)\left\|D_x^\beta f\right\|_{0,p},\quad
 \partial_tu=(\kappa u-\B u-f)/m,\qquad p=1,\infty.
\end{equation}
Iterating the time identity in \eqref{eq:backward-derivative-bounds}
gives, for every integer $r\geq1$ and spatial multi-index $\beta$,
\begin{equation}\label{eq:backward-mixed-derivative-recursion}
 \partial_t^rD_x^\beta u
 =m^{-r}(\kappa-\B)^rD_x^\beta u
 -\sum_{\ell=0}^{r-1}m^{-(r-\ell)}
       (\kappa-\B)^{r-1-\ell}\partial_t^\ell D_x^\beta f.
\end{equation}
All terms have finite $\|\cdot\|_{0,p}$ norms by
\eqref{eq:backward-derivative-bounds}. More directly, extend $f$
smoothly by zero after $T$ and set $g_{r,\beta}:=\partial_t^rD_x^\beta f$.
Changing variables in \eqref{eq:resolvent} gives
$u(t)=m^{-1}\int_0^T\cQ_s f(t+s)\dd s$. Differentiation under
this fixed-interval integral is justified in $L^1$ and in the
uniform norm by the bounded derivatives of $f$ and the contraction
of $\cQ_s$. Hence, for $r\geq0$,
\begin{equation}\label{eq:backward-mixed-derivative-representation}
 \partial_t^rD_x^\beta u(t)
 =\frac1m\int_0^T\cQ_s g_{r,\beta}(t+s)\dd s,
 \qquad
 \left\|\partial_t^rD_x^\beta u\right\|_{0,p}
 \leq\frac Tm\left\|g_{r,\beta}\right\|_{0,p},\quad p=1,\infty.
\end{equation}
Choose $R_0$ containing the spatial support of $f$, hence of each
$g_{r,\beta}$. Let $Y_s$ have the transition law of $\cP_s$ in
\eqref{eq:operators}, so that
$\E[|Y_s|^2]=s\operatorname{tr}(\Sigma)/m$. If $R>R_0$ and
$|y|\leq R_0$, then $|y+Y_s|>R$ implies $|Y_s|>R-R_0$.
Tonelli's theorem, \eqref{eq:backward-mixed-derivative-representation},
and Markov's inequality give
\begin{align}\label{eq:backward-mixed-derivative-tails}
 \sup_{0\leq t\leq T}\int_{\{|x|>R\}}
       \left|\partial_t^rD_x^\beta u(t,x)\right|\dd x
 &\leq\frac{\|g_{r,\beta}\|_{0,1}}m
              \int_0^T\Pp\left(|Y_s|>R-R_0\right)\dd s\nonumber\\
 &\leq\frac{T^2\operatorname{tr}(\Sigma)}{2m^2(R-R_0)^2}
                 \left\|g_{r,\beta}\right\|_{0,1}\longrightarrow0.
\end{align}
This proves the uniform spatial $L^1$ tail bound for every mixed
derivative used in the Taylor expansions below.

The telescoping series in \eqref{eq:K} gives
$\cA_\eps(-\LL_\eps u)=u$. Its last term tends to zero in the weighted
norm of \eqref{eq:renewalbound}. Extending
\eqref{eq:exact-mean-resolvent} to integrable spatial test functions by truncation
therefore gives $-\mu_\eps(\LL_\eps u)=a_\eps \B_\eps u$.
Since $f=-\LL_\kappa u$ and \eqref{eq:mean-duality} gives
$\mu(f)=m\int_{\R^d}\lambda_0(x)u(0,x)\dd x$, we obtain
\begin{equation}\label{eq:biasidentity}
 d_\eps(f)=a_\eps \B_\eps u-m\int_{\R^d}\lambda_0(x)u(0,x)\dd x
                        +\mu_\eps((\LL_\eps-\LL_\kappa)u).
\end{equation}

Fubini's theorem in \eqref{eq:Jpast}, with $S=q-r$, gives
\begin{equation}\label{eq:bias-crossing-representation}
 \B_\eps u=\int_{\R^d}\lambda_0(x)
                 \E\left[\int_0^S u\left(\delta_\eps r,x+\eps\Xi\right)\dd r\right]\dd x.
\end{equation}
The temporal increment and the first spatial increment give
\begin{align}\label{eq:bias-boundary-first-bound}
 \left|\B_\eps u-m\int_{\R^d}\lambda_0(x)u(0,x)\dd x\right|
 &\leq\|\lambda_0\|_\infty\left[
 \frac{\delta_\eps}{2}\E \left[S^2\right]\|\partial_tu\|_{0,1}
+\eps\E\left[S|\Xi|\right]\|\nabla_xu(0,x)\|_1\right].
\end{align}
Here $\E[S|\Xi|]\leq(\E[S^2]\E[|\Xi|^2])^{1/2}<\infty$
by the Cauchy--Schwarz inequality. With conditional symmetry,
$\E[S\Xi]=0$, so Taylor's integral formula gives
\begin{align}\label{eq:bias-boundary-spatial-cancellation}
 &\E\left[S\{u(0,x+\eps\Xi)-u(0,x)\}\right]\nonumber\\
 &\quad=\eps\E[S\Xi]\cdot\nabla u(0,x)
 +\eps^2\E\left[S\int_0^1(1-r)
         \Xi^\top D_x^2u(0,x+r\eps\Xi)\Xi\dd r\right]\nonumber\\
 &\quad=\eps^2\E\left[S\int_0^1(1-r)
         \Xi^\top D_x^2u(0,x+r\eps\Xi)\Xi\dd r\right].
\end{align}
Equation~\eqref{eq:bias-boundary-spatial-cancellation} replaces the
first-order spatial bound
$\eps\E[S|\Xi|]\|\nabla_xu(0)\|_1$ in
\eqref{eq:bias-boundary-first-bound} by the second-order bound
$\delta_\eps\E[S|\Xi|^2]\|D_x^2u(0)\|_1/2$. Together with the
temporal estimate in \eqref{eq:bias-boundary-first-bound}, it gives
\begin{equation}\label{eq:bias-boundary-symmetric-bound}
 \begin{aligned}
 &\left|\B_\eps u-m\int_{\R^d}\lambda_0(x)u(0,x)\dd x\right|\\
 &\quad\leq\frac{\delta_\eps\|\lambda_0\|_\infty}{2}
 \left(\E[S^2]\|\partial_tu\|_{0,1}
              +\E[S|\Xi|^2]\|D_x^2u(0)\|_1\right).
 \end{aligned}
\end{equation}
The mixed moment is finite by \eqref{eq:bias-mixed-moment}.
Also $|1-a_\eps|=\delta_\eps\kappa_\eps$, and
\eqref{eq:Bbounds} bounds $|\B_\eps u|$ uniformly.

For the generator, Taylor's integral formula applied separately to
space and time gives
\begin{align}\label{eq:bias-generator-Taylor-bound}
 \left\|\cK_\eps u-u-\delta_\eps(m\partial_t+\B)u\right\|_{0,1}
 &\leq\frac{\eps^3}{6}\E\left[|\Xi|^3\right]\|D_x^3u\|_{0,1}
+\eps^3\E\left[S|\Xi|\right]\|\partial_t\nabla_xu\|_{0,1}\nonumber\\
 &\quad+\frac{\eps^4}{2}\E\left[S^2\right]\|\partial_t^2u\|_{0,1}
 \leq C_f\eps^3.
\end{align}
The first spatial derivative cancels by $\E[\Xi]=0$; the time-linear
term is expanded once in space. Under conditional symmetry, the odd
spatial terms cancel, and the same computation gives
\begin{equation}\label{eq:bias-generator-symmetric-bound}
 \begin{aligned}
 &\|\cK_\eps u-u-\delta_\eps(m\partial_t+\B)u\|_{0,1}\\
 &\quad\leq\eps^4\left[
 \frac{\E[|\Xi|^4]}{24}\|D_x^4u\|_{0,1}
+\frac{\E[S|\Xi|^2]}2\|\partial_tD_x^2u\|_{0,1}
 +\frac{\E[S^2]}2\|\partial_t^2u\|_{0,1}\right].
 \end{aligned}
\end{equation}
These are finite by Assumption~\ref{ass:bias-moments} and
\eqref{eq:backward-mixed-derivative-representation}.
Using $\LL_\eps u=(\cK_\eps u-u)/\delta_\eps-\kappa_\eps \cK_\eps u$,
we obtain, with $\beta=1$ in general and $\beta=2$ under symmetry,
\begin{align*}
 \left\|(\LL_\eps-\LL_\kappa)u\right\|_{0,1}
 &\leq\delta_\eps^{-1}\left\|\cK_\eps u-u-\delta_\eps(m\partial_t+\B)u\right\|_{0,1}
      +|\kappa_\eps-\kappa|\|u\|_{0,1}
      +\kappa_\eps\|\cK_\eps u-u\|_{0,1}\\
 &\leq C_f\left(\eps^\beta+|\kappa_\eps-\kappa|\right).
\end{align*}
The normalized density bound \eqref{eq:normalized-mean-density-bound} gives
\[
 |\mu_\eps((\LL_\eps-\LL_\kappa)u)|
 \leq T\|\lambda_0\|_\infty\|(\LL_\eps-\LL_\kappa)u\|_{0,1}.
\]
Insert this and \eqref{eq:bias-boundary-first-bound} into
\eqref{eq:biasidentity}. This proves \eqref{eq:biasrate} and its
symmetric improvement, with constants depending on the fixed test function and
the moments in \eqref{eq:basic}, \eqref{eq:bias-moments},
and \eqref{eq:bias-mixed-moment}.
 If $\kappa_\eps=\kappa$,
$\sqrt{n_\eps}\eps^\beta=\eps^{\beta-(\alpha-1)/2}\to0$ precisely
for $1<\alpha<1+2\beta$. Proposition~\ref{prop:centering} gives the
announced PDE-centered CLT ranges.
At $\alpha=1+2\beta$, the identity
$\sqrt{n_\eps}=\eps^{-\beta}$ and Proposition~\ref{prop:centering}
give \eqref{eq:bias-endpoint-scaling} whenever
$\eps^{-\beta}d_\eps\to d^{(\beta)}$.
\end{proof}

\begin{remark}[Bias condition for MDP centering]
For an error $\mathcal O(\eps^\beta+|\kappa_\eps-\kappa|)$, the sufficient PDE-centering condition is
\[
 \frac{\eps^{-(\alpha-1)/2}}{b_\eps}
 \left\{\eps^\beta+|\kappa_\eps-\kappa|\right\}\longrightarrow0.
\]
Together with $b_\eps\to\infty$ and $b_\eps/\sqrt{n_\eps}\to0$,
this condition implies \eqref{eq:MDPbias} by
Proposition~\ref{prop:biasrate}.
\end{remark}

At the first two density thresholds, the leading weak bias can survive in the CLT; the next proposition computes its coefficients and the resulting Gaussian shifts.

\begin{proposition}[Weak bias coefficients and endpoint CLTs]\label{prop:biascoefficients}
Assume Assumption~\ref{ass:bias-moments} and write
\[
 c_{11}:=\E[S\Xi],\quad M_3:=\E\left[\Xi^{\otimes3}\right],\quad
 s_2:=\E\left[S^2\right],\quad M_{12}:=\E\left[S\Xi\Xi^\top\right],\quad M_4:=\E\left[\Xi^{\otimes4}\right].
\]
Tensor contractions below are componentwise. Let $\mathcal T:=m\partial_t+\B$ and
\begin{align}\label{eq:bias-operators}
 \mathcal C_3&:=c_{11}\cdot\partial_t\nabla_x+\frac16M_3:D_x^3,\\
 \mathcal C_4&:=\frac12s_2\partial_t^2
       +\frac12M_{12}:\partial_tD_x^2+\frac1{24}M_4:D_x^4.\nonumber
\end{align}
If $\kappa_\eps=\kappa+\eps\kappa_1+o(\eps)$, then
$\eps^{-1}(\mu_\eps-\mu)\to d^{(1)}$ in $\cD'(E_T)$, where, with $u=\cA f$,
\begin{equation}\label{eq:first-bias}
 d^{(1)}(f):=\int_{\R^d}\lambda_0(x)c_{11}\cdot\nabla_xu(0,x)\dd x
      +\int_0^T\int_{\R^d}\rho(t,x)(\mathcal C_3u-\kappa_1u)(t,x)\dd x\dd t.
\end{equation}
If instead $h(s,z)=h(s,-z)$ and
$\kappa_\eps=\kappa+\eps^2\kappa_2+o(\eps^2)$, then
$\eps^{-2}(\mu_\eps-\mu)\to d^{(2)}$, where
\begin{align}\label{eq:second-bias}
 d^{(2)}(f)
 &:=\int_{\R^d}\lambda_0(x)\left[\frac12s_2\partial_tu(0,x)
             +\frac12M_{12}:D_x^2u(0,x)-\kappa m u(0,x)\right]\dd x\nonumber\\
 &\quad+\int_0^T\int_{\R^d}\rho(t,x)
       [\mathcal C_4u-\kappa\mathcal T u-\kappa_2u](t,x)\dd x\dd t.
\end{align}
Consequently, at $\alpha=3$ in the first case, or $\alpha=5$ in the second case, respectively,
\begin{equation}\label{eq:endpoint-bias-clt}
 \sqrt{n_\eps}(\NN_\eps-\mu)\Longrightarrow Z+d^{(r)},\qquad r=1\ \text{or}\ 2.
\end{equation}
Centering instead at $\mu+\eps^r d^{(r)}$ removes this endpoint shift. At either endpoint the PDE-centered MDP already has the unshifted rate function $I_G$ for every admissible $b_\eps\to\infty$.
\end{proposition}

\begin{proof}
Let $u=\cA f$ and use the derivative bounds
\eqref{eq:backward-derivative-bounds}. Its spatial and temporal
derivatives, extended by zero beyond $T$, have a common $L^1$ modulus
of translation continuity on the time interval. Denote such bounded
moduli by $\omega_{a,b}(r)$ for $\partial_t^aD_x^bu$; each tends to
zero as $r\downarrow0$.
Taylor's integral formula in time gives the exact decomposition
\begin{align*}
 \cK_\eps u(t,x)
 &=\E \left[u(t,x+\eps\Xi)\right]+\eps^2\E[S\partial_tu(t,x+\eps\Xi)]\\
 &\quad+\eps^4\E\left[S^2\int_0^1(1-r)
           \partial_t^2u(t+r\eps^2S,x+\eps\Xi)\dd r\right].
\end{align*}
Replacing the last integral by $\partial_t^2u(t,x)/2$ has error bounded by
\[
\eps^4\E\left[S^2\omega_{2,0}\left(\eps^2S+\eps|\Xi|\right)\right]=o(\eps^4)
 \quad\text{in }\|\cdot\|_{0,1}.
\]
This follows by the dominated convergence theorem with weight $S^2$.
The pure spatial term has expansion
\begin{align*}
 \E\left[u(t,x+\eps\Xi)\right]
 &=u+\frac{\eps^2}{2}\Sigma:D_x^2u
       +\frac{\eps^3}{6}M_3:D_x^3u
       +\frac{\eps^4}{24}M_4:D_x^4u+R_{\eps,1},\\
 \|R_{\eps,1}\|_{0,1}
 &\leq C\eps^4\E\left[|\Xi|^4\omega_{0,4}(\eps|\Xi|)\right]=o(\eps^4).
\end{align*}
For the mixed term,
\begin{align*}
 \eps^2\E[S\partial_tu(t,x+\eps\Xi)]
 &=\eps^2m\partial_tu+\eps^3c_{11}\cdot\partial_t\nabla_xu
       +\frac{\eps^4}{2}M_{12}:\partial_tD_x^2u+R_{\eps,2},\\
 \|R_{\eps,2}\|_{0,1}
 &\leq C\eps^4\E\left[S|\Xi|^2\omega_{1,2}(\eps|\Xi|)\right]=o(\eps^4).
\end{align*}
The weights $S^2$, $|\Xi|^4$, and $S|\Xi|^2$ are integrable by
\eqref{eq:bias-moments} and \eqref{eq:bias-mixed-moment}.
Thus the definitions in \eqref{eq:bias-operators} give
\begin{equation}\label{eq:finite-moment-generator-expansion}
 \cK_\eps u=u+\eps^2\mathcal T u+\eps^3\mathcal C_3u
                            +\eps^4\mathcal C_4u+o(\eps^4)
 \quad\text{in }\|\cdot\|_{0,1}.
\end{equation}
Multiplying by $a_\eps=1-\eps^2\kappa_\eps$ and dividing by $\eps^2$ gives
\begin{equation}\label{eq:bias-controlled-generator-expansion}
 (\LL_\eps-\LL_\kappa)u
 =\eps\mathcal C_3u-(\kappa_\eps-\kappa)u
       +\eps^2\left(\mathcal C_4u-\kappa_\eps\mathcal T u\right)+o(\eps^2).
\end{equation}

For the boundary term, use \eqref{eq:bias-crossing-representation}.
The spatial second-order expansion weighted by $S$, and the first-order
temporal expansion integrated over $0\leq r\leq S$, give
\begin{equation}\label{eq:bias-boundary-coefficient-expansion}
 \B_\eps u=\beta_0+\eps\beta_1+\eps^2\beta_2+o(\eps^2),
\end{equation}
where
\begin{align*}
\beta_0&:=m\int_{\R^d}\lambda_0(x)u(0,x)\dd x,\\
 \beta_1&:=\int_{\R^d}\lambda_0(x)c_{11}\cdot\nabla_xu(0,x)\dd x,\\
 \beta_2&:=\int_{\R^d}\lambda_0(x)
       \left[s_2\partial_tu(0,x)/2+M_{12}:D_x^2u(0,x)/2\right]\dd x.
\end{align*}
Indeed $\int_0^S r\dd r=S^2/2$, and the two boundary remainder bounds are
\[
 C\|\lambda_0\|_\infty\eps^2
 \left\{\E\left[S|\Xi|^2\omega_{0,2}(\eps|\Xi|)\right]
        +\E\left[S^2\omega_{1,0}\left(\eps^2S+\eps|\Xi|\right)\right]\right\}=o(\eps^2).
\]
These estimates keep every derivative on $u$.

If $\kappa_\eps=\kappa+\eps\kappa_1+o(\eps)$,
\eqref{eq:bias-controlled-generator-expansion} and
\eqref{eq:bias-boundary-coefficient-expansion} give
\[
 \eps^{-1}(\LL_\eps-\LL_\kappa)u\to\mathcal C_3u-\kappa_1u,
 \qquad \eps^{-1}(a_\eps \B_\eps u-\beta_0)\to\beta_1.
\]
Under conditional symmetry, $c_{11}=M_3=0$. If
$\kappa_\eps=\kappa+\eps^2\kappa_2+o(\eps^2)$, the corresponding limits are
\[
 \eps^{-2}(\LL_\eps-\LL_\kappa)u
 \to\mathcal C_4u-\kappa\mathcal T u-\kappa_2u,
 \qquad \eps^{-2}(a_\eps \B_\eps u-\beta_0)\to\beta_2-\kappa\beta_0.
\]
The generator limits hold in $\|\cdot\|_{0,1}$.
To replace $\mu_\eps$ by $\mu$ in these limits, let $g_\eps\to g$ in
that norm and let $\chi_R$ be a spatial cutoff. By
\eqref{eq:normalized-mean-density-bound},
\begin{align*}
 |\mu_\eps(g_\eps)-\mu(g)|
 &\leq T\|\lambda_0\|_\infty\|g_\eps-g\|_{0,1}
       +|(\mu_\eps-\mu)(\chi_R g)|\\
 &\quad+2\|\lambda_0\|_\infty
               \int_0^T\int_{\R^d}(1-\chi_R)|g|\dd x\dd t.
\end{align*}
The middle term tends to zero by Proposition~\ref{prop:hydrodynamic-reference};
the last tends to zero as $R\uparrow\infty$ by the derivative tails.
Insertion into \eqref{eq:biasidentity} proves
\eqref{eq:first-bias} and \eqref{eq:second-bias}.

At $\alpha=1+2r$, $\sqrt{n_\eps}=\eps^{-r}$. Therefore
\[
 \sqrt{n_\eps}(\NN_\eps-\mu)
 =Z_\eps+\eps^{-r}d_\eps\Longrightarrow Z+d^{(r)},
 \qquad r=1\text{ or }2,
\]
by Theorem~\ref{thm:clt} and Proposition~\ref{prop:centering}.
Replacing the center $\mu$ by $\mu+\eps^r d^{(r)}$ removes the limiting shift.
For the MDP, $\eps^{-r}d_\eps/b_\eps\to0$ since $b_\eps\to\infty$;
Proposition~\ref{prop:centering} gives the unshifted rate function at these endpoints.
\end{proof}

\section{Heavy-tailed delays: fluctuations, moderate and large deviations}\label{sec:fractional}
This section studies the factorized heavy-tail model in
Assumption~\ref{ass:heavy}, with $\lambda_0\geq0$ bounded and uniformly
continuous and the state spaces of Section~\ref{subsec:spaces}.
The prescribed Poisson past and scaling are specified in \eqref{eq:heavyfamily}
and \eqref{eq:heavyscales}, following \cite{GWZ2026}.
Theorems~\ref{thm:fractional} and \ref{thm:fractional-mdp} establish
the space--time CLT and MDP, respectively. Under the additional condition $0<L<\infty$,
Theorem~\ref{thm:fractional-spatial-ldp} establishes a full LDP for
the space--time measure under the same microscopic scaling.

\begin{assumption}[Factorized heavy-tail high-density model]\label{ass:heavy}
Let $\lambda_0\geq0$ be bounded and uniformly continuous and $a_\eps\in(0,1]$. Assume $h(s,z)=H(s)p(z)$, where $H$ and $p$ are probability densities on $(0,\infty)$ and $\mathbb{R}^{d}$, respectively, $p$ is centered with finite positive definite covariance $\Sigma$, and
\begin{equation}\label{eq:heavyassumptions}
 \overline H(s)\sim c_Hs^{-\eta},\qquad 0<\eta<1,\qquad
 c_\eta=c_H\Gamma(1-\eta),\qquad
 \frac{1-a_\eps}{\eps^{2\eta}}\longrightarrow\kappa\geq0.
\end{equation}
The process $N^{\eps,\eta}$ follows the offspring construction and
positive-time intensity equation \eqref{eq:intensity}, with a Poisson
past $\Pi^{\eps,\eta}$ and the factorized kernel of this assumption.
Every positive-time event descends from a prescribed parent.
The intensity of the prescribed Poisson past and normalized event measure are
\begin{align}\label{eq:heavyfamily}
 q_\eps\lambda_0(\eps^\eta y)\dd r\dd y,
 &\qquad q_\eps=\eps^{2-\eta(\alpha+3-d)},\quad r<0,\nonumber\\
 \ip{\NN_\eps^\eta}{f}
 &=\eps^{\eta(\alpha+3)}\int_0^\infty\int_{\R^d}f\left(\eps^2t,\eps^\eta x\right)N^{\eps,\eta}(\dd t,\dd x).
\end{align}
Set
\begin{equation}\label{eq:heavyscales}
\tau_\eps:=\eps^2,\quad\delta_{\eps,\eta}:=\eps^{2\eta},\quad
 n_{\eps,\eta}:=\eps^{-\eta(\alpha-1)},\quad\alpha>1.
\end{equation}
\end{assumption}
Thus the mass per event is $\delta_{\eps,\eta}^2/n_{\eps,\eta}$.
Here $\tau_\eps$ is the time scale, $\delta_{\eps,\eta}$ is the generator
scale, and $n_{\eps,\eta}$ is the effective population.

The delay density $H$, its tail $\overline H$, and its Laplace transform retain their definitions in \eqref{eq:delay-marginal}; now $\E[S]=\infty$. The offspring displacement is still denoted by $\Xi$, whereas $Z_\eta$ below denotes a Gaussian field.

Let $S_r^\eta$ be a subordinator with Laplace transform
$\E\left[e^{-zS_r^\eta}\right]=e^{-r c_\eta z^\eta}$, independent of a Brownian motion $B_r^\Sigma$ of covariance $r\Sigma$. Define
\begin{equation}\label{eq:fractionalA}
 \cA_\eta f(t,x):=\int_0^\infty e^{-\kappa r}
\E\left[f\left(t+S_r^\eta,x+B_r^\Sigma\right)\right]\dd r,
\end{equation}
with zero extension beyond $T$. In particular,
$\norm{\cA_\eta}_{q,p}\leq(\kappa+c_\eta q^\eta)^{-1}$ for $q>0$ and $p=1,\infty$. The mean $\rho^\eta$ is characterized by
\begin{equation}\label{eq:heavymean}
 c_\eta\,{}^{\rm C}D_t^\eta\rho^\eta=\B\rho^\eta-\kappa\rho^\eta,
 \qquad\rho^\eta(0)=\lambda_0,
\end{equation}
where, for an absolutely continuous scalar function $r$,
\[
{}^{\rm C}D_t^\eta r:=\frac{1}{\Gamma(1-\eta)}\int_0^t(t-s)^{-\eta}r'(s)\dd s.
\]
For the bounded density $\rho^\eta$, equation \eqref{eq:heavymean}
denotes the resolvent solution given by the inverse-stable-clock
representation in \eqref{eq:fractional-inverse-clock-mean}. The connection
between inverse stable subordinators and fractional equations is reviewed
in \cite{MS}; related heavy-tail Hawkes limits are in \cite{JRrough}.
Set
\begin{equation}\label{eq:heavycov}
C_\eta(f,g):=\int_0^T\int_{\R^d}\rho^\eta(t,x)
                 \cA_\eta f(t,x)\cA_\eta g(t,x)\dd x\dd t.
\end{equation}

The inverse stable clock identifies the fractional hydrodynamic density, while the duality in the following lemma 
links that density to the memory-weighted functional of the prescribed past.

\begin{lemma}[Fractional mean and backward duality]
\label{lem:fractional-mean-duality}
Let $E_t^\eta:=\inf\{r\geq0:S_r^\eta>t\}$, independent of the Brownian
motion used in \eqref{eq:fractionalA}. Then
\begin{equation}\label{eq:fractional-inverse-clock-mean}
 \rho^\eta(t,x)=\E\left[e^{-\kappa E_t^\eta}
                \lambda_0\left(x+B_{E_t^\eta}^\Sigma\right)\right]
\end{equation}
is bounded, nonnegative, and tends uniformly to $\lambda_0$ as $t\downarrow0$.
For each continuous compactly supported $w$,
\begin{equation}\label{eq:fractional-duality-physical}
 \int_0^T\int_{\R^d}\rho^\eta(t,x)w(t,x)\dd x\dd t
 =c_H\int_0^T t^{-\eta}\int_{\R^d}\lambda_0(x)\cA_\eta w(t,x)\dd x\dd t.
\end{equation}
These assertions hold for the bounded uniformly continuous initial
density in Assumption~\ref{ass:heavy}.
\end{lemma}

\begin{proof}
Equation~\eqref{eq:fractional-inverse-clock-mean} gives
$0\leq\rho^\eta\leq\|\lambda_0\|_\infty$.
Let $\omega_0$ be the modulus of uniform continuity of $\lambda_0$.
For $a,h>0$, conditioning on $E_t^\eta$ and using the Brownian second
moment give
\begin{align*}
\Pp\left(\left|B_{E_t^\eta}^\Sigma\right|>a\right)
&\leq\Pp\left(E_t^\eta>h\right)+h\operatorname{tr}(\Sigma)/a^2,\\
 \left\|\rho^\eta(t)-\lambda_0\right\|_\infty
&\leq\omega_0(a)+3\|\lambda_0\|_\infty\Pp\left(E_t^\eta>h\right)
       +2\|\lambda_0\|_\infty h\operatorname{tr}(\Sigma)/a^2
       +\kappa h\|\lambda_0\|_\infty.
\end{align*}
For fixed $h>0$, $\Pp(E_t^\eta>h)=\Pp(S_h^\eta\leq t)\to0$:
$\Pp(S_h^\eta=0)=\lim_{z\to\infty}e^{-hc_\eta z^\eta}=0$.
Taking $t\downarrow0$, then $h\downarrow0$, then $a\downarrow0$
proves the uniform initial limit.

For $z>0$ and $r\geq0$, Tonelli's theorem gives
\[
 \int_0^\infty e^{-zt}\Pp\left(E_t^\eta>r\right)\dd t
 =\E\left[\int_{S_r^\eta}^\infty e^{-zt}\dd t\right]
 =z^{-1}e^{-rc_\eta z^\eta}.
\]
Differentiating its right side as a tail of a measure in $r$ gives
\begin{equation}\label{eq:inverse-clock-Laplace-law}
 \int_0^\infty e^{-zt}\Pp(E_t^\eta\in\dd r)\dd t
 =c_\eta z^{\eta-1}e^{-rc_\eta z^\eta}\dd r.
\end{equation}
Let $\cP_r^{\Sigma}=\exp(r\B)$ be the heat semigroup with covariance
$r\Sigma$, and define the operator-valued potential
\[
 \mathcal U_\eta(\dd u)g
 =\int_0^\infty e^{-\kappa r}\Pp(S_r^\eta\in\dd u)\cP_r^{\Sigma}g\dd r.
\]
Its Laplace transform is
$\widehat{\mathcal U}_\eta(z)=(\kappa+c_\eta z^\eta-\B)^{-1}$,
as follows by integration of the semigroup. By
\eqref{eq:inverse-clock-Laplace-law},
\begin{equation}\label{eq:fractional-mean-Laplace-resolvent}
 \int_0^\infty e^{-zt}\rho^\eta(t)\dd t
 =c_\eta z^{\eta-1}\widehat{\mathcal U}_\eta(z)\lambda_0.
\end{equation}
The kernel $b(t)=c_Ht^{-\eta}$ has Laplace transform
$c_\eta z^{\eta-1}$. The zero-killing temporal marginal of
$\mathcal U_\eta$ has transform $1/(c_\eta z^\eta)$ and hence density
$u^{\eta-1}/(c_\eta\Gamma(\eta))$.
Thus its convolution with $b$ is locally integrable, because
\[
 \frac{c_H}{c_\eta\Gamma(\eta)}
 \int_0^t(t-u)^{-\eta}u^{\eta-1}\dd u=1.
\]
Uniqueness of the Laplace transform in
\eqref{eq:fractional-mean-Laplace-resolvent} gives the operator identity
\begin{equation}\label{eq:fractional-mean-potential-convolution}
 \rho^\eta(t)=\int_0^t b(t-u)\mathcal U_\eta(\dd u)\lambda_0.
\end{equation}

For nonnegative compact $w$, use Tonelli's theorem in
\eqref{eq:fractional-mean-potential-convolution}, set $s:=t+u$, and use
symmetry of $\cP_r^{\Sigma}$. This gives
\begin{align*}
 \int_0^T\int_{\R^d}\rho^\eta(s,x)w(s,x)\dd x\dd s
 &=\int_0^T b(t)\int_{\R^d}\lambda_0(x)
       \left[\int_0^{T-t}\mathcal U_\eta(\dd u)w(t+u)\right](x)\dd x\dd t\\
 &=c_H\int_0^Tt^{-\eta}\int_{\R^d}\lambda_0(x)\cA_\eta w(t,x)\dd x\dd t,
\end{align*}
where the last identity is \eqref{eq:fractionalA}. To justify the
same calculation for signed $w$, apply positivity and the
finite-horizon bound for \eqref{eq:fractionalA} to $|w|$:
\begin{align}\label{eq:fractional-duality-absolute-bound}
 c_H\int_0^Tt^{-\eta}\int_{\R^d}\lambda_0(x)
                    \cA_\eta|w|(t,x)\dd x\dd t
 &\leq\frac{c_HT^{1-\eta}}{1-\eta}
       \|\lambda_0\|_\infty\left\|\cA_\eta|w|\right\|_{0,1}\nonumber\\
 &\leq\frac{c_HT^{1-\eta}e^{zT}}
                   {(1-\eta)(\kappa+c_\eta z^\eta)}
       \|\lambda_0\|_\infty\|w\|_{0,1}<\infty,\qquad z>0.
\end{align}
Therefore, Fubini's theorem extends the identity to signed $w$, proving
\eqref{eq:fractional-duality-physical}. Finally
\eqref{eq:fractional-mean-Laplace-resolvent} gives
\[
 c_\eta\left[z^\eta\widehat\rho^\eta(z)-z^{\eta-1}\lambda_0\right]
 =(\B-\kappa)\widehat\rho^\eta(z),
\]
which verifies the resolvent interpretation of \eqref{eq:heavymean}.
The proof is complete.
\end{proof}

\subsection{The fractional CLT and MDP}

The fractional resolvent and mean duality replace their finite-mean counterparts in both the Gaussian covariance and the quadratic moderate-deviation rate.

\begin{theorem}[Fractional fluctuations]\label{thm:fractional}
Under Assumption~\ref{ass:heavy}, let
$\mu_\eps^\eta=\E\NN_\eps^\eta$. Then
\[
Z_{\eps,\eta}:=\sqrt{n_{\eps,\eta}}\left(\NN_\eps^\eta-\mu_\eps^\eta\right)
 \Longrightarrow Z_\eta\quad\text{in }\cD'(E_T),
\]
where
\begin{equation}\label{eq:fractionalGaussian}
Z_\eta(f)=\int_0^T\int_{\R^d}\sqrt{\rho^\eta(t,x)}\,\cA_\eta f(t,x)W(\dd t,\dd x)
\end{equation}
is a centered Gaussian field with covariance \eqref{eq:heavycov}.
The convergence holds for every $0<\eta<1$ in the weak distribution
topology on $\cD'(E_T)$.
\end{theorem}

At the moderate scale, the covariance in \eqref{eq:heavycov}
determines the quadratic control energy of the fluctuation field.

\begin{theorem}[Fractional moderate deviations]\label{thm:fractional-mdp}
Under Assumption~\ref{ass:heavy}, let $\mu_\eps^\eta=\E\NN_\eps^\eta$.
For every $b_\eps\to\infty$ with
$b_\eps/\sqrt{n_{\eps,\eta}}\to0$, the distributions
\[
 Y_{\eps,\eta}:=\frac{\sqrt{n_{\eps,\eta}}}{b_\eps}
                 (\NN_\eps^\eta-\mu_\eps^\eta)
\]
obey an MDP at speed $b_\eps^2$, with good rate function
\begin{equation}\label{eq:fractionalrate}
 I_\eta(z)=\inf_{\zeta\in L^2(E_T)}\left\{
 \begin{aligned}
 &\frac12\norm{\zeta}_{L^2(E_T)}^2:\\
 &z(f)=\int_0^T\int_{\R^d}\sqrt{\rho^\eta(t,x)}\cA_\eta f(t,x)\zeta(t,x)\dd x\dd t,\\
 &\text{for every }f\in\cD(E_T)
 \end{aligned}\right\}.
\end{equation}
The MDP holds for every $0<\eta<1$ in the weak
distribution topology on $\cD'(E_T)$.
\end{theorem}

\begin{proof}[Proof of Theorem~\ref{thm:fractional}]
Write $\delta=\delta_{\eps,\eta}=\tau_\eps^\eta$ in this proof.
The tail assumption \eqref{eq:heavyassumptions} gives, for $t>0$,
\begin{equation}\label{eq:fractional-preparation-tail-bound}
 b_\eps(t):=\delta^{-1}\overline H(t/\tau_\eps)
 \leq Ct^{-\eta},\qquad b_\eps(t)\to c_Ht^{-\eta}.
\end{equation}
To see the uniform bound, use $\overline H(s)\leq Cs^{-\eta}$ for
large $s$ and $\overline H(s)\leq1$ on a bounded interval. Integration
by parts in the temporal Laplace transform gives
\begin{align}\label{eq:fractional-delay-transform-limit}
 \frac{1-\widehat H(q\tau_\eps)}\delta
 &=q\int_0^\infty e^{-qt}b_\eps(t)\dd t
 \longrightarrow qc_H\int_0^\infty e^{-qt}t^{-\eta}\dd t
 =c_\eta q^\eta,\qquad q>0.
\end{align}
Domination in this limit follows from \eqref{eq:fractional-preparation-tail-bound}.
Let $\widehat p(\xi):=\E_p\left[ e^{\mathrm{i}\xi\cdot\Xi}\right]$. Centering and the finite
spatial covariance give
$\widehat p(\sqrt\delta\xi)=1-\delta\xi^\top\Sigma\xi/2+o(\delta)$.
Factorization therefore gives
\[
 \E\left[e^{-q\tau_\eps S+\mathrm{i}\sqrt\delta\xi\cdot\Xi}\right]
 =1-\delta\left(c_\eta q^\eta+\xi^\top\Sigma\xi/2\right)+o(\delta).
\]

Let $\cK_{\eps,\eta}f(t,x):=\E\left[f\left(t+\tau_\eps S,x+\sqrt\delta\Xi\right)\right]$ and
$\cA_{\eps,\eta}:=\delta\sum_{j=0}^\infty a_\eps^j\cK_{\eps,\eta}^j$.
For independent pairs $(S_\ell,\Xi_\ell)$ with the offspring law, let
$T_j:=\sum_{\ell=1}^jS_\ell$, $X_j:=\sum_{\ell=1}^j\Xi_\ell$,
and $T_0:=X_0:=0$. The renewal measure
$U_{\eps,\eta}:=\delta\sum_{j=0}^\infty a_\eps^j
\operatorname{Law}(\tau_\eps T_j,\sqrt\delta X_j)$ has Laplace--Fourier transform
\begin{equation}\label{eq:fractional-renewal-transform}
 \frac\delta{1-a_\eps\widehat H(q\tau_\eps)\widehat p(\sqrt\delta\xi)}
 \longrightarrow\frac1{\kappa+c_\eta q^\eta+\xi^\top\Sigma\xi/2},
\end{equation}
and weighted operator bounds
\begin{equation}\label{eq:fractional-renewal-operator-bound}
 \|\cA_{\eps,\eta}\|_{q,p}
 \leq\frac{\delta}{1-a_\eps\widehat H(q\tau_\eps)}
 \longrightarrow\frac1{\kappa+c_\eta q^\eta},
 \qquad
 \|\cA_\eta\|_{q,p}\leq\frac1{\kappa+c_\eta q^\eta},
 \quad p=1,\infty.
\end{equation}
Here the limit follows from \eqref{eq:fractional-delay-transform-limit}
and \eqref{eq:heavyassumptions}, while the last inequality follows
from \eqref{eq:fractionalA}. Fix $q>0$ and a sufficiently small
$\eps_0>0$. Passing from weighted to unweighted norms gives the
common bound
\begin{equation}\label{eq:fractional-common-resolvent-bound}
 \begin{aligned}
 \max\{\|\cA_{\eps,\eta}\|_{0,p},\|\cA_\eta\|_{0,p}\}
 &\leq M_T,\qquad 0<\eps\leq\eps_0,\quad p=1,\infty,\\
 M_T&:=e^{qT}\max\left\{
       \sup_{0<\eps\leq\eps_0}
             \frac{\delta_{\eps,\eta}}{1-a_\eps\widehat H(q\tau_\eps)},
       \frac1{\kappa+c_\eta q^\eta}\right\}<\infty.
 \end{aligned}
\end{equation}
The limiting potential measure is
\[
 U_\eta(\dd u,\dd z)=\int_0^\infty e^{-\kappa r}
               \operatorname{Law}\left(S_r^\eta,B_r^\Sigma\right)(\dd u,\dd z)\dd r,
\]
whose convolution is \eqref{eq:fractionalA}.
For $p_\eps=\widehat H(q\tau_\eps)$, independence of delays and
displacements gives the spatial moment bound
\[
 \int_0^\infty\int_{\R^d}e^{-qu}|z|^2U_{\eps,\eta}(\dd u,\dd z)
=\delta^2\E\left[|\Xi|^2\right]\sum_{j=1}^\infty ja_\eps^jp_\eps^j
 =\frac{\delta^2a_\eps p_\eps\E[|\Xi|^2]}{(1-a_\eps p_\eps)^2}\leq C_q.
\]
For $q,L>0$, the temporal tail satisfies
\begin{align}\label{eq:fractional-renewal-temporal-tail}
 \int_L^\infty\int_{\R^d}e^{-qu}U_{\eps,\eta}(\dd u,\dd z)
 &\leq e^{-qL/2}\int_0^\infty\int_{\R^d}
                  e^{-qu/2}U_{\eps,\eta}(\dd u,\dd z)\nonumber\\
 &=e^{-qL/2}\frac\delta{1-a_\eps\widehat H(q\tau_\eps/2)}.
\end{align}
These bounds give tightness of $e^{-qu}U_{\eps,\eta}$;
\eqref{eq:fractional-renewal-transform} identifies every weak limit
with $e^{-qu}U_\eta$. The finite uniform-net convolution argument in
Lemma~\ref{lem:renewal} then gives uniform convergence on continuous
sources vanishing at $T$. For compact sources its spatial moment bound
also gives
\[
 \left\|\cA_{\eps,\eta}f-\cA_\eta f\right\|_{0,1}
 \leq |B_R|\left\|\cA_{\eps,\eta}f-\cA_\eta f\right\|_\infty+C_{f,T}R^{-2}\to0,
\]
by taking $\eps\downarrow0$ and then $R\uparrow\infty$.

The offspring equation \eqref{eq:positivecluster} uses $\delta$ as
the rescaling of $F-1$. Its nonlinearity satisfies
$G_\delta'(z)=z/(1+\delta z)$.
For $R_f:=2M_T\|f\|_\infty$, choose the source radius so that
$2M_TR_f\leq1/2$ and $\delta R_f\leq1/2$.
The bounds on $G_\delta$ used in \eqref{eq:cluster-contraction}
and the two operator norms in
\eqref{eq:fractional-common-resolvent-bound} give
\begin{align}\label{eq:fractional-cluster-contraction}
 \left\|\cA_{\eps,\eta}(f+G_\delta(w))\right\|_\infty
 &\leq M_T\|f\|_\infty+M_TR_f^2\leq R_f,\nonumber\\
 \left\|\cA_{\eps,\eta}(G_\delta(w)-G_\delta(\widetilde w))\right\|_{0,p}
 &\leq2M_TR_f\|w-\widetilde w\|_{0,p}
 \leq\frac12\|w-\widetilde w\|_{0,p},
\end{align}
for $\|w\|_\infty,\|\widetilde w\|_\infty\leq R_f$ and $p=1,\infty$.
Thus the same fixed-point construction gives
\begin{equation}\label{eq:fractional-cluster-norm-bounds}
\left\|v_{\eps,\eta,f}\right\|_\infty\leq2M_T\|f\|_\infty,\qquad
\left\|v_{\eps,\eta,f}\right\|_{0,1}\leq2M_T\|f\|_{0,1},
\end{equation}
and these solutions converge in both norms to
$v_{\eta,f}=\cA_\eta(f+v_{\eta,f}^2/2)$.
The finite-generation domination in the proof of Lemma~\ref{lem:analytic}
identifies them with positive and signed cluster moments;
its proof uses the common resolvent bounds and hence applies to these scales.

Integrating the factorized offspring law against the past gives
\begin{equation}\label{eq:heavyB}
 \B_{\eps,\eta}w=\int_0^T b_\eps(t)\int_{\R^d}\lambda_0(x)
          \E_p\left[w\left(t,x+\eps^\eta\Xi\right)\right]\dd x\dd t.
\end{equation}
Thus
$|\B_{\eps,\eta}w|\leq C_T\|\lambda_0\|_\infty\|w\|_{0,1}$.
If $w_\eps\to w$ in this norm and $w$ is continuous in $L^1$,
\begin{align*}
 \left|\B_{\eps,\eta}w_\eps-\B_\eta w\right|
 &\leq C_T\|\lambda_0\|_\infty\|w_\eps-w\|_{0,1}\\
 &\quad+\|\lambda_0\|_\infty\int_0^T b_\eps(t)
                    \E_p\left[\left\|w\left(t,\cdot+\eps^\eta\Xi\right)-w(t)\right\|_1\right]\dd t\\
 &\quad\quad+\|\lambda_0\|_\infty\int_0^T
              \left|b_\eps(t)-c_Ht^{-\eta}\right|\|w(t)\|_1\dd t\longrightarrow0,
\end{align*}
where $\B_\eta w:=c_H\int_0^Tt^{-\eta}\int_{\R^d}\lambda_0w\dd x\dd t$.
The last two limits follow from translation continuity and the
integrable dominating density in \eqref{eq:fractional-preparation-tail-bound}.

Define $\cJ_{\eps,\eta}$ by \eqref{eq:Jpast} with time scale $\tau_\eps$
and spatial scale $\eps^\eta$. The unit-population parent density in
macroscopic space is $\tau_\eps\delta^{-2}\lambda_0$. Hence,
\begin{equation}\label{eq:heavy-parent-cgf}
 \Lambda_{\eps,\eta}(f)=\frac{\tau_\eps}{\delta^2}
 \int_{-\infty}^0\int_{\R^d}\lambda_0(x)
       \left[e^{a_\eps\delta \cJ_{\eps,\eta}v_{\eps,\eta,f}(r,x)}-1\right]\dd x\dd r.
\end{equation}
The identity
$\int_{-\infty}^0\int_{\R^d}\lambda_0\cJ_{\eps,\eta}w\dd x\dd r
=(\delta/\tau_\eps)\B_{\eps,\eta}w$ and
$|\cJ_{\eps,\eta}w|\leq\|w\|_\infty$ give the explicit remainder bound
\begin{align}\label{eq:fractional-parent-transform-remainder}
 \left|\Lambda_{\eps,\eta}(f)-a_\eps \B_{\eps,\eta}v_{\eps,\eta,f}\right|
 &\leq\frac{\tau_\eps}{2}e^{\delta\|v_{\eps,\eta,f}\|_\infty}
       \|v_{\eps,\eta,f}\|_\infty
       \int_{-\infty}^0\int_{\R^d}\lambda_0\cJ_{\eps,\eta}\left|v_{\eps,\eta,f}\right|\dd x\dd r
 \nonumber\\
 &\leq C\delta\|v_{\eps,\eta,f}\|_\infty\|v_{\eps,\eta,f}\|_{0,1}.
\end{align}
Let $C_{\eps,\eta}(f,g):=D^2\Lambda_{\eps,\eta}(0)(f,g)$.
Differentiating \eqref{eq:heavy-parent-cgf} using
\eqref{eq:cluster-first-second-derivatives} at the fractional scales gives
\begin{align}\label{eq:fractional-cumulant-derivatives}
 D\Lambda_{\eps,\eta}(0)f&=a_\eps \B_{\eps,\eta}\cA_{\eps,\eta}f,\\
 C_{\eps,\eta}(f,g)
 &=a_\eps \B_{\eps,\eta}\cA_{\eps,\eta}
                     [(\cA_{\eps,\eta}f)(\cA_{\eps,\eta}g)]+R_{\eps,f,g},\nonumber\\
 \left|R_{\eps,f,g}\right|&\leq C\delta\left\|\cA_{\eps,\eta}f\right\|_\infty
                                      \left\|\cA_{\eps,\eta}g\right\|_{0,1}\to0.\nonumber
\end{align}
Lemma~\ref{lem:fractional-mean-duality} identifies
$\B_\eta \cA_\eta w=\int_0^T\int_{\R^d}\rho^\eta w\dd x\dd t$.
Consequently the mean derivative tends to $\rho^\eta\dd t\dd x$,
and the second derivative tends to $C_\eta(f,g)$.
The constant microscopic past is a supersolution for the mean as
in \eqref{eq:meanbound}; the change of variables in
\eqref{eq:heavyfamily} gives
$\mu_\eps^\eta\leq\|\lambda_0\|_\infty\dd t\dd x$.

Cauchy's integral formula applied to the uniform analytic neighborhood
and \eqref{eq:fractional-parent-transform-remainder} gives, for fixed
compact spatial support $K$,
\begin{equation}\label{eq:fractional-cumulant-expansion}
 \Lambda_{\eps,\eta}(zf)=z\mu_\eps^\eta(f)
      +\frac12z^2C_{\eps,\eta}(f,f)+\mathcal O_K\left(|z|^3\|f\|_\infty^3\right).
\end{equation}
As in \eqref{eq:cauchy-coefficient-bound} and
\eqref{eq:centered-cgf-bound}, the same coefficient estimates give
\begin{equation}\label{eq:fractional-centered-cgf-bound}
 0\leq\Lambda_{\eps,\eta}(\theta f)-\theta\mu_\eps^\eta(f)
                    \leq C_K\theta^2\|f\|_\infty^2,
 \qquad |\theta|\|f\|_\infty\leq r.
\end{equation}
Here \eqref{eq:fractional-cumulant-expansion} holds for complex $z$
with $|z|\|f\|_\infty\leq r$, while
\eqref{eq:fractional-centered-cgf-bound} holds for real bounded
measurable $f$ supported in $[0,T]\times K$ and real $\theta$.
The source radius $r>0$ and the small-$\eps$ cutoff are independent
of $K$; $C_K$ depends on its volume and not on derivatives of $f$.
The Poisson parent formula \eqref{eq:heavy-parent-cgf} is linear in
the parent intensity. Multiplying that intensity by $n=n_{\eps,\eta}$
therefore gives the fractional counterpart of \eqref{eq:superposition}:
\begin{equation}\label{eq:fractional-superposition}
 \log\E\left[e^{n_{\eps,\eta}\NN_\eps^\eta(f)}\right]
 =n_{\eps,\eta}\Lambda_{\eps,\eta}(f).
\end{equation}
Taking $z=\mathrm{i}/\sqrt n$ in \eqref{eq:fractional-cumulant-expansion}
and using \eqref{eq:fractional-superposition} yields
\begin{equation}\label{eq:fractional-clt-characteristic}
 \log\E\left[e^{\mathrm{i}\sqrt n(\NN_\eps^\eta-\mu_\eps^\eta)(f)}\right]
 =-\frac12C_{\eps,\eta}(f,f)+\mathcal O_{K,f}\left(n^{-1/2}\right)
 \longrightarrow-\frac12C_\eta(f,f).
\end{equation}
Applying \eqref{eq:fractional-clt-characteristic} to
$f=\sum_{i=1}^k\theta_i f_i$ gives all joint Gaussian limits.
For the localized Fourier test functions in \eqref{eq:localized-fourier-norm},
exact Poisson superposition and differentiation of
\eqref{eq:fractional-centered-cgf-bound} at $\theta=0$ give
\[
\E\left[\left|Z_{\eps,\eta}\left(f_{j,k}^{(r)}\right)\right|^2\right]
=C_{\eps,\eta}\left( f_{j,k}^{(r)},f_{j,k}^{(r)}\right)\leq2C_j.
\]
Consequently, after increasing $C_j$,
\begin{equation}\label{eq:fractional-fourier-second-moment}
 \sup_{\substack{0<\eps\leq\eps_0,\ k\in\mathbb Z^{d+1}\\r\in\{1,2\}}}
       \E\left[\left|Z_{\eps,\eta}\left(f_{j,k}^{(r)}\right)\right|^2\right]\leq C_j.
\end{equation}
Using \eqref{eq:fractional-fourier-second-moment} in the summation
and Fourier tail estimates in the tightness part of
Lemma~\ref{lem:distribution-tightness} proves tightness of $Z_{\eps,\eta}$.
Let $\cTG_\eta f:=\sqrt{\rho^\eta}\cA_\eta f$. The mean bound in
Lemma~\ref{lem:fractional-mean-duality} and \eqref{eq:fractionalA} give
\[
 \left\|\cTG_\eta f\right\|_{L^2(E_T)}^2
 \leq T\|\lambda_0\|_\infty\|\cA_\eta f\|_\infty
                              \|\cA_\eta f\|_{0,1}<\infty.
\]
It\^o's isometry identifies the white-noise integral
\eqref{eq:fractionalGaussian} with the centered Gaussian field of
covariance $C_\eta$. Tightness and the joint limits from
\eqref{eq:fractional-clt-characteristic} identify every subsequential
limit with this field, completing the CLT proof.
\end{proof}

\begin{proof}[Proof of Theorem~\ref{thm:fractional-mdp}]
Let $n=n_{\eps,\eta}$. The cumulant expansion
\eqref{eq:fractional-cumulant-expansion}, proved above, and exact
Poisson superposition give, for $z=b_\eps/\sqrt n$,
\begin{equation}\label{eq:fractional-mdp-cgf}
b_\eps^{-2}\log\E\left[e^{b_\eps^2Y_{\eps,\eta}(f)}\right]
 =\frac12C_{\eps,\eta}(f,f)+\mathcal O_{K,f}\left(b_\eps/\sqrt n\right)
 \longrightarrow\frac12C_\eta(f,f).
\end{equation}
For $f=\sum_{i=1}^k\theta_i f_i$, this convergence is uniform on
compact subsets of $\R^k$. The exponential tilt
\eqref{eq:finite-projection-tilt} and its finite-dimensional upper
bound therefore give the projected MDP with covariance matrix
$\Gamma=(C_\eta(f_i,f_j))_{1\leq i,j\leq k}$, including singular $\Gamma$.

For the localized Fourier test functions in
\eqref{eq:localized-fourier-norm}, exact Poisson superposition and
\eqref{eq:fractional-centered-cgf-bound} give
\[
 \log\E\left[e^{\pm b_\eps^2Y_{\eps,\eta}(f_{j,k}^{(r)})}\right]
 \leq b_\eps^2C_j,
 \qquad k\in\mathbb Z^{d+1},\quad r\in\{1,2\}.
\]
The small-$\eps$ cutoff is common to $j,k,r$ because the source
radius in \eqref{eq:fractional-centered-cgf-bound} is independent
of the support. Chernoff's inequality now yields
\begin{equation}\label{eq:fractional-fourier-coordinate-tail}
\Pp\left(\left|Y_{\eps,\eta}\left(f_{j,k}^{(r)}\right)\right|>R_j(1+|k|)\right)
 \leq2e^{-b_\eps^2[R_j(1+|k|)-C_j]}.
\end{equation}
The sums over $j,k,r$ and the Fourier tail estimates in the
exponential-tightness part of Lemma~\ref{lem:distribution-tightness},
with \eqref{eq:fractional-fourier-coordinate-tail} replacing
\eqref{eq:fourier-coordinate-tail}, prove exponential tightness.
The cylinder lower bound and compact-cover upper bound in the
proof of Theorem~\ref{thm:mdp}, applied to
\eqref{eq:fractional-mdp-cgf}, then yield the full MDP with rate
\[
 \widetilde I_\eta(z)=\sup_{f\in\cD(E_T)}
              \left\{z(f)-\frac12\left\|\cTG_\eta f\right\|_2^2\right\},
 \qquad \cTG_\eta f=\sqrt{\rho^\eta}\cA_\eta f.
\]
If $\widetilde I_\eta(z)<\infty$, optimizing its defining supremum
over scalar multiples of $f$ gives
\[
 |z(f)|^2\leq2\widetilde I_\eta(z)\left\|\cTG_\eta f\right\|_2^2.
\]
The Riesz representation theorem on
$\mathsf H_\eta=\overline{\cTG_\eta\cD(E_T)}^{L^2(E_T)}$ gives a
unique $\zeta_0\in\mathsf H_\eta$ with
$z(f)=\ip{\cTG_\eta f}{\zeta_0}_2$. Completing the square gives
\[
 z(f)-\frac12\left\|\cTG_\eta f\right\|_2^2
 =\frac12\|\zeta_0\|_2^2
       -\frac12\left\|\cTG_\eta f-\zeta_0\right\|_2^2,
\]
so density of $\cTG_\eta\cD(E_T)$ in $\mathsf H_\eta$ gives
$\widetilde I_\eta(z)=\|\zeta_0\|_2^2/2=I_\eta(z)$ as in
\eqref{eq:fractionalrate}. Every other representing control
$\zeta$ satisfies $\zeta-\zeta_0\perp\mathsf H_\eta$ and has no
smaller norm. Conversely, any representing $\zeta\in L^2(E_T)$
bounds the supremum by $\|\zeta\|_2^2/2$, which proves the equality
also when the rate function is infinite. Finally, for $M\geq0$,
\[
 \left\{z:I_\eta(z)\leq M\right\}
 =\left\{f\mapsto\ip{\cTG_\eta f}{\zeta}_2:
                  \zeta\in L^2(E_T),\ \|\zeta\|_2^2\leq2M\right\}.
\]
The Hilbert ball is weakly compact and each evaluation of its image
is weakly continuous. Hence the sublevel sets are compact in the
weak distribution topology, proving goodness of $I_\eta$.
\end{proof}

\begin{remark}[The fractional noise equation]
In dual distribution notation, \eqref{eq:fractionalGaussian} solves
\[
 \left(c_\eta D_{0+}^{\eta}-\B+\kappa\right)Z_\eta=\sqrt{\rho^\eta}\,\dot W
\]
with zero past extension, where $D_{0+}^{\eta}=\partial_t I_{0+}^{1-\eta}$ is the causal Riemann--Liouville derivative, with
$I_{0+}^{1-\eta}r(t)=\Gamma(1-\eta)^{-1}\int_0^t(t-s)^{-\eta}r(s)\dd s$,
extended by duality to causal distributions. The mean equation
\eqref{eq:heavymean} uses a Caputo derivative with initial value
$\lambda_0$. The fluctuations in \eqref{eq:fractionalGaussian} are
white-noise integrals in physical time with the deterministic backward
resolvent $\cA_\eta$, and their covariance is \eqref{eq:heavycov}.
\end{remark}

\subsection{The full fractional space--time LDP}\label{subsec:fractional-full}

The macroscopic deviations of the space--time measure are governed by
a fractional control energy. The lower bound uses a regularization
that preserves the fractional source generated by the prescribed past;
the stable clock smooths time and the Brownian motion smooths space.
Throughout this subsection Assumption~\ref{ass:heavy} holds and
$0<L=\int_{\R^d}\lambda_0(x)\dd x<\infty$.

For a test function extended by zero after $T$, define
\begin{align}\label{eq:fractional-full-generator}
 D_+^\eta\phi(t,x)
 &:=\frac{\eta}{\Gamma(1-\eta)}\int_0^\infty
       \frac{\phi(t,x)-\phi(t+s,x)}{s^{1+\eta}}\dd s,\nonumber\\
 \mathscr L_\eta\phi&:=\B\phi-c_\eta D_+^\eta\phi-\kappa\phi,
 \qquad
 \beta_\eta(\dd t,\dd x):=c_Ht^{-\eta}\lambda_0(x)\dd t\dd x.
\end{align}
The causal adjoint of $-\mathscr L_\eta$ is
$\mathscr H_\eta=c_\eta D_{0+}^\eta-\B+\kappa$, with zero extension
to negative times. Let $\cR_\eta:=\cA_\eta^*$ denote the forward
potential, so that
\begin{equation}\label{eq:fractional-forward-potential}
 \ip{\cR_\eta\gamma}{f}=\ip\gamma{\cA_\eta f},\qquad
 \cR_\eta\beta_\eta=\rho^\eta\dd t\dd x.
\end{equation}
It maps finite signed measures to $L^1(\Omega_T)$ densities.
All $L^1$ norms in this subsection are over $\Omega_T$ unless a
different domain is specified. Define
\begin{align}\label{eq:fractional-variational-action}
 J_\eta^{\rm loc}(\nu)
 &:=\sup_{\phi\in C_c^\infty(\Omega_T)}
   \left\{-\nu\left(\mathscr L_\eta\phi+\phi^2/2\right)-\beta_\eta(\phi)\right\},
 \nonumber\\
 \mathcal J_{\eta,\lambda_0}(\nu)
 &:=\begin{cases}J_\eta^{\rm loc}(\nu),&\nu(\Omega_T)<\infty,\\
                 +\infty,&\nu(\Omega_T)=\infty.
   \end{cases}
\end{align}

The next theorem concerns the space--time event measure in
Theorems~\ref{thm:fractional} and \ref{thm:fractional-mdp}, at its
macroscopic deviation scale.

\begin{theorem}[Fractional large deviations]
\label{thm:fractional-spatial-ldp}
Under Assumption~\ref{ass:heavy} and $0<\int_{\mathbb{R}^{d}}\lambda_{0}(x)\dd x<\infty$,
for every $0<\eta<1$, the measures $\NN_\eps^\eta$ satisfy a full
LDP on $\mathsf E$ with the vague topology, speed $n_{\eps,\eta}$,
and good rate function $\mathcal J_{\eta,\lambda_0}$.
Its finite domain consists of measures $q\dd t\dd x$, where
$q\geq0$ belongs to $L^1(\Omega_T)$ and there is a
$j\in L^1(\Omega_T)$ such that
\begin{equation}\label{eq:fractional-controlled-potential}
 q=\rho^\eta+\cR_\eta j,\qquad
 \int_{\Omega_T}\frac{(j(t,x))^2}{q(t,x)}\dd t\dd x<\infty.
\end{equation}
Here $j^2/q=0$ when $j=q=0$ and is $+\infty$ when $q=0$, $j\ne0$.
The residual $j=\mathscr H_\eta q-\beta_\eta$ is unique, and
\begin{equation}\label{eq:fractional-full-energy}
 \mathcal J_{\eta,\lambda_0}(q\dd t\dd x)
       =\frac12\int_{\Omega_T}\frac{(j(t,x))^2}{q(t,x)}\dd t\dd x.
\end{equation}
\end{theorem}

The source $\beta_\eta$ in \eqref{eq:fractional-controlled-potential}
is part of the rate function domain. For sufficiently regular $q$ with
$q(0)=\lambda_0$, its residual can also be written
$c_\eta\,{}^{\rm C}D_t^\eta q-\B q+\kappa q$.
The causal formulation \eqref{eq:fractional-controlled-potential}
specifies the domain without assuming a pointwise initial trace.

We first identify finite action and its compact sublevels; this also
provides the integrability required by the regularization below.
\begin{lemma}[Fractional finite-action representation and coercivity]
\label{lem:fractional-action-domain}
For finite $\nu$, equations
\eqref{eq:fractional-controlled-potential}--\eqref{eq:fractional-full-energy}
are equivalent to $J_\eta^{\rm loc}(\nu)<\infty$.
There are constants $c,C>0$, depending only on the fixed model and
$T$, such that
\begin{equation}\label{eq:fractional-action-mass-bound}
 \mathcal J_{\eta,\lambda_0}(\nu)\geq c\nu(\Omega_T)-C.
\end{equation}
The extended rate function in \eqref{eq:fractional-variational-action} is good.
\end{lemma}

\begin{proof}
Let $\ell_\nu(\phi):=-\nu(\mathscr L_\eta\phi)-\beta_\eta(\phi)$.
Optimizing over real multiples of $\phi$ gives
\[
 |\ell_\nu(\phi)|^2\leq2J_\eta^{\rm loc}(\nu)\nu(\phi^2).
\]
Density of compact smooth functions in $L^2(\nu)$ and the Riesz
representation theorem give a unique $v\in L^2(\nu)$ with
$\ell_\nu(\phi)=\nu(v\phi)$ and
$J_\eta^{\rm loc}(\nu)=\nu(v^2)/2$. Since $\nu$ is finite,
$v\nu$ is a finite signed measure. Let $\gamma=\beta_\eta+v\nu$.
The potential in \eqref{eq:fractional-forward-potential} satisfies
\begin{equation}\label{eq:fractional-forward-potential-L1-bound}
 \left\|\cR_\eta\gamma\right\|_1\leq M_T|\gamma|(\Omega_T),\qquad
 M_T=\frac{T^\eta}{c_\eta\Gamma(1+\eta)}.
\end{equation}
Indeed its transition kernel has a density for every positive
genealogical time. The stable-clock law in \eqref{eq:fractionalA}
gives the zero-killing potential transform
\[
 \int_0^\infty e^{-zt}
       \left[\int_0^\infty\Pp\left(S_r^\eta\in\dd t\right)\dd r\right]
 =\int_0^\infty e^{-rc_\eta z^\eta}\dd r
 =\frac1{c_\eta z^\eta},\qquad z>0.
\]
Its density is $t^{\eta-1}/(c_\eta\Gamma(\eta))$.
Dropping killing and integrating this density gives
\[
 \int_0^\infty e^{-\kappa r}\Pp\left(S_r^\eta\leq T\right)\dd r
 \leq\int_0^T\frac{t^{\eta-1}}{c_\eta\Gamma(\eta)}\dd t
 =M_T,
\]
which proves the claimed $L^1$ potential bound.
For $f\in C_c^\infty(\Omega_T)$, the function $\cA_\eta f$ is
smooth, vanishes after the temporal support of $f$, and has bounded
spatial derivatives. Apply the weak equation to
$\chi_R\cA_\eta f$, where $\chi_R$ is a spatial cutoff. The
identity $-\mathscr L_\eta\cA_\eta f=f$ and the product rule give
\[
 \left|\nu(f)-\gamma(\cA_\eta f)\right|
 \leq C_f\left(R^{-1}+R^{-2}\right)\nu(\Omega_T)+o_R(1)\longrightarrow0.
\]
The term $o_R(1)$ comes from dominated convergence against the finite
measures $\gamma$ and $\nu$. Consequently $\nu=\cR_\eta\gamma$,
so $\nu=q\dd t\dd x$ and $v\nu=j\dd t\dd x$ with $j=vq$.
This proves \eqref{eq:fractional-controlled-potential} and
\eqref{eq:fractional-full-energy}. Conversely that potential equation
gives $\ell_\nu(\phi)=\int_{\Omega_T}j\phi\dd t\dd x$, and square completion
and density in $L^2(q\dd t\dd x)$ give the same energy.

To prove \eqref{eq:fractional-action-mass-bound}, choose smooth
nonincreasing $g_k:[0,T]\to[0,1]$, equal to one on
$[0,T-2/k]$ and zero on $[T-1/k,T]$, for $k>2/T$.
Let $u_k:=\cA_\eta g_k$, with no spatial dependence. Positivity of
the potential gives $0\leq u_k\leq M_Tg_k$ and
$-\mathscr L_\eta u_k=g_k$. Spatial cutoffs and finite mass extend
the action inequality to $a u_k$. For $0<a\leq M_T^{-2}$,
\[
 J_\eta^{\rm loc}(\nu)
 \geq\nu\left(a g_k-a^2u_k^2/2\right)-a\beta_\eta(u_k)
 \geq\frac a2\nu(g_k)-aM_T\beta_\eta(\Omega_T).
\]
Letting $k\to\infty$ proves the mass bound. The function
$J_\eta^{\rm loc}$ is lower semicontinuous, since
$\mathscr L_\eta\phi+\phi^2/2$ is continuous with compact support
in $\Omega_T$. Every sublevel of the extended rate function therefore lies
in a fixed vaguely compact set $\{\nu:\nu(\Omega_T)\leq R\}$
and is closed there. This proves the goodness of the rate function.
\end{proof}

The next lemma is the approximation step for the full lower bound. It
approximates every finite-action measure while keeping the source
generated by the prescribed past exactly equal to $\beta_\eta$.

\begin{lemma}[Source-preserving fractional approximation]
\label{lem:fractional-recovery}
If $\mathcal J_{\eta,\lambda_0}(q\dd t\dd x)<\infty$, there are
real $v_k\in C_c^\infty(E_T)$ such that the nonnegative solutions
\begin{equation}\label{eq:fractional-bounded-control}
 q^{v_k}=\rho^\eta+\cR_\eta(v_k q^{v_k})
\end{equation}
satisfy
\begin{equation}\label{eq:fractional-recovery-conclusion}
 \|q^{v_k}-q\|_1\longrightarrow0,\qquad
 \frac12\int_{\Omega_T}v_k^2q^{v_k}\dd t\dd x
       \longrightarrow\mathcal J_{\eta,\lambda_0}(q\dd t\dd x).
\end{equation}
\end{lemma}
\begin{proof}
Let $p_h^\eta$ be the density of $S_h^\eta$, and let $g_h^\Sigma$
be the centered Gaussian density of covariance $h\Sigma$.
The characteristic function of $S_h^\eta$ has modulus
$\exp\left\{-hc_\eta\cos(\pi\eta/2)|\xi|^\eta\right\}$; its integrability,
also after multiplication by any polynomial, gives a bounded smooth
density. For $h>0$ and $f\in L^1(\Omega_T)$, extended by zero
outside $[0,T)$ in the time variable, define the forward killed semigroup by
\begin{equation}\label{eq:fractional-smoothing-semigroup}
 (\mathsf Q_h f)(t,x):=e^{-\kappa h}\int_0^t p_h^\eta(t-s)
             \int_{\R^d}g_h^\Sigma(x-y)f(s,y)\dd y\dd s.
\end{equation}
Set $\mathsf Q_0=\Id$. It is a positive contraction on $L^1(\Omega_T)$, strongly continuous
as $h\downarrow0$. Independent increments of $S^\eta$ and
$B^\Sigma$ give the semigroup identity, and
\eqref{eq:fractionalA} and \eqref{eq:fractional-forward-potential}
give its potential representation:
\begin{equation}\label{eq:fractional-semigroup-potential-identity}
 \mathsf Q_h\mathsf Q_r=\mathsf Q_{h+r},\qquad
 \cR_\eta f=\int_0^\infty\mathsf Q_r f\dd r,\qquad
 \mathsf Q_h\cR_\eta f
 =\cR_\eta\mathsf Q_h f
 =\int_h^\infty\mathsf Q_r f\dd r .
\end{equation}
These equalities hold in $L^1(\Omega_T)$ for $f\in L^1(\Omega_T)$;
positivity and \eqref{eq:fractional-forward-potential-L1-bound} give
\[
 \int_0^\infty\|\mathsf Q_r f\|_1\dd r
 \leq\|\cR_\eta|f|\|_1
 \leq\frac{T^\eta}{c_\eta\Gamma(1+\eta)}\|f\|_1<\infty,
\]
which justifies the integrals and Fubini's theorem in
\eqref{eq:fractional-semigroup-potential-identity}.
In this proof write $\rho=\rho^\eta$ and let $j$ be the residual
in \eqref{eq:fractional-controlled-potential}. Since
$\rho=\cR_\eta\beta_\eta$ by \eqref{eq:fractional-forward-potential},
subtracting the last integral in
\eqref{eq:fractional-semigroup-potential-identity} from the second
one, with $f(t,x)=c_Ht^{-\eta}\lambda_0(x)$, gives the
nonnegative correction
\begin{equation}\label{eq:fractional-positive-source-correction}
 d_h:=\rho-\mathsf Q_h\rho
       =\int_0^h\mathsf Q_r\beta_\eta\dd r\geq0.
\end{equation}
Define
\begin{equation}\label{eq:fractional-smoothed-pair}
 q_h:=\mathsf Q_hq+d_h,\qquad j_h:=\mathsf Q_hj.
\end{equation}
Then $q_h=\rho+\cR_\eta j_h$, $q_h\geq0$, and
$(q_h,j_h)\to(q,j)$ in $L^1\times L^1$. The Cauchy--Schwarz
inequality for the positive kernel in
\eqref{eq:fractional-smoothing-semigroup} gives
\begin{equation}\label{eq:fractional-smoothing-energy}
 \int_{\Omega_T}\frac{j_h^2}{q_h}\dd t\dd x
 \leq\int_{\Omega_T}\frac{(\mathsf Q_hj)^2}{\mathsf Q_hq}\dd t\dd x
 \leq\int_{\Omega_T}\mathsf Q_h\left(j^2/q\right)\dd t\dd x
 \leq\int_{\Omega_T}\frac{j^2}{q}\dd t\dd x.
\end{equation}

We next truncate space. Choose $0\leq\chi_R\leq1$, equal to one
on $B_R$ and zero outside $B_{2R}$, with
\begin{equation}\label{eq:fractional-flat-cutoff}
 \frac{|\nabla\chi_R|^2}{\chi_R(1-\chi_R)}\leq C R^{-2},
 \qquad
 \frac{|D^2\chi_R|^2}{\chi_R(1-\chi_R)}\leq C R^{-4}.
\end{equation}
The ratios are zero where the numerator and denominator vanish.
Such cutoffs result by rescaling a smooth radial cutoff flat at both
ends of its transition interval, for example one constructed from
$a(s)=e^{-1/s}\one_{\{s>0\}}$.
Let $a_h:=\mathsf Q_hq$, $b_h:=\mathsf Q_h\rho$, and define
\begin{align}\label{eq:fractional-space-truncated-pair}
 q_{h,R}&:=\rho+\chi_R(a_h-b_h)
          =\chi_Rq_h+(1-\chi_R)\rho,\nonumber\\
 j_{h,R}&:=\chi_Rj_h-e_{h,R},\nonumber\\
 e_{h,R}&:=(\B\chi_R)(a_h-b_h)
             +(\Sigma\nabla\chi_R)\cdot\nabla(a_h-b_h).
\end{align}
The product rule for $\mathscr H_\eta$ gives
$q_{h,R}=\rho+\cR_\eta j_{h,R}$, with the same source
$\beta_\eta$. Moreover
$q_{h,R}=\chi_R a_h+(1-\chi_R)b_h+d_h\geq0$.

The Gaussian density in \eqref{eq:fractional-smoothing-semigroup}
and its spatial derivative are explicitly
\begin{equation}\label{eq:fractional-Gaussian-kernel-derivative}
 g_h^\Sigma(z)
 =\frac{\exp\{-z^\top\Sigma^{-1}z/(2h)\}}
        {(2\pi h)^{d/2}(\det\Sigma)^{1/2}},
 \qquad
 \Sigma^{1/2}\nabla g_h^\Sigma(z)
 =-\frac{\Sigma^{-1/2}z}{h}g_h^\Sigma(z).
\end{equation}
For $f\geq0$ in $L^1(\Omega_T)$, differentiating
\eqref{eq:fractional-smoothing-semigroup} under the integral is
justified by boundedness of $p_h^\eta$ and $\nabla g_h^\Sigma$.
The Cauchy--Schwarz inequality with respect to the positive measure
$e^{-\kappa h}p_h^\eta(t-s)g_h^\Sigma(x-y)f(s,y)\dd y\dd s$
gives, for $0<t<T$,
\begin{equation}\label{eq:fractional-semigroup-fisher-pointwise}
 \frac{|\Sigma^{1/2}\nabla\mathsf Q_hf(t,x)|^2}
      {\mathsf Q_hf(t,x)}
 \leq e^{-\kappa h}\int_0^t p_h^\eta(t-s)
       \int_{\R^d}\frac{|\Sigma^{-1/2}(x-y)|^2}{h^2}
            g_h^\Sigma(x-y)f(s,y)\dd y\dd s .
\end{equation}
The ratio is zero when $\mathsf Q_hf(t,x)=0$.
Since $g_h^\Sigma$ has covariance $h\Sigma$,
\[
 \int_{\R^d}\frac{|\Sigma^{-1/2}z|^2}{h^2}
                       g_h^\Sigma(z)\dd z=\frac dh.
\]
After integration in $x$, the temporal kernel left in
\eqref{eq:fractional-semigroup-fisher-pointwise} is
\begin{equation}\label{eq:fractional-temporal-smoothing-kernel}
 k_h^\eta(t,s):=e^{-\kappa h}p_h^\eta(t-s)\one_{\{0<s<t<T\}},
 \qquad
 \int_s^T k_h^\eta(t,s)\dd t
 =e^{-\kappa h}\Pp(S_h^\eta<T-s)\leq1
 \quad(0<s<T).
\end{equation}
Tonelli's theorem now yields
\begin{align}\label{eq:fractional-semigroup-fisher}
 \int_{\Omega_T}
       \frac{|\Sigma^{1/2}\nabla\mathsf Q_hf|^2}{\mathsf Q_hf}
          \dd t\dd x
 \leq\frac dh\int_0^T\int_{\R^d} f(s,y)
                  \left[\int_s^T k_h^\eta(t,s)\dd t\right]\dd y\dd s
\leq\frac dh\|f\|_1 .
\end{align}
To derive the cutoff energy estimate, split the residual in
\eqref{eq:fractional-space-truncated-pair} into
\begin{equation}\label{eq:fractional-cutoff-residual-parts}
 e_{h,R}^{(0)}:=(\B\chi_R)(a_h-b_h),\qquad
 e_{h,R}^{(1)}:=(\Sigma\nabla\chi_R)\cdot\nabla(a_h-b_h),\qquad
 e_{h,R}=e_{h,R}^{(0)}+e_{h,R}^{(1)}.
\end{equation}
On the transition region $0<\chi_R<1$, the denominator satisfies
$q_{h,R}\geq\chi_Ra_h+(1-\chi_R)b_h$ by
\eqref{eq:fractional-space-truncated-pair}. Therefore,
\begin{align}\label{eq:fractional-cutoff-pointwise-bounds}
 \frac{\chi_R(1-\chi_R)|a_h-b_h|^2}{q_{h,R}}
 &\leq2\big((1-\chi_R)a_h+\chi_Rb_h\big)
 \leq2(a_h+b_h),\nonumber\\
 \frac{\chi_R(1-\chi_R)
             |\Sigma^{1/2}\nabla(a_h-b_h)|^2}{q_{h,R}}
 &\leq2(1-\chi_R)\frac{|\Sigma^{1/2}\nabla a_h|^2}{a_h}
       +2\chi_R\frac{|\Sigma^{1/2}\nabla b_h|^2}{b_h}.
\end{align}
These inequalities follow by splitting each squared difference and
using $q_{h,R}\geq\chi_Ra_h$ and $q_{h,R}\geq(1-\chi_R)b_h$
for its two terms. The ratios are zero at simultaneous zeros.
Using \eqref{eq:fractional-flat-cutoff} in
\eqref{eq:fractional-cutoff-pointwise-bounds} now gives the two
explicit integral bounds
\begin{align}\label{eq:fractional-cutoff-component-energy}
 \int_{\Omega_T}\frac{|e_{h,R}^{(0)}|^2}{q_{h,R}}\dd t\dd x
 &\leq CR^{-4}\int_{\Omega_T}(a_h+b_h)\dd t\dd x,\nonumber\\
 \int_{\Omega_T}\frac{|e_{h,R}^{(1)}|^2}{q_{h,R}}\dd t\dd x
 &\leq CR^{-2}\int_{\Omega_T}
       \left(\frac{|\Sigma^{1/2}\nabla a_h|^2}{a_h}
             +\frac{|\Sigma^{1/2}\nabla b_h|^2}{b_h}\right)\dd t\dd x.
\end{align}
For the second inequality, write
$e_{h,R}^{(1)}=(\Sigma^{1/2}\nabla\chi_R)\cdot
\Sigma^{1/2}\nabla(a_h-b_h)$ and apply the Cauchy--Schwarz inequality.
More precisely, let $C_\chi$ be a common constant in the two
inequalities of \eqref{eq:fractional-flat-cutoff}, and use the
Frobenius norm for $D^2\chi_R$. The inequalities
\[
 |\B\chi_R|^2\leq\frac14\|\Sigma\|_{\rm F}^2\left|D^2\chi_R\right|^2,
 \qquad
 \left|\Sigma^{1/2}\nabla\chi_R\right|^2
       \leq\|\Sigma\|_{\rm op}|\nabla\chi_R|^2
\]
show that the constants $C$ in the first and second lines of
\eqref{eq:fractional-cutoff-component-energy} can be taken as
$C_0=C_\chi\|\Sigma\|_{\rm F}^2/2$ and
$C_1=2C_\chi\|\Sigma\|_{\rm op}$, respectively. They depend only on
the fixed cutoff profile and $\Sigma$, and are independent of $h,R$.
Since $\mathsf Q_h$ is an $L^1$ contraction,
$\|a_h\|_1+\|b_h\|_1\leq\|q\|_1+\|\rho\|_1$.
Combining \eqref{eq:fractional-cutoff-component-energy} with
\eqref{eq:fractional-semigroup-fisher} and
$|e_{h,R}|^2\leq2\left|e_{h,R}^{(0)}\right|^2+2\left|e_{h,R}^{(1)}\right|^2$ proves
\begin{equation}\label{eq:fractional-cutoff-energy}
 \int_{\Omega_T}\frac{e_{h,R}^2}{q_{h,R}}\dd t\dd x
 \leq C\left(R^{-4}+h^{-1}R^{-2}\right)
                 (\|q\|_1+\|\rho\|_1).
\end{equation}
Since $q_{h,R}\geq\chi_Rq_h$, for every $u>0$,
\begin{equation}\label{eq:fractional-truncated-cost}
 \int_{\Omega_T}\frac{j_{h,R}^2}{q_{h,R}}\dd t\dd x
 \leq(1+u)\int_{\Omega_T}\frac{j_h^2}{q_h}\dd t\dd x
       +\left(1+u^{-1}\right)\int_{\Omega_T}\frac{e_{h,R}^2}{q_{h,R}}\dd t\dd x.
\end{equation}
Also
$\|q_{h,R}-q_h\|_1\to0$ as $R\to\infty$. First taking
$R\to\infty$ and then $u\downarrow0$ proves approximation of $q_h$
with no increase of cost.

We verify that $v_{h,R}=j_{h,R}/q_{h,R}$ is bounded, including
near physical time zero. Finite nonzero $L$ and positivity of the
Gaussian kernel imply, for each compact spatial set $K$,
\[
 \inf_{x\in K,\,\, h/4\leq r\leq h/2}
            \left(g_r^\Sigma*\lambda_0\right)(x)>0.
\]
Because $c_Hs^{-\eta}\geq c_HT^{-\eta}$ for $0<s\leq T$,
\eqref{eq:fractional-positive-source-correction} gives
\begin{equation}\label{eq:fractional-correction-lower}
 \inf_{x\in K}d_h(t,x)
       \geq c_{K,h,T}\Pp\left(S_{h/2}^\eta\leq t\right)>0,
       \qquad 0<t\leq T.
\end{equation}
The last probability is positive: splitting the subordinator into
jumps above a small threshold and smaller jumps gives a positive
probability of no large jump and total small-jump mass below $t$.
The convolution identity $p_h^\eta=p_{h/2}^\eta*p_{h/2}^\eta$
gives
\[
 \sup_{0<s\leq t}p_h^\eta(s)
 \leq\left\|p_{h/2}^\eta\right\|_\infty\Pp\left(S_{h/2}^\eta\leq t\right).
\]
Consequently, for $f\in L^1(\Omega_T)$ and $|\alpha|\leq1$,
\begin{equation}\label{eq:fractional-smoothed-function-bound}
 \sup_{x\in\R^d}\left|D_x^\alpha\mathsf Q_h f(t,x)\right|
 \leq C_{h,\alpha}\|f\|_1\Pp\left(S_{h/2}^\eta\leq t\right).
\end{equation}
Apply \eqref{eq:fractional-smoothed-function-bound} to $j$, $q$,
and $\rho$, and use $q_{h,R}\geq d_h$ and
\eqref{eq:fractional-correction-lower} on $K=\overline B_{2R}$.
This proves $\|v_{h,R}\|_\infty<\infty$. Its spatial support is
compact, and it is continuous for $t>0$.

For any bounded $v$, the equation $q^v=\rho+\cR_\eta(vq^v)$ has
a unique $L^1$ solution. In the norm
$\|f\|_{1,-z}=\int_{\Omega_T}e^{-zt}|f(t,x)|\dd t\dd x$,
\begin{equation}\label{eq:fractional-controlled-contraction}
 \|\cR_\eta(vf)\|_{1,-z}
       \leq\frac{\|v\|_\infty}{\kappa+c_\eta z^\eta}
                       \|f\|_{1,-z}.
\end{equation}
Choose $z$ so the coefficient is less than one. 
For positivity, set $M=\|v\|_{\infty}$ and let $\mathcal{R}_{\eta,M}$ 
be the forward potential with killing $\kappa+M$. 
The resolvent identity rewrites the equation as
\[
q^{v}=\mathcal{R}_{\eta,M}\beta_{\eta}+\mathcal{R}_{\eta,M}((v+M)q^{v}).
\]
Choose $z$ still larger, if necessary, so that $2M/(\kappa+M+c_{\eta}z^{\eta})<1$. 
Its convergent Neumann series is nonnegative, since both $\beta_{\eta}$ and $v+M$ are nonnegative.

Approximate $v_{h,R}$ pointwise on $E_T$ by
$v_{h,R,k}\in C_c^\infty(E_T)$ with a common uniform bound.
This can be done by temporal cutoffs and mollification. Subtraction
in \eqref{eq:fractional-controlled-contraction} and the dominated
convergence theorem give
\[
 \left\|q^{v_{h,R,k}}-q_{h,R}\right\|_1\longrightarrow0,
 \qquad
 \int_{\Omega_T}v_{h,R,k}^2q^{v_{h,R,k}}\dd t\dd x
       \longrightarrow\int_{\Omega_T}\frac{j_{h,R}^2}{q_{h,R}}\dd t\dd x.
\]
Combining this with \eqref{eq:fractional-smoothing-energy} and
\eqref{eq:fractional-truncated-cost}, and choosing a diagonal sequence
as $h\downarrow0$ and $R,k\to\infty$, gives the required convergence
and the upper bound on the limiting costs. Lower semicontinuity from
Lemma~\ref{lem:fractional-action-domain} gives equality of the costs
in \eqref{eq:fractional-recovery-conclusion}.
\end{proof}

The microscopic exponential identity yields the upper bound, while a
bounded change of reproduction intensity realizes the controls in
Lemma~\ref{lem:fractional-recovery}.
\begin{lemma}[Fractional exponential bounds and controlled limits]
\label{lem:fractional-micro-bounds}
For every vague closed $F\subset\mathsf E$,
\begin{equation}\label{eq:fractional-full-upper}
 \limsup_{\eps\downarrow0}\frac1{n_{\eps,\eta}}
       \log\Pp(\NN_\eps^\eta\in F)
       \leq-\inf_{\nu\in F}\mathcal J_{\eta,\lambda_0}(\nu).
\end{equation}
For every $v\in C_c^\infty(E_T)$ and every vague open neighborhood
$O$ of $q^v\dd t\dd x$,
\begin{equation}\label{eq:fractional-controlled-lower}
 \liminf_{\eps\downarrow0}\frac1{n_{\eps,\eta}}
       \log\Pp(\NN_\eps^\eta\in O)
       \geq-\frac12\int_{\Omega_T}v^2q^v\dd t\dd x.
\end{equation}
\end{lemma}
\begin{proof}
Write $\delta=\delta_{\eps,\eta}$, $\tau=\tau_\eps$,
$n=n_{\eps,\eta}$ and
$\mathscr L_{\eps,\eta}=(a_\eps\cK_{\eps,\eta}-\Id)/\delta$.
We first bound the mass of the space--time measure exponentially.
Choose $c>0$ in the common source neighborhood constructed in the
proof of Theorem~\ref{thm:fractional}, and let
$f_k(t,x)=c\one_{[0,T]}(t)\chi_k(x)$, where the compact spatial
cutoffs $\chi_k$ increase to one. The contraction estimate
\eqref{eq:fractional-cluster-contraction} with the common resolvent bound
\eqref{eq:fractional-common-resolvent-bound} gives
$0\leq v_{\eps,\eta,f_k}\leq V$ with $V$ independent of $k$
and sufficiently small $\eps$. By \eqref{eq:heavy-parent-cgf},
\eqref{eq:heavyB}, and \eqref{eq:fractional-preparation-tail-bound},
\begin{align*}
 \Lambda_{\eps,\eta}(f_k)
 \leq a_\eps e^{\delta V}
               \B_{\eps,\eta}v_{\eps,\eta,f_k}
               \leq e^{\delta V}LV\int_0^T b_\eps(t)\dd t
 \leq C_TLV e^{\delta V}\leq C.
\end{align*}
The exact superposition identity \eqref{eq:fractional-superposition}
and the monotone convergence theorem, applied as $f_k\uparrow c\one_{\Omega_T}$, give
\begin{equation}\label{eq:fractional-mass-exponential-moment}
 \log\E\left[e^{nc\NN_\eps^\eta(\Omega_T)}\right]\leq nC.
\end{equation}
Consequently, Chernoff's inequality yields
\begin{equation}\label{eq:fractional-mass-exp-tail}
 \Pp\left(\NN_\eps^\eta(\Omega_T)>R\right)\leq e^{-n(cR-C)}.
\end{equation}

For a compact smooth $\phi$, integration against the delay tail gives
\[
 \frac{\E[\phi(t+\tau S,x)-\phi(t,x)]}{\delta}
       =\int_0^T b_\eps(s)\partial_t\phi(t+s,x)\dd s
       \longrightarrow-c_\eta D_+^\eta\phi(t,x),
\]
uniformly for $0\leq t<T$ and $x\in\R^d$. Indeed
\begin{align*}
 \sup_{\substack{0\leq t<T\\x\in\R^d}}
 \left|\frac{\E[\phi(t+\tau S,x)-\phi(t,x)]}{\delta}
                         +c_\eta D_+^\eta\phi(t,x)\right|
\leq\|\partial_t\phi\|_\infty
           \int_0^T\left|b_\eps(s)-c_Hs^{-\eta}\right|\dd s\longrightarrow0,
\end{align*}
by \eqref{eq:fractional-preparation-tail-bound}. Factorization and
the centered spatial Taylor expansion give, on the full domain,
\begin{equation}\label{eq:fractional-global-generator}
 \left\|\mathscr L_{\eps,\eta}\phi-\mathscr L_\eta\phi\right\|_\infty
       \longrightarrow0.
\end{equation}
The spatial Taylor remainder is bounded by
$\frac12\E_p\left[|\Xi|^2\omega\left(\sqrt\delta|\Xi|\right)\right]$ for a bounded
modulus of continuity of $D_x^2\phi$; the temporal shift of $\B\phi$
converges uniformly to $\B\phi$.

Use \eqref{eq:expmart} with these time and space scales, and let
$\cI_{\eps,\eta}$ be \eqref{eq:pastcompensator-definition} with
$\delta,n,\tau,\eps^\eta$ and the heavy-tail past. The algebra of
\eqref{eq:exactmartfunctional} gives
\begin{equation}\label{eq:fractional-exact-martingale}
 n^{-1}\log\mathscr E_{\eps,\eta}(\phi)
 =\NN_\eps^\eta\left(-\mathscr L_{\eps,\eta}\phi
                   -a_\eps\cK_{\eps,\eta}c_\delta(\phi)\right)
       -\delta^{-1}\cI_{\eps,\eta}(\phi)
       -\cI_{\eps,\eta}(c_\delta(\phi)),
\end{equation}
where $c_\delta(z)=(e^{\delta z}-1-\delta z)/\delta^2$.
Its conditional expectation given the past is at most one by the
stopping and Fatou argument in Lemma~\ref{lem:martingale-past}.
Let $X_\eps(\phi)=\delta^{-1}\cI_{\eps,\eta}(\phi)$.
The parent intensity in microscopic time and macroscopic space is
$n\tau\delta^{-2}\lambda_0(x)\dd r\dd x$ on
$(-\infty,0)\times\R^d$. Applying the centered Poisson exponential
identity used in \eqref{eq:past-centered-exponential-bound} with this
intensity gives, 
\[
\E[X_\eps(\phi)]
=a_\eps\B_{\eps,\eta}\phi\longrightarrow\beta_\eta(\phi),
\]
and for every real $\theta$,
\begin{align}\label{eq:fractional-past-superexp}
 &\frac1n\log\E\left[e^{n\theta(X_\eps(\phi)-\E[X_\eps(\phi)])}\right]
 \nonumber\\
 &=\frac\tau{\delta^2}\int_{-\infty}^0\int_{\R^d}\lambda_0(x)
   \left[e^{\theta a_\eps\delta\cJ_{\eps,\eta}\phi(r,x)}-1
                    -\theta a_\eps\delta\cJ_{\eps,\eta}\phi(r,x)\right]
                    \dd x\dd r\nonumber\\
 &\leq\frac{\theta^2a_\eps^2\delta}{2}
       e^{|\theta|\delta\|\phi\|_\infty}
       \|\phi\|_\infty\B_{\eps,\eta}|\phi|\nonumber\\
 &\leq C_\phi\theta^2\delta e^{|\theta|\delta\|\phi\|_\infty},
\end{align}
where we used $|e^y-1-y|\leq y^2e^{|y|}/2$ and
$|\cJ_{\eps,\eta}\phi|^2\leq\|\phi\|_\infty\cJ_{\eps,\eta}|\phi|$.
The crossing identity following \eqref{eq:heavy-parent-cgf} is
\[
 \int_{-\infty}^0\int_{\R^d}\lambda_0(x)
       \cJ_{\eps,\eta}|\phi|(r,x)\dd x\dd r
 =\frac\delta\tau\B_{\eps,\eta}|\phi|,
\]
and \eqref{eq:heavyB} bounds $\B_{\eps,\eta}|\phi|$ uniformly for
small $\eps$. These bounds also justify exhaustion of the parent
region in the Poisson formula.
Thus, $X_\eps(\phi)\simeq\beta_\eta(\phi)$ and
$\cI_{\eps,\eta}(c_\delta(\phi))\simeq0$ at speed $n$.
Taylor's formula and \eqref{eq:fractional-global-generator} show
that the event-measure integrand in
\eqref{eq:fractional-exact-martingale} converges uniformly to
$-\mathscr L_\eta\phi-\phi^2/2$.
On $\{\NN_\eps^\eta(\Omega_T)\leq R\}$ its pairing error is at
most $R$ times that uniform error. Exponential Markov's inequality
therefore gives the local upper bound for each affine functional
in \eqref{eq:fractional-variational-action}. A finite cover gives
the upper bound on each vague compact subset of
$\{\nu:\nu(\Omega_T)\leq R\}$. The mass bound
\eqref{eq:fractional-mass-exp-tail}, followed by $R\to\infty$,
proves \eqref{eq:fractional-full-upper}.

For the lower bound, keep the past law unchanged and multiply the
positive-time intensity by
$e^{\delta v(\tau t,\eps^\eta x)}$. Write $\Pp^v$ for this law
and $\NN_\eps^{\eta,v}$ for its event measure. With $M=\|v\|_\infty$, domination by the constant multiplier
$b=a_\eps e^{\delta M}$ gives the same temporal renewal equation
as \eqref{eq:tilted-mean-renewal}, with time scale $\tau$. Its positive
renewal iteration gives
\begin{align}\label{eq:fractional-controlled-scalar-mean}
 R_\eps(t)
 =\sum_{j=0}^\infty b^{j+1}
       \Pp\left(\tau T_j\leq t<\tau T_{j+1}\right)
       =1+(b-1)\sum_{j=0}^\infty b^j\Pp\left(\tau T_j\leq t\right).
\end{align}
The second equality follows by writing each probability in the first
line as $\Pp(\tau T_j\leq t)-\Pp(\tau T_{j+1}\leq t)$ and shifting
the summation index. The bound below justifies the convergent sums.
Choose $z>0$ with $\kappa+c_\eta z^\eta>M$. By
\eqref{eq:heavyassumptions} and the delay transform limit in
\eqref{eq:fractional-renewal-transform},
\begin{align}\label{eq:fractional-controlled-renewal-bound}
 1-b\widehat H(z\tau)
 &=\delta\left(\kappa+c_\eta z^\eta-M+o(1)\right),\nonumber\\
 \sum_{j=0}^\infty b^j\Pp\left(\tau T_j\leq t\right)
 &\leq e^{zT}\sum_{j=0}^\infty\left(b\widehat H(z\tau)\right)^j
 =\frac{e^{zT}}{1-b\widehat H(z\tau)},\qquad 0\leq t\leq T.
\end{align}
Since $|b-1|\leq C_M\delta$ for small $\eps$, equations
\eqref{eq:fractional-controlled-scalar-mean}--\eqref{eq:fractional-controlled-renewal-bound}
give $\sup_{0\leq t\leq T}R_\eps(t)\leq C_{T,M}$.
Domination of the spatial profile by $\|\lambda_0\|_\infty$ gives
the pointwise mean-density bound, while integration of the same
positive mean iteration over space replaces this bound by $L$:
\begin{align}
 \E^v[\NN_\eps^{\eta,v}(\dd t,\dd x)]
 &\leq C_{T,M}\|\lambda_0\|_\infty\dd t\dd x,
 \label{eq:fractional-controlled-mean-density-bound}\\
 \E^v[\NN_\eps^{\eta,v}(\Omega_T)]
 &\leq C_{T,M}LT.
 \label{eq:fractional-controlled-total-mean}
\end{align}
Both bounds hold for $0<\eps\leq\eps_0$, with constants independent
of $\eps$ and of the generation cutoff.

To construct the tilted law, let
$D_{\eps,T}:=(0,T/\tau)\times\R^d$ and retain the prescribed past
$\Pi^{\eps,\eta}$ with its original law. Independently of this past,
let $(\mathsf P_j)_{j\geq1}$ be independent Poisson random measures
on $D_{\eps,T}\times(0,\infty)$, each with intensity
$\dd t\dd x\dd u$. Use the filtration generated by the whole
prescribed past and by all these Poisson measures up to the current
time. Let
\[
 \vartheta_\eps(t,x):=a_\eps e^{\delta v(\tau t,\eps^\eta x)}
       \leq b,\qquad C_k:=(-k,0)\times B_k\quad(k\geq1).
\]
For parents in $C_k$, define the generation intensities and counting
measures recursively, for $j\geq1$, by
\begin{align}\label{eq:fractional-generation-recursion}
 \gamma_{k,1}(t,x)
 &:=\vartheta_\eps(t,x)\int_{-k}^0\int_{B_k}
       h(t-r,x-y)\Pi^{\eps,\eta}(\dd r,\dd y),\nonumber\\
 G_{k,j}(\dd t,\dd x)
 &:=\int_0^\infty\one_{\{u\leq\gamma_{k,j}(t,x)\}}
                      \mathsf P_j(\dd t,\dd x,\dd u),\nonumber\\
 \gamma_{k,j+1}(t,x)
 &:=\vartheta_\eps(t,x)\int_{(0,t)}\int_{\R^d}
                        h(t-r,x-y)G_{k,j}(\dd r,\dd y).
\end{align}
Here $G_{k,1}$ counts children of the prescribed parents, and
$G_{k,j+1}$ counts children of generation $j$. The restriction
$r<t$ makes each $\gamma_{k,j}$ predictable. For fixed $k,J$, the
past has finitely many atoms in $C_k$ almost surely, and each event
has at most $b$ expected children on $D_{\eps,T}$; hence the first
$J$ generations have finite expected counts. The common Poisson
measures and positivity of the kernel give
$G_{k,j}\leq G_{k+1,j}$ by induction on $j$. Define
\begin{equation}\label{eq:fractional-finite-generation-definition}
 \begin{aligned}
 N_{k,J}&:=\sum_{j=1}^J G_{k,j},\qquad N_{k,0}:=0,\\
 N_J^{\eps,\eta,v}(A)&:=\lim_{k\to\infty}N_{k,J}(A),
       \qquad A\in\mathcal B(D_{\eps,T}),\qquad N_0^{\eps,\eta,v}:=0.
 \end{aligned}
\end{equation}
These limits are increasing in $J$. Integrating the nonnegative mean
iteration of \eqref{eq:fractional-generation-recursion} over space,
using the normalization in \eqref{eq:heavyfamily}, gives
\begin{align}\label{eq:fractional-truncated-generation-mean}
 \frac{\delta^2}{n}\E^v[N_{k,J}(D_{\eps,T})]
 &\leq L\int_0^T\sum_{j=0}^{J-1}b^{j+1}
             \Pp(\tau T_j\leq t<\tau T_{j+1})\dd t\nonumber\\
 &\leq L\int_0^T R_\eps(t)\dd t\leq C_{T,M}LT.
\end{align}
The last inequality follows from
\eqref{eq:fractional-controlled-scalar-mean}--\eqref{eq:fractional-controlled-renewal-bound}.
The monotone convergence theorem in $k$ and the mass per event
$\delta^2/n$ therefore give
\begin{equation}\label{eq:fractional-generation-nonexplosion}
 \sup_{J\geq1}\E^v\left[
       N_J^{\eps,\eta,v}\left((0,T/\tau)\times\R^d\right)\right]
 \leq\frac n{\delta^2}C_{T,M}LT<\infty.
\end{equation}
Let $N^{\eps,\eta,v}(A):=\lim_{J\to\infty}N_J^{\eps,\eta,v}(A)$
for $A\in\mathcal B(D_{\eps,T})$. Applying the monotone convergence
theorem to \eqref{eq:fractional-generation-nonexplosion} yields
\begin{align}\label{eq:fractional-limit-count-bound}
 \E^v[N^{\eps,\eta,v}(D_{\eps,T})]
 =\lim_{J\to\infty}\E^v[N_J^{\eps,\eta,v}(D_{\eps,T})]
       \leq\frac n{\delta^2}C_{T,M}LT<\infty,
\end{align}
and
\begin{align*}
 \Pp^v\big(N^{\eps,\eta,v}(D_{\eps,T})<\infty\big)=1.
\end{align*}
For the finite truncation $N_{k,J}$, the compensator density is
$\sum_{j=1}^J\gamma_{k,j}$. The recursion
\eqref{eq:fractional-generation-recursion}, followed by the
monotone convergence theorem in $k$ and then in $J$, gives
\begin{align}\label{eq:fractional-generation-intensity-limit}
 \lambda_J(t,x)
 &:=\lim_{k\to\infty}\sum_{j=1}^J\gamma_{k,j}(t,x)\nonumber\\
 &=\vartheta_\eps(t,x)\bigg[
       \int_{-\infty}^0\int_{\R^d}h(t-r,x-y)
                                  \Pi^{\eps,\eta}(\dd r,\dd y)+\int_{(0,t)}\int_{\R^d}h(t-r,x-y)
                                  N_{J-1}^{\eps,\eta,v}(\dd r,\dd y)\bigg]
       \nonumber\\
 &\underset{J\to\infty}{\uparrow}
 \lambda^{\eps,\eta,v}(t,x)
 :=e^{\delta v(\tau t,\eps^\eta x)}
       \lambda^{\eps,\eta}(t,x;N^{\eps,\eta,v},\Pi^{\eps,\eta}),
\end{align}
where $\lambda^{\eps,\eta}(t,x;N,\Pi)$ is the offspring intensity
in \eqref{eq:intensity} with the kernel of Assumption~\ref{ass:heavy}.
The increasing limit in \eqref{eq:fractional-generation-intensity-limit}
is predictable. For every bounded nonnegative predictable $F$,
compensation for each $\mathsf P_j$ and the monotone convergence
theorem give
\begin{align}\label{eq:fractional-limit-compensation}
 \E^v\left[\int_{D_{\eps,T}}F(t,x)N^{\eps,\eta,v}(\dd t,\dd x)\right]
 &=\lim_{J\to\infty}\lim_{k\to\infty}
       \E^v\left[\int_{D_{\eps,T}}F(t,x)N_{k,J}(\dd t,\dd x)\right]
       \nonumber\\
 &=\lim_{J\to\infty}\lim_{k\to\infty}
       \E^v\left[\int_{D_{\eps,T}}F(t,x)
                    \sum_{j=1}^J\gamma_{k,j}(t,x)\dd t\dd x\right]
       \nonumber\\
 &=\E^v\left[\int_{D_{\eps,T}}F(t,x)
                    \lambda^{\eps,\eta,v}(t,x)\dd t\dd x\right].
\end{align}
Taking $F=1$ in \eqref{eq:fractional-limit-compensation} and using
\eqref{eq:fractional-limit-count-bound} gives
\begin{equation}\label{eq:fractional-finite-compensator}
 \E^v\left[\int_0^{T/\tau}\int_{\R^d}
                     \lambda^{\eps,\eta,v}(t,x)\dd x\dd t\right]
 =\E^v[N^{\eps,\eta,v}(D_{\eps,T})]
 \leq\frac n{\delta^2}C_{T,M}LT<\infty.
\end{equation}
Thus the constructed process has the required predictable intensity
and is nonexplosive on the observation horizon.

For a bounded compactly supported deterministic $g$, let
\begin{equation}\label{eq:fractional-controlled-equation-martingale}
 \mathcal M_\eps^v(g)
 :=\frac\delta n\int_0^{T/\tau}\int_{\R^d}g(\tau t,\eps^\eta x)
       \left[N^{\eps,\eta,v}(\dd t,\dd x)
                         -\lambda^{\eps,\eta,v}(t,x)\dd t\dd x\right].
\end{equation}
The $L^2$ isometry for integration against the compensated counting
measure, together with \eqref{eq:fractional-controlled-mean-density-bound},
gives
\begin{align}\label{eq:fractional-controlled-martingale-isometry}
 \E^v\left[\left|\mathcal M_\eps^v(g)\right|^2\right]
 &=\frac{\delta^2}{n^2}\E^v\left[
       \int_0^{T/\tau}\int_{\R^d}\left(g(\tau t,\eps^\eta x)\right)^2
                          \lambda^{\eps,\eta,v}(t,x)\dd x\dd t\right]
       \nonumber\\
 &=\frac1n\E^v\left[\NN_\eps^{\eta,v}(g^2)\right]
 \leq\frac{C_{T,M}\|\lambda_0\|_\infty}{n}
                 \int_{\Omega_T}(g(t,x))^2\dd t\dd x.
\end{align}
The exact weak identity
\eqref{eq:controlled-weak-microscopic-equation} becomes
\[
 \NN_\eps^{\eta,v}(-\mathscr L_{\eps,\eta}g-v_\eps g)
 =\delta^{-1}\cI_{\eps,\eta}(g)
       +\mathcal M_\eps^v(g)-\delta\mathcal M_\eps^v(v_\eps g),
 \qquad v_\eps=\left(1-e^{-\delta v}\right)/\delta.
\]
The mean bounds \eqref{eq:fractional-controlled-mean-density-bound}--\eqref{eq:fractional-controlled-total-mean},
 \eqref{eq:fractional-controlled-martingale-isometry},
\eqref{eq:fractional-global-generator}, and
\eqref{eq:fractional-past-superexp} imply that every subsequential
vague limit is a finite nonnegative measure satisfying
$-\nu(\mathscr L_\eta g+vg)=\beta_\eta(g)$.
The potential argument in Lemma~\ref{lem:fractional-action-domain}
and the uniqueness estimate
\eqref{eq:fractional-controlled-contraction} identify it with
$q^v\dd t\dd x$. Thus
\begin{equation}\label{eq:fractional-controlled-lln}
 \NN_\eps^{\eta,v}\longrightarrow q^v\dd t\dd x
       \quad\text{in probability, vaguely.}
\end{equation}

The stopped likelihood for this intensity change is normalized;
letting the stopping level increase preserves normalization because
the tilted process has finite expected event count on the support
of $v$. Under $\Pp^v$ it gives
\begin{align}\label{eq:fractional-likelihood-cost}
 \frac1n\log\frac{\dd\Pp}{\dd\Pp^v}
 &=-\mathcal M_\eps^v(v)
       -\NN_\eps^{\eta,v}(d_\delta(v))
       +\delta\mathcal M_\eps^v(d_\delta(v))\nonumber\\
 &\longrightarrow-\frac12\int_{\Omega_T}v^2q^v\dd t\dd x,
 \qquad d_\delta(z)=\frac{\delta z-1+e^{-\delta z}}{\delta^2}.
\end{align}
For $|z|\leq M$, Taylor's integral formula gives
\[
 d_\delta(z)=z^2\int_0^1(1-r)e^{-r\delta z}\dd r,
 \qquad |d_\delta(z)|\leq\frac12e^{\delta M}|z|^2.
\]
Consequently \eqref{eq:fractional-controlled-martingale-isometry},
applied to $g=v$ and $g=d_\delta(v)$, gives
\begin{align}\label{eq:fractional-likelihood-martingale-bounds}
 \E^v\left[\left|\mathcal M_\eps^v(v)\right|^2\right]&\leq C_v/n\longrightarrow0,\nonumber\\
 \E^v\left[\left|\delta\mathcal M_\eps^v(d_\delta(v))
 \right|^2\right]
 &\leq C_v\delta^2/n\longrightarrow0.
\end{align}
Here $C_v$ is independent of small $\eps$, since both integrands
are supported in $\supp v$ and have a common uniform bound.
Also $d_\delta(v)\to v^2/2$ uniformly on that support. The mean
bound \eqref{eq:fractional-controlled-mean-density-bound} and the
controlled limit \eqref{eq:fractional-controlled-lln} therefore give
the convergence asserted in \eqref{eq:fractional-likelihood-cost}.
Restricting \eqref{eq:fractional-likelihood-cost} to its
probability-tending-to-one event and to
$\{\NN_\eps^{\eta,v}\in O\}$ proves
\eqref{eq:fractional-controlled-lower}.
This completes the proof.
\end{proof}

\begin{proof}[Proof of Theorem~\ref{thm:fractional-spatial-ldp}]
Lemma~\ref{lem:fractional-micro-bounds} gives the closed-set upper
bound. For an open $O$ and $\nu=q\dd t\dd x\in O$ with finite
rate, choose the approximating sequence in
Lemma~\ref{lem:fractional-recovery}. For every sufficiently large
fixed $k$, its density belongs to $O$, so
\eqref{eq:fractional-controlled-lower} gives
\[
 \liminf_{\eps\downarrow0}n_{\eps,\eta}^{-1}
       \log\Pp(\NN_\eps^\eta\in O)
 \geq-\frac12\int_{\Omega_T}v_k^2q^{v_k}\dd t\dd x.
\]
Let $k\to\infty$ and then take the infimum over $\nu\in O$.
This proves the full lower bound. The domain and goodness for the rate function follow
from Lemma~\ref{lem:fractional-action-domain}; zero energy forces
$j=0$ in \eqref{eq:fractional-controlled-potential}, and hence
$q=\rho^\eta$.
\end{proof}

The full fractional action also yields the space--time Gaussian
energy on smooth perturbations of the hydrodynamic density.
\begin{corollary}[Quadratic consistency of the fractional action]
\label{cor:fractional-action-tangent}
For $z\in C_c^\infty(E_T)$, let $r_z=\mathscr H_\eta z$ with
zero past extension. Then
\begin{equation}\label{eq:fractional-action-tangent}
 \lim_{h\to0}h^{-2}\mathcal J_{\eta,\lambda_0}
                   \left((\rho^\eta+hz)\dd t\dd x\right)
 =\frac12\int_{\Omega_T}\frac{r_z^2}{\rho^\eta}\dd t\dd x
 =I_\eta(z\dd t\dd x).
\end{equation}
\end{corollary}
\begin{proof}
The residual of $\rho^\eta+hz$ is $hr_z$ by
\eqref{eq:fractional-forward-potential}. The function $r_z$ is
bounded, supported in a compact spatial set, and zero before a
positive time. On this set up to $T$, $\rho^\eta$ has a positive
lower bound by the inverse-clock representation. For small $|h|$,
\eqref{eq:fractional-full-energy} therefore gives
\[
 h^{-2}\mathcal J_{\eta,\lambda_0}\left((\rho^\eta+hz)\dd t\dd x\right)
       =\frac12\int_{\Omega_T}\frac{r_z^2}{\rho^\eta+hz}\dd t\dd x
       \longrightarrow\frac12\int_{\Omega_T}\frac{r_z^2}{\rho^\eta}\dd t\dd x.
\]
The causal identity $z=\cR_\eta r_z$ identifies the control
$\zeta=r_z/\sqrt{\rho^\eta}$ in \eqref{eq:fractionalrate}.
Uniqueness of the distributional residual identifies its energy
with $I_\eta$ and proves \eqref{eq:fractional-action-tangent}.
\end{proof}


\section{Conclusion}\label{sec:conclusion}

In this paper, we established the CLT, the MDP throughout the moderate range, and a full
space--time event-measure LDP for critical and nearly critical Hawkes
processes with a prescribed Poisson past. For finite-mean delays, positive
cluster-transform estimates yield the Gaussian covariance and moderate
deviation energy. Under finite nonzero initial mass, the LDP has a good
quadratic branching-control rate function with the prescribed initial trace;
the density-path representation of finite-action measures and approximation by smooth controls 
identify its domain.
We also derived moment-based conditions and bias corrections for replacing
the exact mean by the hydrodynamic density.

For factorized heavy-tailed power-law delays, the corresponding
CLT, MDP, and finite-mass LDP retain the source generated by the prescribed
past in the fractional evolution. Source-preserving regularization gives
approximation of every finite-action measure, and the quadratic expansion of
the action agrees with the Gaussian energy. For infinite initial mass,
compact-kernel exponential bounds hold; extending approximation by smooth controls 
to this setting remains open.

\subsection*{Acknowledgements}

Lingjiong Zhu would like to thank Prof. S.R.S. Varadhan 
for helpful discussions and suggestions.

\bibliographystyle{alpha}
\bibliography{bibtex}
\end{document}